\documentclass[final]{siamltex1213}

\usepackage[utf8]{inputenc}
\usepackage{graphicx}
\usepackage{amssymb}
\usepackage{amsmath}
\usepackage{bm}
\usepackage{url}
\usepackage{geometry}

\def\fatu{\mathbf{u}}

\def\bfatu{\bar{\mathbf{u}}}
\def\fatv{\mathbf{v}}

\def\fate{\mathbf{e}}
\def\fatf{\mathbf{f}}
\def\fatg{\mathbf{g}}
\def\fats{\mathbf{s}}
\def\fatx{\bm{x}}
\def\fatn{\bm{n}}
\def\faty{\mathbf{y}}
\def\fatw{\mathbf{w}}
\def\fatmu{\bm{\mathcal{\mu}}}

\def\fatui{\mathbf{u}_{\infty}}

\def\tfatui{\tilde{\mathbf{u}}_{\infty}}
\def\fatvi{\mathbf{v}_{\infty}}

\def\tfatvi{\tilde{\mathbf{v}}_{\infty}}
\def\fatwi{\mathbf{w}_{\infty}}

\def\tfatwi{\tilde{\mathbf{w}}_{\infty}}

\def\ordo{\mathcal{O}}

\def\bgamma{\bar{\gamma}}
\def\pOmega{\partial\Omega}

\def\fatA{\bm{A}}
\def\fatB{\bm{B}}

\def\fatD{\bm{D}}
\def\fatE{\bm{E}}
\def\fatF{\bm{F}}
\def\fatG{\bm{G}}
\def\fatH{\bm{H}}
\def\fatI{\bm{I}}
\def\fatK{\bm{K}}
\def\fatL{\bm{L}}
\def\fatS{\bm{S}}
\def\fatR{\bm{R}}
\def\fatP{\bm{P}}
\def\fatQ{\bm{Q}}
\def\fatT{\bm{T}}

\def\fatU{\bm{U}}

\def\fatM{\bm{M}}

\def\fatnu{\boldsymbol{\nu}}
\def\fatLambda{\boldsymbol{\Lambda}}
\def\calu{\boldsymbol{\xi}}
\def\calui{{\calu}_{\infty}}
\def\tcalui{\tilde{\calu}_{\infty}}
\def\calv{\boldsymbol{\eta}}
\def\calV{\mathcal{V}}

\numberwithin{equation}{section}
\numberwithin{table}{section}
\numberwithin{figure}{section}

\title{Mesoscopic modeling of stochastic reaction-diffusion kinetics
  in the subdiffusive regime}

\author{Emilie Blanc, Stefan Engblom, Andreas Hellander, Per L{\"o}tstedt\thanks{Division of Scientific Computing,
			Department of Information Technology,
			Uppsala University, P. O. Box 337, SE-75105 Uppsala, Sweden.
(\email{emilie.blanc@it.uu.se}, {stefane@it.uu.se}, {andreas.hellander@it.uu.se}, {perl@it.uu.se}). 
}}

\begin{document}

\maketitle

\begin{abstract}
  Subdiffusion has been proposed as an explanation of various kinetic
  phenomena inside living cells. In order to fascilitate large-scale
  computational studies of subdiffusive chemical processes, we extend
  a recently suggested mesoscopic model of subdiffusion into an
  accurate and consistent reaction-subdiffusion computational
  framework.
  Two different possible models of chemical reaction are revealed and
  some basic dynamic properties are derived. In certain cases those
  mesoscopic models have a direct interpretation at the macroscopic
  level as fractional partial differential equations in a bounded time
  interval. Through analysis and numerical experiments we estimate the
  macroscopic effects of reactions under subdiffusive mixing. The
  models display properties observed also in experiments: for a short
  time interval the behavior of the diffusion and the reaction is
  ordinary, in an intermediate interval the behavior is anomalous, and
  at long times the behavior is ordinary again.
  
\end{abstract}






%
%



\section{Introduction}

Quantitative models of reaction-diffusion systems are important tools
to theoretically and computationally study the dynamics of
intracellular control systems. Macromolecules such as proteins, mRNA
and DNA interact in regulatory pathways in which each chemical
reaction occurs in a certain part of the cell. Molecules are
transported via diffusion or active transport to arrive at the
subcellular location needed to fulfill their function. As an example,
in many gene regulatory pathways, proteins called transcription
factors will diffuse and search for the correct binding site on DNA,
where they modulate the translation of DNA into mRNA, which in turn
will be translated into proteins. Apart from the spatial dynamics,
realistic models need to account for large fluctuations in the copy
numbers of the species due to the small reaction volume and thus the
small number of molecules of key species such as transcription
factors. To that end, spatial stochastic models based on a Markov
process formalism are popular due to their high level of biological
realism compared to macroscopic partial differential equations (PDE),
with only a moderate increase in computational complexity and
cost. Such mesoscopic models, based on the so called
Reaction-Diffusion Master Equation (RDME), have recently been used to
address different biological phenomena such as regulation of cell
division in \emph{E. coli} \cite{FaEl}, yeast polarization
\cite{Lawson2013} and genetic oscillators \cite{Sturrock2013}.

A fundamental assumption in the RDME models is that molecules are
point particles. Thus, the model does not account for molecular volume
exclusion, something that can lead to anomalous diffusion in the
crowded compartments of living cells. Especially on 2D membranes,
crowding and anomalous diffusion behavior can be expected to have an
impact on the reaction kinetics.
For diffusion limited reactions, subdiffusion often results in a
slower decay \cite{YLL07}. A reversible ligand binding bimolecular
reaction on a 2D membrane is studied in \cite{SOULA14a} where
anomalous diffusion is simulated with a continuous-time random walk
(CTRW). The steady state distribution of the bound complex depends
critically on the anomalous diffusion and its parameters.  An
irreversible bimolecular reaction is simulated in \cite{BESB13} with
ordinary diffusion and a modified non-Fickian diffusion which is
anomalous for intermediate time but tends to ordinary diffusion in the
long term. The distribution of the reactant is sensitive to the type
of diffusion. The conclusion in \cite{Saxton02} is also that
subdiffusion changes the reaction rate in a bimolecular reaction.

To capture effects of subdiffusion due to molecular crowding in
stochastic reaction diffusion simulations, detailed Brownian Dynamics
(BD) models based on a hard-sphere assumption can be used instead of
the mesoscopic RDME model, since they account for the volume exclusion
of molecules. However, this comes at the price of a very large
increase in computational cost, making simulations on the timescales
of interest in systems biology (minutes to hours) challenging. Thus
there is a need to investigate less computationally demanding
approximations of mesoscopic subdiffusion. A promising approach for
mesoscopic anomalous diffusion simulation was suggested recently by
Mommer and Lebiedz \cite{MOMMER09}, in which a stochastic CTRW model
\cite{MONTROLL65} is approximated by an internal states model where
all transitions between states have exponentially distributed waiting
times. Hence, the subdiffusive process is simulated as a coupled
reaction-diffusion system with ordinary diffusion, making it easily
implemented in software frameworks based on the RDME \cite{urdme,
  mesoRD, steps}.

Starting from the viewpoint that the mesoscopic model arises as an
approximation of molecular crowding on the microscopic level, the
internal states model can be viewed as a mathematical means to arrive
at an approximate, faster simulation method. On a more general level,
such internal states models arise also from the need to model for
example conformational changes of macromolecules. Many proteins exist
in a number of different states due to e.g. ligand binding,
methylation, and conformational changes \cite{Bray03, EndyBrent01}. In
practical modeling, if each state is to be represented as a separate
chemical species, an exponential growth in model complexity due to a
combinatorial explosion in the number of states may be the
consequence. To handle this from a practical point of view, special
rule-based modeling languages have been developed, such as BioNetGen
Language \cite{blinov2004bionetgen} and PySB \cite{pysb}. This
combinatorial explosion in states motivates the development of
computational methods to simulate models with these types of molecules
more efficiently \cite{SBSK}. In this case, we start with an internal
states model naturally arising from the biological model, and are then
interested in how this model can be approximated on a
phenomenological, macroscopic level as a fractional PDE (FPDE).

Molecules transported by subdiffusion move slower than with ordinary
diffusion \cite{METZLER00} which could be an effect of crowding in a
biological cell. The governing equation for subdiffusion is a FPDE
with a fractional time derivative and the Laplacian as the diffusive
space operator. Chemical reactions can also be included in this
framework.  Three different models for monomolecular reactions with
subdiffusion are derived and compared with analytical solutions in
\cite{HLW06}. In one of the models where the fractional derivative
acts only on the diffusion, the concentration of the species can
become negative, making it less suitable.  In the preferred model,
where a special fractional derivative is derived for the diffusion,
the equation for the homogeneous solution is recovered without a
fractional derivative.  The model where the fractional derivative is
applied to both the diffusion and the reaction is used in \cite{YAL04}
for a bimolecular reaction. A model with internal states is derived in
\cite{Shkilev09} for the propagation of a reaction front in an
inhomogeneous medium. The particles move by diffusion and react with
each other at the front. Macroscopic equations are obtained after
summation over the internal states. Analytical results show that the
speed of the front depends on the type of reaction and the particle
distribution between the internal states. That same model is suggested
in \cite{Shkilev14} to explain the results in \cite{SOULA14a}.

In this paper, we develop theory to extend the mesoscopic anomalous
diffusion internal states model \cite{MOMMER09} to account also for
chemical reactions, making it possible to model general
reaction-diffusion processes. In particular, we develop theory for the
connection between the model framework that thus emerges and existing
models of anomalous reaction-diffusion processes on the macroscopic,
FPDE level. As we will see, the proposed reactive internal states
model is quite general and can, depending on the involved parameters,
result in different mean-field equations. Some of these mesoscopic
equations have a simple interpretation at the macroscopic level as a
FPDE but in general this is not possible. In numerical simulations,
the reaction-diffusion systems have the same behavior as observed in
many experiments: in a short time interval after start we see ordinary
behavior, then there is an anomalous phase, and finally the behavior
is ordinary again but with different diffusion and reaction
coefficients. These conclusions are supported theoretically in
\cite{HofFra,JEON11,KUSUMI05,Saxton07}. The reaction-diffusion
equations encompass many of the known models of subdiffusive reaction
systems, making it promising as a general modeling framework.

The rest of the paper is organized as follows. In Section
\ref{sec:background}, we review some basic properties of stochastic
models of diffusion and how they give rise to a phenomenological,
macroscopic FPDE in the thermodynamic limit. Here, we will also see
how the addition of chemical reactions on the macroscopic level can
lead to two different, possible FPDE models.  We introduce an internal
states model of reaction-subdiffusion systems in
Section~\ref{sec:internalstates} and analyze its mean field
properties.  In Section 4, we relate the models in
Sections~\ref{sec:background} and \ref{sec:internalstates} to each
other for linear and non-linear chemical reactions, and in Setion 5,
numerical experiments are presented. Finally, the paper is concluded
in Section 6.


\section{Fractional partial differential equations as limits of continuous-time random walks}
\label{sec:background}

In this section, we recall some well-known facts about Brownian and, in
particular, subdiffusive random movements. When approached in the
proper macroscopic limit, FPDEs emerge as a convenient mathematical model for subdiffusion. Terms are added
to the FPDEs to model chemical reactions. 






The CTRW model was introduced by Montroll and Weiss for hopping transport on a disordered lattice \cite{MONTROLL65}. In this model in continuous space, the particle is assumed to traverse the space by a series of jumps. The displacement and the waiting time to perform the next jump are drawn from a given probability density function (PDF) $\Psi(\fatx,t)$. We assume that the jump length PDF $\lambda(\fatx)$ and waiting time PDF $\psi(t)$ are independent random variables. Consequently, $\Psi(\fatx,t)$ is written
\begin{equation}
\Psi(\fatx,t) = \psi(t)\,\lambda(\fatx).
\label{eq:pdf_ctrw}
\end{equation}
Different diffusion processes can be categorized by the expected waiting time
\begin{equation}
\tau^* = \int_0^{\infty}t\,\psi(t)\,dt,
\label{eq:time_charac_ctrw}
\end{equation}
and the jump length variance
\begin{equation}
\Sigma^2 = \int_{\mathbb{R}^d}\|\fatx\|^2_2\,\lambda(\fatx)\,d\fatx,
\label{eq:length_charac_ctrw}
\end{equation}
where $d$ is the dimension of the embedding space. If both $\Sigma^2$ and $\tau^*$ are finite, the long-time limit corresponds to Brownian motion. A diverging $\tau^*$ with finite $\Sigma^2$ gives rise to subdiffusion. On the contrary, a diverging $\Sigma^2$ with a finite $\tau^*$ induces superdiffusion \cite{HofFra,METZLER00}, which is beyond the scope of this article.

\subsection{Brownian motion}

We consider a Gaussian jump length PDF
\begin{equation}
\lambda(\fatx) = \frac{1}{\left(4\,\pi\,\sigma^2\right)^{d/2}}\,e^{-\|\fatx\|^2_2/(4\,\sigma^2)}
\label{eq:jump_length_gaussian}
\end{equation}
and a Poissonian waiting time PDF
\begin{equation}
\psi(t) = \frac{1}{\tau}\,e^{-t/\tau}.
\label{eq:waiting_time_brownian}
\end{equation}
Equations \eqref{eq:time_charac_ctrw}--\eqref{eq:length_charac_ctrw} lead to $\Sigma^2 = 2\,\sigma^2 < \infty$ and $\tau^* = \tau <  \infty$. Since $\Sigma^2$ and $\tau^*$ are finite, the long-time limit thus corresponds to Brownian motion. At the macroscopic scale, we recover the classical diffusion equation \cite{HofFra,METZLER00} for the concentration $U$ of the chemical species $A$
\begin{equation}
\frac{\partial U}{\partial t} = D\,\Delta U,
\label{eq:PDE_brownian}
\end{equation}
with 
\begin{equation}
D = \frac{\sigma^2}{\tau} = \frac{\Sigma^2}{2\,\tau}.
\label{eq:diffusion_coef_brownian}
\end{equation}
The linear time dependence of the mean squared displacement
\begin{equation}
\left\langle \|\fatx\|^2_2(t) \right\rangle  = 2\,d\,D\,t
\label{eq:MSD_brownian}
\end{equation}
is characteristic of Brownian motion.

\subsection{Subdiffusion}

We consider now instead a Gaussian jump length PDF \eqref{eq:jump_length_gaussian} and a waiting time PDF \cite{KLAFTER12,MARQUEZLAGO12}
\begin{equation}
\psi(t) = \frac{t^{\alpha-1}}{\tau^{\alpha}}\,E_{\alpha,\alpha}\left(-\left(\frac{t}{\tau}\right)^{\alpha}\right),\quad 0<\alpha<1,
\label{eq:waiting_time_subdiffusion}
\end{equation}
where
\begin{equation}
E_{\alpha,\alpha}\left(-\left(\frac{t}{\tau}\right)^{\alpha}\right) = \sum\limits_{k=0}^{\infty} \frac{\left(-(t/\tau)^{\alpha}\right)^k}{\Gamma(\alpha\,k + \alpha)}
\label{eq:generalized_mittag_leffler}
\end{equation}
is the generalized Mittag-Leffler function. The Mittag-Leffler waiting time PDF has been observed experimentally, from polymer rheology \cite{GLOCKLE91}, over ligand rebinding to proteins \cite{GLOCKLE95} and protein conformation dynamics \cite{YANG03}, to financial market time series \cite{MAINARDI00}. Equations \eqref{eq:time_charac_ctrw}--\eqref{eq:length_charac_ctrw} lead to $\Sigma^2 = 2\,\sigma^2 < \infty$ and $\tau^* = \infty$, which give rise to subdiffusion. At the macroscopic scale, we obtain the FPDE \cite{HofFra,METZLER00}
\begin{equation}
\frac{\partial^{\alpha} U}{\partial t^{\alpha}} = K_{\alpha}\,\Delta U,
\label{eq:FPDE}
\end{equation}
with
\begin{equation}
K_{\alpha} = \frac{\sigma^2}{\tau^{\alpha}} = \frac{\Sigma^2}{2\,\tau^{\alpha}}.
\label{eq:Kalpha}
\end{equation}
At the boundary $\pOmega$ of the domain $\Omega$, the molecules are reflected back implying homogeneous Neumann
boundary conditions at $\pOmega$. The operator $\frac{\partial^{\alpha}}{\partial t^{\alpha}}$ involved in \eqref{eq:FPDE} is a Caputo fractional derivative in time of order $\alpha$, generalizing the usual derivative. It is defined as \cite{CAPUTO67,MILLER93}
\begin{equation}
\frac{\partial^{\alpha} U}{\partial t^{\alpha}} = \int_0^t \frac{(t-\tau)^{-\alpha}}{\Gamma(1-\alpha)}\,\frac{dU}{dt}(\tau)\,d\tau.
\label{eq:caputo_fractional_derivative}
\end{equation}
The mean squared displacement is given by the power law
\begin{equation}
\left\langle \|\fatx\|^2_2(t) \right\rangle = \frac{2\,d\,K_{\alpha}}{\Gamma(1+\alpha)}\,t^{\alpha}.
\end{equation}

In the one dimensional (1D) case ($d=1$) in free space, the Green's function of the FPDE \eqref{eq:FPDE} is \cite{METZLER00}
\begin{equation}
U(x,t) = \frac{1}{\sqrt{4\,\pi\,K_{\alpha}\,t^{\alpha}}}\,H^{2,0}_{1,2}\left[ \frac{x^2}{4\,K_{\alpha}\,t^{\alpha}}\left|\begin{array}{ll}
\displaystyle \left( 1-\frac{\alpha}{2},\alpha\right)  & \\
[10pt]
(0,1) & \left( \frac{1}{2},1\right) 
\end{array}\right.
\right],
\label{eq:sol_FPDE1_fox}
\end{equation}
where the Fox function $H^{2,0}_{1,2}$ is defined in Appendix~\ref{annexe:special_functions}. In the particular case $\alpha = 1/2$, the Green's function \eqref{eq:sol_FPDE1_fox} can be rewritten
\begin{equation}
U(x,t) = \frac{1}{\sqrt{8\,\pi^3\,K_{1/2}\,t^{1/2}}}\,G^{3,0}_{0,3}\left[ \left(\frac{x^2}{16\,K_{1/2}\,t^{1/2}}\right)^2 \left| \; \left( 0,\frac{1}{4},\frac{1}{2}\right) \right.
\right],
\label{eq:sol_FPDE1_meijer}
\end{equation}
where the Meijer-G function $G^{3,0}_{0,3}$ is also defined in Appendix~\ref{annexe:special_functions}. The Green's function \eqref{eq:sol_FPDE1_meijer} will be used in numerical experiments in Section~\ref{sec:num_exp} as a reference solution.

\subsection{Adding chemical reactions}
\label{subsec:reaction_modeling}
In the following 
we discuss how chemical reactions can be added to the macroscopic subdiffusion model \eqref{eq:FPDE} for the cases of annihilation, reversible isomerization, and reversible bimolecular reactions. As we will see, in the first two cases two different FPDE models are possible, while only one of them turns out to be well defined for bimolecular association.  
\subsubsection{Annihilation}\label{subsubsec:reaction_modeling_annihilation}

We first consider one species $A$ and the annihilation process $A \mathop{\longrightarrow}\limits^{k_\ast} \emptyset$. Two different FPDEs can be used to model an annihilation process. 
If a constant proportion of walkers are removed instantaneously at the start of each step then the long-time asymptotic limit yields a fractional reaction-diffusion equation with a fractional order temporal derivative operating both on the standard diffusion term and on the linear reaction kinetics term \cite{HLW06,HORNUNG05,SEKI03}:
\begin{equation}
\frac{\partial U}{\partial t} = \frac{\partial^{1-\alpha}}{\partial t^{1-\alpha}}\left(K_{\alpha}\,\Delta U - k_\ast\,U\right).
\label{eq:FPDE_decayI}
\end{equation}
In what follows, \eqref{eq:FPDE_decayI} is referred to as \emph{model I}. The total amount $\bar{U}$ of $A$ in a bounded domain $\Omega$ with boundary $\partial \Omega$ is defined as
\begin{equation}
\bar{U} = \int_{\Omega} U(\fatx,t)\,d\Omega.
\label{eq:Na}
\end{equation}
The equation \eqref{eq:FPDE_decayI} is integrated over $\Omega$
\begin{equation}
\frac{d \bar{U}}{d t} = \frac{\partial^{1-\alpha}}{\partial t^{1-\alpha}}\left(K_{\alpha}\,\int_{\Omega}\Delta U\,d\Omega - k_{\ast}\,\bar{U}\right).
\label{eq:PDE_Na}
\end{equation}
Since
\begin{equation}
\int_{\Omega}\Delta U\,d\Omega = \int_{\pOmega}\fatn\cdot\nabla U\,dS = 0
\label{eq:neumannBC}
\end{equation}
for Neumann conditions at $\pOmega$ with normal $\fatn$, $\bar{U}$ satisfies the fractional ordinary differential equation (ODE)
\begin{equation}
\frac{d^{\alpha}\bar{U}}{dt^{\alpha}} = -k_{\ast}\,\bar{U}.
\label{eq:PDE_Na_modelI}
\end{equation}
Consequently, the total amount of $A$ is also affected by the subdiffusion. This behavior is called {\it anomalous kinetics}. In the 1D case in free space, the Green's function of \eqref{eq:FPDE_decayI} is \cite{HLW06}
\begin{equation}
U(x,t) = \frac{1}{\sqrt{4\,\pi\,K_{\alpha}\,t^{\alpha}}}\,\sum\limits_{j=0}^{\infty}\frac{(-k_\ast\,t^{\alpha})^j}{j!}\,H^{2,0}_{1,2}\left[ \frac{x^2}{4\,K_{\alpha}\,t^{\alpha}}\left|\begin{array}{ll}
\displaystyle \left( 1-\frac{\alpha}{2} + \alpha\,j,\alpha\right)  & \\
[10pt]
(0,1) & \left( \frac{1}{2} + j,1\right) 
\end{array}\right.
\right].
\label{eq:sol_FPDE_decayI_fox}
\end{equation}
In the particular case of $\alpha = 1/2$, \eqref{eq:sol_FPDE_decayI_fox} can be rewritten as
\begin{equation}
U(x,t) = \frac{1}{\sqrt{8\,\pi^3\,K_{1/2}\,t^{1/2}}}\,\sum\limits_{j=0}^{\infty}\frac{(-2\,k_\ast\,t^{1/2})^j}{j!}\,G^{3,0}_{0,3}\left[ \left(\frac{x^2}{16\,K_{1/2}\,t^{1/2}}\right)^2 \left| \; \left( 0,\frac{1}{4} + \frac{j}{2}, \frac{1}{2}\right)\right.
\right].
\label{eq:sol_FPDE_decayI_meijer}
\end{equation}

If instead the walkers are removed at a constant per capita rate during the waiting time between steps then the long time asymptotic limit has a standard linear reaction kinetics term but a fractional order temporal derivative operating on a nonstandard diffusion term \cite{HLW06,LOMHOLT07,SOKOLOV06,YADAV06}:
\begin{equation}
\frac{\partial U}{\partial t} = K_{\alpha}\,e^{-k_\ast\,t}\,\frac{\partial^{1-\alpha}}{\partial t^{1-\alpha}}\left(e^{k_\ast\,t}\,\Delta U\right) - k_\ast\,U.
\label{eq:FPDE_decayII}
\end{equation}
This case is referred to as \emph{model II} in what follows. Integrating \eqref{eq:FPDE_decayII} and using \eqref{eq:neumannBC} leads to
\begin{equation}
\begin{array}{ll}
\displaystyle \frac{\partial \bar{U}}{\partial t} & \displaystyle = K_{\alpha}\,e^{-k_\ast\,t}\,\frac{\partial^{1-\alpha}}{\partial t^{1-\alpha}}\left(e^{k_\ast\,t}\,\int_{\Omega}\Delta U\,d\Omega\right) - k_\ast\,\bar{U}= - k_\ast\,\bar{U}.
\end{array}
\label{eq:PDE_Na_modelII}
\end{equation}
Ordinary kinetics is recovered: the total amount of $A$ is not affected by the subdiffusion. The change of variables $\tilde{U}=e^{k_\ast\,t}U$ in \eqref{eq:FPDE_decayII} leads to the FPDE \eqref{eq:FPDE}. In the one dimensional case, the Green's solution of \eqref{eq:FPDE_decayII} is then \cite{HLW06, METZLER00}
\begin{equation}
U(x,t) = \frac{e^{-k_\ast\,t}}{\sqrt{4\,\pi\,K_{\alpha}\,t^{\alpha}}}\,H^{2,0}_{1,2}\left[ \frac{x^2}{4\,K_{\alpha}\,t^{\alpha}}\left|\begin{array}{ll}
\displaystyle \left( 1-\frac{\alpha}{2},\alpha\right)  & \\
[10pt]
(0,1) & \left( \frac{1}{2},1\right) 
\end{array}\right.
\right].
\label{eq:sol_FPDE_decayII_fox}
\end{equation}

\subsubsection{Monomolecular reactions}
\label{sec:analytical_mono}

Let two species $A$ and $B$ undergo the monomolecular reversible
isomerization reaction $A
\mathop{\rightleftharpoons}\limits^{k_\ast}_{\ell_\ast} B$. We assume
that the jump length variance and the time scale involved in
\eqref{eq:Kalpha} are the same for each species:
$\Sigma_A^2=\Sigma_B^2$ and $\tau_A=\tau_B=\tau$. Consequently, the
diffusion coefficients are also identical $K_{A\alpha} = K_{B\alpha} =
K_{\alpha}$. The concentrations of $A$ and $B$ are $U$ and $V$,
respectively. As in the case of annihilation, two different FPDEs can
be used to model this reaction. Model I is here written
\begin{equation}
\frac{\partial \fatU}{\partial t} = K_{\alpha}\,\frac{\partial^{1-\alpha}}{\partial t^{1-\alpha}}\left(\Delta \fatU - \fatR\,\fatU\right),
\label{eq:FPDE_monoI}
\end{equation}
with $\fatU = (U\;V)^T$ and $\fatR = \left( \begin{array}{rr}k_\ast & -\ell_\ast\\-k_\ast & \ell_\ast\end{array}\right)$. The reaction matrix $\fatR$ is diagonalizable, $\fatR=\fatP\,\fatLambda_r\,\fatP^{-1}$, with
\begin{equation}
\fatP = \left( \begin{array}{rc}
1 & \ell_\ast \\ -1 & k_\ast
\end{array} \right),\quad \fatLambda_r = \left( \begin{array}{cc}
k_\ast + \ell_\ast & 0 \\ 0 & 0
\end{array} \right).
\label{eq:reaction_matrix}
\end{equation}
By setting $\tilde{\fatU} = \fatP^{-1}\,\fatU$, the governing system of evolution equations decouples to
\begin{equation}
\frac{\partial \tilde{\fatU}}{\partial t} = K_{\alpha}\,\frac{\partial^{1-\alpha}}{\partial t^{1-\alpha}}\left(\Delta \tilde{\fatU} - \fatLambda_r\,\tilde{\fatU}\right).
\label{eq:FPDE_monoI_w}
\end{equation}
In the 1D case without boundaries, the Green's function of the FPDE \eqref{eq:FPDE_monoI_w} can be computed using \eqref{eq:sol_FPDE1_fox} and \eqref{eq:sol_FPDE_decayI_fox}.

Model II is written
\begin{equation}
\frac{\partial \fatU}{\partial t} = K_{\alpha}\,e^{-\fatR\,t}\,\frac{\partial^{1-\alpha}}{\partial t^{1-\alpha}}\left(\,e^{\fatR\,t}\,\Delta \fatU\right) - \fatR\,\fatU,
\label{eq:FPDE_monoII}
\end{equation}
The change of variables $\tilde{\fatU}=e^{\fatR\,t}\fatU$ in \eqref{eq:FPDE_monoII} leads to the FPDE \eqref{eq:FPDE}.
The Green's function of the FPDE \eqref{eq:FPDE_monoII} in 1D can be computed using \eqref{eq:sol_FPDE_decayII_fox}.

\subsubsection{Bimolecular reactions}

Finally, we consider three species $A$, $B$ and $C$ and the reversible bimolecular reaction $A + B \mathop{\rightleftharpoons}\limits^{k_\ast}_{\ell_\ast} C$. As previously, we assume that the jump length variance and the time scale involved in \eqref{eq:Kalpha} are the same for all species, which implies that the diffusion coefficients are also identical. The concentrations of $A$, $B$ and $C$ are $U$, $V$ and $W$, respectively.
In the bimolecular case, due to the nonlinearity, only model I is defined. It is written
\begin{equation}
\left\lbrace 
\begin{array}{lll}
\displaystyle \frac{\partial U}{\partial t} & = & \displaystyle K_{\alpha}\,\frac{\partial^{1-\alpha}}{\partial t^{1-\alpha}}\left(\Delta U - k_\ast\,U\,V + \ell_\ast\,W\right),\\
[10pt]
\displaystyle \frac{\partial V}{\partial t} & = & \displaystyle K_{\alpha}\,\frac{\partial^{1-\alpha}}{\partial t^{1-\alpha}}\left(\Delta V - k_\ast\,U\,V + \ell_\ast\,W\right),\\
[10pt]
\displaystyle \frac{\partial W}{\partial t} & = & \displaystyle K_{\alpha}\,\frac{\partial^{1-\alpha}}{\partial t^{1-\alpha}}\left(\Delta W + k_\ast\,U\,V - \ell_\ast\,W\right).
\end{array}
\right. 
\label{eq:FPDE_bimoI}
\end{equation}
No analytical solution is available in this case. Nevertheless, the quantity $U-V$ satisfies the FPDE \eqref{eq:FPDE}. Consequently, the analytical solution of $U-V$ is given by \eqref{eq:sol_FPDE1_fox}.

\section{The internal states model}\label{sec:internalstates}

In this section, we introduce an internal states model of subdiffusion on the mesoscale. The molecules jump on a lattice and can react 
with each other when they are in the same lattice cell. Diffusion is here modeled by discrete jumps with a waiting time PDF $\psi(t)$ 
and a jump length depending on the nearest neighbors in the lattice.

\subsection{Mesoscopic diffusion and reactions}
In the stochastic mesoscopic model, the state of the system at time $t$ is given by the number of molecules of each one of the $N$ species in the
state vector $\faty(t)\in\mathbb{Z}_+^N$. The domain is partitioned into non-overlapping voxels or compartments $\calV_i, \; i=1,\ldots,M,$ with a mesh or a lattice.  
Each voxel has a vertex at $\fatx_i$
in the center. The component $y_{jk}$ of $\faty_k$ is the number of molecules of species $j$ in $\calV_k$. The states are changed 
randomly at discrete time points $t_\ell$
depending on chemical reactions or a jump to a neighboring voxel due to diffusion. Between $t_\ell$ and $t_{\ell+1}$, $\faty$ is constant. 
The PDF is $p(\faty, t)$ for the system to be in state $\faty$ at time $t$, and satisfies a master
equation \cite{EnFeHeLo, VanKampen}, commonly referred to as the Reaction-Diffusion Master Equation (RDME) in the context of chemical kinetics.

The jump coefficient $\lambda_{ji}$ for a diffusive jump of a molecule $A$ 
from $\calV_j$ to $\calV_i$ is determined from the coefficients of a discretized diffusion equation 
in \cite{EnFeHeLo}. The rate for a jump is $\lambda_{ji}y_j$ where $y_j$ is the copy number of $A$ in $\calV_j$.
After the diffusive jump, $y_j$ and $y_i$ are updated: $y_j:=y_j-1,\; y_i:=y_i+1$.

A bimolecular reaction between the molecular species $A$ and $B$
produces a $C$ molecule
in a voxel in
\begin{equation}
A +B\mathop{\rightarrow}^{\kappa} C.
\label{eq:bimol_IS}
\end{equation}
The propensity of the reaction is $\kappa$ times the copy numbers $y$ and $z$ of $A$ and $B$. If the reaction takes place at $t_\ell$, then the state
before the reaction $\faty(t^-_\ell)$ in the voxel is changed to $\faty(t^+_\ell)=\faty(t^-_\ell)+\fatnu_r$ immediately after the reaction $r$, where $\fatnu_r$ is the 
stoichiometric vector associated with the reaction. In \eqref{eq:bimol_IS}, $y:=y-1,\; z:=z-1,$ and the copy number of $C$ increases by one.

A realization of the chemical system is generated by the stochastic simulation algorithm (SSA) \cite{gillespie} in a Monte Carlo method. In the limit of large 
copy numbers, the mean values of the concentrations of the species converge to the solution of the reaction rate equations \cite{Kur70}. These equations
are stated in Section~\ref{subsec:mean_properties}.

In the next section, a master equation is derived for the PDF of a system with internal states. Internal states of the species are introduced in \cite{MOMMER09} to model subdiffusion. We will extend their model to include reactions between the species.

\subsection{A generalized master equation}
\label{subsec:genmastereq}

As a generalization of the classical chemical master equation (CME), we let a chemical species $A$ have a number $N$ of internal states $A_l,\, l=1,\ldots,N$. The state $A_l$ can model a conformational state that can be difficult to observe and that is more or less hidden from a practical point of view as in a hidden Markov model. On the molecular level, the hidden states can model different things, for example, they can represent different geometrical configurations of the molecule or the macromolecule can have different small molecules attached to it at different locations \cite{Bray03, EndyBrent01}. 

We now consider a general event $q$ which could be a diffusive jump, a change of internal states or a reaction.
Assume that the molecule in the state $j$ has a PDF for the waiting time
with density $\psi_{qj}(t)$ for the event $q$ and that the molecule, due to the event, then changes state to
$i$ with probability $\pi_{qij}$ such that $\sum_i \pi_{qij}=1$. Although formally  $\psi_{qj}(t)$ here is slightly different from the one 
in Section \ref{sec:background}, the notation is retained since the observed effects are the same.  
Usually, the waiting time is assumed to be exponentially distributed and dependent on reaction rates and diffusion propensities.

The Laplace transform of a function $f(t)$ is denoted by
$\tilde{f}(s)$. An auxiliary function $\phi_{qj}(t)$ is defined by its
Laplace transform involving the Laplace transform of the waiting time
PDF $\tilde{\psi}_{qj}(s)$
\begin{equation}
  \tilde{\phi}_{qj}(s)=s\,\tilde{\psi}_{qj}(s)/(1-\tilde{\psi}_{qj}(s))
\label{eq:Laplace}
\end{equation}
Gillespie shows in \cite{Gill77} (see also \cite{HausKehr, MONTROLL73}) that there is a generalized master equation 
for the PDF $p_i$ of state $i$  
\begin{equation}
\displaystyle{\frac{\partial p_i(t)}{\partial t} = \sum_q\int_0^t \sum_j \left(\pi_{qij}\,\phi_{qj}(t-t')\,p_j(t')
                                              -\pi_{qji}\,\phi_{qi}(t-t')\,p_i(t')\right)\, dt'}
\label{eq:ME}
\end{equation}
If $\psi_{qi}(t)=a_{qi}\exp(-a_{qi} t)$, then $\phi_{qi}(t)=a_{qi}\,\delta(t)$ where $\delta(t)$ is the Dirac delta and $a_{qi}$ 
is the reaction propensity. Then we obtain the usual master equation
\begin{equation}
\begin{array}{rl}
\displaystyle{\frac{\partial p_i(t)}{\partial t} = \sum_q\sum_j (\pi_{qij}\,a_{qj}\,p_j(t)-\pi_{qji}\,a_{qi}\,p_i(t))= \sum_q\sum_j (\pi_{qij}\,a_{qj}\,p_j(t))-a_{qi}\,p_i(t)}.
\end{array}
\label{eq:ME2}
\end{equation}


Let state $i$ be $\faty$ and let state $j$ be $\faty-\fatnu_q$ in \eqref{eq:ME2} and identify $a_q(\faty-\fatnu_q)$ with $\pi_{qij}\,a_{qj}$
and $a_{qi}$ with $a_q(\faty)$ for event $q$.
Let $Q$ be the
total number of events. If the waiting times are exponentially distributed then we recover the usual CME \cite{VanKampen}
\begin{equation}
\displaystyle{\frac{\partial p(\faty, t)}{\partial t} = \sum_{q=1}^Q \left(a_q(\faty-\fatnu_q)\,p\,(\faty-\fatnu_q,t)
                                              -a_q(\faty)\,p(\faty,t)\right)}.
\label{eq:GME2}
\end{equation}

\begin{table}[h!]
\centering
\begin{tabular}{|l|l|c|c|}
\hline
& event & probability& propensity \\
\hline
Change of internal state & $A_j\rightarrow A_i$ & $\pi_{ij}=\mu_i$ & $a_j=y_j/\tau_j$ \\
Diffusion from voxel $k$ to voxel $\ell$ & $A_{jk}\rightarrow A_{j\ell}$ & $\pi_{\ell k}=\lambda_{\ell k}/\lambda_k$ &  
$a_j=\sigma_A^2\,\lambda_k\,y_{jk}/\tau_j$  \\
Monomolecular reaction & $A_j\rightarrow B_i$ & $\pi_{ij}=\kappa_{ji}/\kappa_j$ & $a_j=\kappa_j\,y_{j}/\tau_j$  \\
Bimolecular reaction & $A_i+B_j\rightarrow C_k$ & $\pi_{ij}=\kappa_{ijk}/\kappa_{ij}$ & 
$a_{ij}=\kappa_{ij}\,y_{i}\,z_{j}/\tau_{ij}$ \\
Production of one molecule & $\emptyset\rightarrow A_i$ & $\pi_{i0}=\mu_{i}$ & $a_{i}=k/\tau_{i}$  \\
Annihilation & $A_i\rightarrow \emptyset$ & $\pi_{0i}=1$ & $a_{i}=k\,y_j/\tau_{i}$   \\
\hline
\end{tabular}
\caption{Jump probabilities and propensities in \eqref{eq:ME2} for different events.}
\label{tab:coefficients}
\end{table}
We consider six different elementary events listed in Table~\ref{tab:coefficients}.
The coefficients in \eqref{eq:ME2} and waiting times in \eqref{eq:Laplace} defining the events are also found in Table~\ref{tab:coefficients}.
The copy numbers of the species $A_j$ and $B_i$ are $y_j$ and $z_i$, respectively. 
The total rates away from the state in the table to all other
states that can be reached are 
\[
    \sum_i \mu_i=1,\; \lambda_k=\sum_\ell \lambda_{\ell k},\; \kappa_j=\sum_i \kappa_{ji},\; \kappa_{ij}=\sum_k \kappa_{ijk}
\]
ensuring that $\sum_i\pi_{ij}=1$. The waiting time distribution in the state is $\psi_j(t)=a_j\,\exp(-a_j\,t)$
when one molecule is transformed in Table~\ref{tab:coefficients}.
The time scale $\tau_j$ of the transformations depends on the internal state $j$ of $A$. 
We will find that if $\tau_i\ne \tau_j$ when $i\ne j$ 
then the annihilation, diffusion, reaction, or production is anomalous at the macroscopic level. If $\tau_i= \tau$ then the transformation
between the states is the ordinary one with the same waiting time for all states. 

The scale $\tau_{ij}$ of the bimolecular reaction
depends on the states of $A$ and $B$. By splitting the reaction rate into $\theta\kappa_{ijk}/\tau_i+(1-\theta)\kappa_{ijk}/\tau_j$ with
$0\leqslant\theta\leqslant 1$, the time constant is $\tau_{ij}=1/(\theta\tau_i^{-1}+(1-\theta)\tau_j^{-1})$. 
This $\tau_{ij}$ is in agreement with the time scale obtained from \cite{CoKi}.
Let $D_{Ai}$ and $D_{Bi}$ be the diffusion coefficients of species $A$ and $B$ in state $i$, $k_b$ the microscopic reaction rate, and $\rho_r$ the reaction radius.
Then the reaction rate 
$\kappa_{ijk}$ in \cite{CoKi} follows from Smoluchowski's rate law and is a function of $D_{Ai}$ and $D_{Bj}$ 
\begin{equation}
\displaystyle \kappa_{ijk}   = \frac{k_b\,4\,\pi\,(D_{Ai} + D_{Bj})\,\rho_r}{k_b + 4\,\pi\,(D_{Ai} + D_{Bj})\,\rho_r}.
\label{eq:CollinsKimball}
\end{equation}
Let the diffusion coefficients of species $A$ and $B$ depend on the internal state such that $D_{Ai}=\sigma_A^2/\tau_i$ and $D_{Bj}=\sigma_B^2/\tau_j$ \eqref{eq:diffusion_coef_brownian}.
Then with $\theta=\sigma_A^2/(\sigma_A^2+\sigma_B^2)$ and for diffusion limited systems with a large $k_b$ compared to $D_{Ai}+D_{Bj}$, $\kappa_{ijk}$ is approximated by
\begin{equation}
\displaystyle \kappa_{ijk}   \approx 4\,\pi\,(\sigma_A^2+\sigma_B^2)\left(\frac{\theta}{\tau_i}+\frac{1-\theta}{\tau_j}\right).
\label{eq:CollinsKimball2}
\end{equation}
When $k_b$ is small then the influence of the diffusion disappears in \eqref{eq:CollinsKimball} 
and $\kappa_{ijk}\approx k_b$ without dependence on the states $i$ and $j$.

With the master equation \eqref{eq:ME2} or \eqref{eq:GME2}, we can derive reaction rate equations approximately satisfied by the 
mean values as in \cite{VanKampen, MOMMER09}.

\subsection{Mean-field properties of the internal states model}
\label{subsec:mean_properties}

The mean values of the concentrations of the species approximately satisfy a system of ODEs often denoted the
reaction rate equations. These equations are derived from the PDF in the
master equation \eqref{eq:ME2} or \eqref{eq:GME2}, see \cite{VanKampen}. When the reactions are such that
all propensities are linear in the chemical system, see Table~\ref{tab:coefficients}, then the solutions to the equations are
the exact mean values. They are approximations if there is a bimolecular reaction in the chemical system.
  
\subsubsection{Diffusion}

Let us first examine a system with one molecular species and $N$ internal states. This is the problem investigated in \cite{MOMMER09}.
We allow changes of internal state and diffusion between two voxels but without chemical reactions. The
mean values of the concentrations of the $N$ states $\fatu_i(t)\in\mathbb{R}^N$ in voxel $i$ are the solution of
\begin{equation}
\displaystyle \frac{\partial \fatu_i}{\partial t}=\sigma^2\left(\sum_{j=1}^{n_i} \frac{\lambda_{ij}}{\tau_i}\,\fatu_j-\frac{\lambda_{i}}{\tau_i}\,\fatu_i\right)+\fatA\,\fatu_i
                                                 =\frac{\sigma^2}{\tau_i}\left(\sum_{j=1}^{n_i} \lambda_{ij}\,\fatu_j-\lambda_{i}\,\fatu_i\right)+\fatA\,\fatu_i.
\label{eq:Mommer_diff}
\end{equation}
The number of voxels directly connected to $\calV_i$ is $n_i$. Consequently, a diffusive jump between $\calV_j$ and $\calV_i$ is possible.
The elements of $\fatA$ follow from Table~\ref{tab:coefficients} and are $A_{ij}=\mu_i/\tau_j,\; i\ne j,$ and $A_{ii}=(\mu_i-1)/\tau_i$. 
Let $\fatT$ be a diagonal matrix with $T_{ii}=1/\tau_i$, $\fatmu$ a vector with non-negative components $\mu_i$ such that $\fate^T\fatmu=1$ and $\fate$ a vector with $e_i = 1$ for all $i$. Then $\fatA$ in \eqref{eq:Mommer_diff} can be written
\begin{equation}
   \fatA=(\fatmu\,\fate^T-\fatI)\,\fatT.
\label{eq:Aexpr}
\end{equation}
The nullspace consists of one vector $ \tfatui$ where 
\begin{equation}
    \tfatui=\fatT^{-1}\,\fatmu.
\label{eq:uinfexpr}
\end{equation}
The left eigenvector of $\fatA$ corresponding to eigenvalue $\lambda_1=0$ is $\fate$ such that
\begin{equation}
\fate^T\,\fatA=0.
\label{eq:kernelA}
\end{equation}
The diffusion jump coefficients in \eqref{eq:Mommer_diff} are derived such that the Laplacian is approximated in voxel $i$
\begin{equation}
   \Delta \fatu(\fatx, t)\approx\sum_{j=1}^{n_i} \lambda_{ij}\,\fatu_j-\lambda_{i}\,\fatu_i.
\label{eq:Laplaceapprox}
\end{equation}
On a Cartesian mesh with equal mesh spacing $h$, $\lambda_{ji}=1/h^2$ and $\lambda_i=2d/h^2$ where $d$ is the dimension. On an unstructured mesh,
the coefficients can be derived by a finite element method as in \cite{EnFeHeLo}. With a continuous $\fatu(\fatx, t)$ in space, the equation
approximated by \eqref{eq:Mommer_diff} is 
\begin{equation}
\displaystyle \frac{\partial \fatu}{\partial t}=\fatD\,\Delta\fatu+\fatA\,\fatu,\quad \fatD=\sigma^2\,\fatT.
\label{eq:analyt_diff}
\end{equation}
The boundary conditions are of Neumann type at the boundary to preserve the total concentration of the species as in \eqref{eq:neumannBC}.
The analysis is simplified in this section if we consider the solution $\fatu(\fatx, t)$ of \eqref{eq:analyt_diff} instead of the
discrete solutions $\fatu_i(t)$ in the voxels in \eqref{eq:Mommer_diff}.

With positive initial data, it is easy
to see that there is a unique positive steady-state solution $\tfatui$ as $t \to
\infty$. This solution is space independent. We write this in the normalized form
\begin{align}\label{eq:liminfu}
  \lim_{t\rightarrow \infty}\fatu(\fatx, t) &= \tfatui=u\fatui,
\end{align}
where 
$\|\fatui\|_{1} = 1$ and $u = \|\fatu(x,0)\|_{1}$ by the
preservation of mass. Since $\Delta \fatui=0$ we have $\fatA\,\fatui=0$ and $\fatui$ is given by \eqref{eq:uinfexpr} and \eqref{eq:liminfu}. 

Expand the solution of \eqref{eq:analyt_diff} in a cosine series in 1D
\begin{equation}
   \fatu(x, t)=\sum_{\omega=0}^\infty \fatu_{\omega}(t)\,\cos(\omega x),
\label{eq:cosineu}
\end{equation}
with $x$ in $[0, 2\pi]$ and insert into \eqref{eq:analyt_diff}. Then each mode $\fatu_\omega$ satisfies
\begin{equation}
   \frac{\partial \fatu_{\omega}(t)}{\partial t}=(-\omega^2\,\sigma^2\,\fatT+\fatA)\,\fatu_\omega(t).
\label{eq:omegaeq}
\end{equation}
The solution to the equation is 
\begin{equation}
   \fatu_{\omega}(t)=\fatS\,\exp(\fatLambda t)\,\fatS^{-1}\,u_{\omega}(t),
\label{eq:omegaeqsol}
\end{equation}
where $\fatS(\omega)=(\fats_1, \fats_2,\ldots, \fats_N)$ is the eigenvector matrix of $-\omega^2\,\sigma^2\,\fatT+\fatA$
and the eigenvalues $\lambda_j(\omega)$ are on the diagonal of $\fatLambda$.
By Gerschgorin's theorem for the eigenvalues of a matrix, the eigenvalues 
all satisfy $\Re \lambda_j\leqslant 0$ for $\omega=0$ and $\Re \lambda_j< 0$ for $\omega>0$. Thus, as $t$ increases
all modes vanish except for one mode $\fats_1(0)=\fatui$ at $\omega=0$ with eigenvalue $\lambda_1= 0$. 

The concentrations of the internal states are summed at the macroscopic level. Then $U=\fate^T\fatu$ satisfies
\begin{equation}
\begin{array}{rl}
   \displaystyle{\frac{\partial U(x, t)}{\partial t}}&=\fate^T(\sigma^2\,\fatT\,\Delta\fatu +\fatA\,\fatu)=\sigma^2\,\fate^T\,\fatT\,\Delta\fatu=
    \sum_{\omega} -\omega^2\,\sigma^2\,\fate^T\,\fatT\,\fatu_{\omega}(t)\cos(\omega x)\\
     &=\bgamma(x, t)\,\fate^T\,\Delta \fatu=\bgamma(x, t)\,\Delta U,
\end{array}
\label{eq:Ueq}
\end{equation}
where
\begin{equation}
\displaystyle  \bgamma(x, t)=\sigma^2\,\frac{\sum\limits_{\omega=0}^{\infty} -\omega^2\,\fate^T\,\fatT\,\fatu_{\omega}(t)\cos(\omega x)}
                                 {\sum\limits_{\omega=0}^{\infty} -\omega^2\,\fate^T\,\fatu_{\omega}(t)\cos(\omega x)}.
\label{eq:gamma}
\end{equation}
The macroscopic $U$ satisfies a diffusion equation with a diffusion coefficient varying in space and time.
For large $t$, the dominant mode in the spatially non-constant part of the solution is damped by the eigenvalue
$\lambda_1=\max\limits_{j,\omega}\lambda_j(\omega)<0$ at $\omega_1$. Then $\fatu_{\omega}\approx \fats_1(\omega_1)\,e^{\lambda_1 t}$ and the steady-state macroscopic diffusion coefficient in \eqref{eq:gamma} is 
\begin{equation}
\displaystyle  \bgamma\approx\sigma^2\,\frac{\omega_1^2\,\fate^T\,\fatT\,\fats_1\,e^{\lambda_1 t}\cos(\omega_1 x)}
                                 {\omega_1^2\,\fate^T\,\fats_1\,e^{\lambda_1 t} \cos(\omega_1 x)}
                           =\sigma^2\,\frac{\fate^T\,\fatT\,\fats_1}{\fate^T\,\fats_1}=\sigma^2\,\frac{\sum\limits_{i=1}^N s_{1i}/\tau_{i}}{\sum\limits_{i=1}^N s_{1i}}.
\label{eq:gammaappr}
\end{equation}

\subsubsection{A reversible reaction}\label{sec:revreact}

Consider next the simple case of a reversible isomerization,
\begin{align}
  A_i \mathop{\rightleftharpoons}\limits_{\lambda_{ji}}^{\kappa_{ij}} B_j,
  \label{eq:unimolecular_reaction}
\end{align}
which is to be understood in the sense that, for example, the rate for the $j$th internal state of a $B$-molecule to transform into
the $i$th state of an $A$-molecule is $\lambda_{ji}$. These are thus two monomolecular reactions with reaction rates in Table~\ref{tab:coefficients}.

We assume that the  variance of the jump length, the diagonal matrix $\fatT$ and the vector $\fatmu$ are the same for all species with identical matrices $\fatD$ and $\fatA$.
Taking the mean as in \eqref{eq:Mommer_diff} we readily arrive at the coupled set of PDEs
\begin{align}
  \frac{\partial \fatu}{\partial t} &= 
  \fatD\, \Delta \fatu + \fatA\,\fatu-\fatK_1\,\fatu+\fatL_2\,\fatv,
  \label{eq:RDPDE_internalstates1} \\
  \frac{\partial \fatv}{\partial t} &= 
  \fatD\, \Delta \fatv + \fatA\,\fatv+\fatK_2\,\fatu-\fatL_1\,\fatv,
  \label{eq:RDPDE_internalstates2}
\end{align}
for some matrices $\fatK_{1}$, $\fatK_{2}$, and $\fatL_{1}$,
$\fatL_{2}$ whose precise form we now determine. Define positive rate matrices $\fatK$ and $\fatL$ with $K_{ij} = \kappa_{ij}$ and $L_{ij} = \lambda_{ij}$. From the
prescription \eqref{eq:unimolecular_reaction} we find that in
\eqref{eq:RDPDE_internalstates1}--\eqref{eq:RDPDE_internalstates2},
the $i$th state is affected by the reaction terms
\begin{align}
  (-\fatK_1\,\fatu+\fatL_2\,\fatv)_{i} &= 
  -\sum\limits_{j=1}^N \kappa_{ij}\,u_{i}+\sum\limits_{j=1}^N \lambda_{ji}\,v_{j}, \\
  (+\fatK_2\,\fatu-\fatL_1\,\fatv)_{i} &= 
  +\sum\limits_{j=1}^N \kappa_{ji}\,u_{j}-\sum\limits_{j=1}^N \lambda_{ij}\,v_{i}.
\end{align}
Identifying terms we obtain that
\begin{align}
  \label{eq:formK}
  \fatK_{1} &= \diag(\fatK\,\fate), \qquad 
  \fatK_{2} = \fatK^{T}, \\
  \label{eq:formL}
  \fatL_{1} &= \diag(\fatL\,\fate), \qquad
  \fatL_{2} = \fatL^{T},
\end{align}
where the notation $\diag(\fatf)$ denotes a diagonal matrix with $f_i$ on the diagonal.

We are interested in the stable, space independent solutions to 
\eqref{eq:RDPDE_internalstates1}--\eqref{eq:RDPDE_internalstates2}. Introduce
\begin{align}
  \fatB &= \left(\begin{array}{cc}
    \fatA-\fatK_1 &\fatL_2\\
    \fatK_2 & \fatA-\fatL_1
  \end{array}\right).
  \label{eq:Bdef}
\end{align}
By construction and using \eqref{eq:kernelA} we have the
crucial properties that
\begin{align}
  \label{eq:crucial1}
  B_{ii} &= -\displaystyle\sum_{k=1 \atop k\not = i}^{M} B_{ki}, \\
  \label{eq:crucial2}
  B_{ij} &\geqslant 0, \qquad i \not = j. \\
  \intertext{with $M=2N$. Also, by inspection $\fatB$ is \emph{irreducible}, i.e. there
    is no permutation matrix $\fatQ$ such that}
  \label{eq:crucial3}
  \fatQ \fatB \fatQ^{-1} &= \left( \begin{array}{cc}
    \fatE & \fatF \\
    \bm{0} & \fatG \end{array} \right).
\end{align}
Taken together, $\fatB$ is a \emph{$\mathbb{W}$-matrix} and this
provides us with certain general stability properties.

\begin{lemma}[\emph{$\mathbb{W}$-matrix lemma}]
  \label{lem:Wmatrix}
  Suppose that a real, irreducible matrix $\fatH\in\mathbb{R}^{M\times
    M}$ satisfies \eqref{eq:crucial1}--\eqref{eq:crucial3} (with
  $\fatB$ replaced by $\fatH$). Then the system of ODEs
  \begin{align}
    \label{eq:BODE}
    \calu'(t) &= \fatH\,\calu(t)
  \end{align}
  has a single stable equilibrium solution $\calui$ as $t \to
  \infty$. Moreover, if initial data $\calu(0)$ with positive mass
  $\fate^T\calu(0) > 0$ is given, then $\fate^T\calu(t)=\fate^T\calu(0) > 0$ 
  for all $t>0$.
\end{lemma}

This particular formulation is discussed in detail in
\cite[Chap.~V.3]{VanKampen} and we note that it can also be shown to
follow from the Perron-Frobenius theorem. To get some further insight
into the stability we present some alternative arguments as outlined
in \cite[Chap.~V.9, p.~129]{VanKampen}.

\begin{proof}
  Consider the adjoint problem
  \begin{align}
    \calv'(t) &= \fatH^{T}\,\calv(t), \\
    \intertext{such that}
    \label{eq:adjrel}
    \calv(0)^{T}\,\calu(t) &= \calv(t)^{T}\,\calu(0). \\
    \intertext{Using \eqref{eq:crucial1} we see that}
    \eta'_{i}(t) &= \sum_{k=1 \atop k \not = i}^{M}  H_{ki}\,(\eta_{k}-\eta_{i})(t).
  \end{align}
  By the irreducibility of $\fatH$, the largest element in $\calv$
  decreases and the smallest element increases with time such that an
  all-constant vector $\calv_\infty$ emerges in the limit. With $\calv(0) = \fate_{i}$,
  the $i$th unit vector, we recover from \eqref{eq:adjrel} the unique
  equilibrium solution $\tcalui$. Also, with $\calv(0) = \fate$ we have $\calv(t)=\fate$ and
  mass is a preserved quantity.
\end{proof}


Apply the lemma to $\fatB$ in \eqref{eq:Bdef} with $\calu^T=(\fatu^T\, , \fatv^T)$. Then it follows that
there is a steady state $(\tfatui\, , \tfatvi)$ and the initial mass $\fate^T\fatu(0)+\fate^T\fatv(0)$ is conserved for all $t$.

We now turn our attention to the equivalent reaction rates as induced by the
subdiffusive reactions \eqref{eq:unimolecular_reaction}. We write the unique equilibrium solution $\tfatui, \tfatvi$ in the normalized form $\tfatui = u\,\fatui$, $\tfatvi = v\,\fatvi$, where $\left|\left| \fatui \right|\right|_1 = \left|\left| \fatvi \right|\right|_1 = 1$. Using
\eqref{eq:formK}--\eqref{eq:formL} we readily find in good agreement
that
\begin{equation}\label{eq:kldef}
  k_{eq} := \fate^{T}\,\fatK_{1}\,\fatui = \fate^{T}\,\fatK_{2}\,\fatui, \quad
  l_{eq} := \fate^{T}\,\fatL_{1}\,\fatvi = \fate^{T}\,\fatL_{2}\,\fatvi. 
\end{equation}
Also, to mention just one example,
\begin{equation}\label{eq:kldef2}
  k_{eq}\,\fate^{T}\,\tilde{\fatu}_\infty = \fate^{T}\,\fatK_{1}\,\tilde{\fatu}_\infty,
\end{equation}
which can be understood as a kind of consistency result; at
steady state the equivalent reaction rate applied to the sum of the internal
states gives the same result as the sum of the individual subdiffusive
rates.

To conclude, using the subdiffusion steady state solutions, which
remain valid also when the coupling reactions
\eqref{eq:unimolecular_reaction} are `turned on', we can read off the
equivalent reaction rates as a function of the corresponding subdiffusive
rates. Insert $\fatu(\fatx, t)=u(\fatx, t)\,\fatui$ and $\fatv(\fatx, t)=v(\fatx, t)\,\fatvi$ into \eqref{eq:RDPDE_internalstates1} and \eqref{eq:RDPDE_internalstates2}.
Then we have
\begin{align}\label{eq:RDPDE_internalstates3}
  \fatui\,\frac{\partial u}{\partial t} &= 
  \fatD\,\fatui\,\Delta u + u\,\fatA\,\fatui-u\,\fatK_1\,\fatui+v\,\fatL_2\,\fatvi,
   \\
  \fatvi\,\frac{\partial v}{\partial t} &= 
  \fatD\,\fatvi\,\Delta v + v\,\fatA\,\fatvi+u\,\fatK_2\,\fatui-v\,\fatL_1\,\fatvi.
\end{align}
Multiply by $\fate^T$ and in the long-time limit we thus arrive at the familiarly looking macroscopic
reaction-diffusion PDE
\begin{align}\label{eq:RDPDE_internalstates3_inf1}
  \displaystyle{\frac{\partial u}{\partial t}} & = \gamma_u\, \Delta u - k_{eq}\,u + l_{eq}\,v, \\
  \displaystyle{\frac{\partial v}{\partial t}}  &= \gamma_v\, \Delta v + k_{eq}\,u - l_{eq}\,v, \label{eq:RDPDE_internalstates3_inf2}
\end{align}
with $\gamma_u=\fate^T\,\fatD\,\fatui$ and $\gamma_v=\fate^T\,\fatD\,\fatvi$.

\subsubsection{A bimolecular reaction}

We finally consider the case of a reversible dimerization,
\begin{align}
  A_i + B_j \mathop{\rightleftharpoons}\limits_{\lambda_{ijk}}^{\kappa_{ijk}} C_k.
  \label{eq:bimolecular_reaction}
\end{align}
As in Section~\ref{sec:revreact}, the jump length variance, the diagonal matrix $\fatT$ and the vector $\fatmu$ are the same for each species
and $\fatD$ and $\fatA$ are also identical.
Due to the nonlinearities it is inconvenient to write this in a matrix
form as in
\eqref{eq:RDPDE_internalstates1}--\eqref{eq:RDPDE_internalstates2}. Instead,
for each internal state $i$ we have the mean-field equations,
\begin{align}
  \label{eq:RDPDE2_internalstates1}
  \frac{\partial u_i}{\partial t} &= D_i\, \Delta u_i+\sum_{j=1}^N A_{ij}\,u_j
  -\sum_{j,k=1}^N K_{ijk}\,u_i\,v_j+\sum_{j,k=1}^N L_{ijk}\,w_k,\\
  \label{eq:RDPDE2_internalstates2}
  \frac{\partial v_i}{\partial t} &= D_i\, \Delta v_i+\sum_{j=1}^N A_{ij}\,v_j
  -\sum_{j,k = 1}^N K_{jik}\,u_j\,v_i+\sum_{j,k=1}^N L_{jik}\,w_k,\\
  \label{eq:RDPDE2_internalstates3}
  \frac{\partial w_i}{\partial t} &= D_i\, \Delta w_i+\sum_{j=1}^N A_{ij}\,w_j
  +\sum_{j,k=1}^N K_{jki}\,u_j\,v_k-\sum_{j,k=1}^N L_{jki}\,w_i,
\end{align}
where for readability we write $K_{ijk} = \kappa_{ijk}$ and $L_{ijk} =
\lambda_{ijk}$.
The difference $u_i-v_i$ satisfies
\begin{equation}\label{eq:uvdiff}
   \displaystyle \frac{\partial (u_i-v_i)}{\partial t} =  D_i\Delta (u_i-v_i) + \sum_{j=1}^N A_{ij} (u_i-v_i) - \sum_{j,k=1}^N K_{ijk}u_i v_j-K_{jik}u_j v_i.
\end{equation} 
If $K_{ijk}=K_{jik}$ and $L_{ijk}=L_{jik}$ in \eqref{eq:uvdiff} and $u_i(\fatx, 0)=v_i(\fatx, 0)$  for all $i$, then it follows that $u_i(\fatx, t)=v_i(\fatx, t)$ for $t>0$.

Assume that there is a positive steady state solution $(\tfatui\, , \tfatvi\, , \tfatwi)$ to \eqref{eq:RDPDE2_internalstates1}--\eqref{eq:RDPDE2_internalstates3} 
when $t\rightarrow \infty$. Let $(\tfatui, \tfatvi, \tfatwi)=(u\,\fatui\, , v\,\fatvi\, , w\,\fatwi)$ where $\|\fatui\|_1=\| \fatvi\|_1=\| \fatwi\|_1=1$.
Small perturbations around the equilibrium solution are evolved by the
Jacobian of the system. For
\eqref{eq:RDPDE2_internalstates1}--\eqref{eq:RDPDE2_internalstates3},
it can be written in the form
\begin{align}
  \label{eq:B2def}
  \fatB &= \left(\begin{array}{ccc}
    \fatA-\fatK_{11}&-\fatK_{12} &\fatL_1 \\
    -\fatK_{21} & \fatA-\fatK_{22} & \fatL_2 \\
    \fatK_{31} & \fatK_{32} & \fatA-\fatL_3
  \end{array}\right).
\end{align}
We identify after some tedious work
\begin{align}
  \fatK_{11} &= \diag \; \sum_{j,k=1}^N K_{ijk}\,v_{\infty j}, \quad\label{eq:mat_bimolecular}
  \fatK_{22} = \diag \; \sum_{j,k=1}^N K_{jik}\,u_{\infty j}, \quad
  \fatL_{3} = \diag \; \sum_{j,k=1}^N L_{jki}, \\
  (\fatL_{1})_{ij} &= \sum_{k=1}^N L_{ikj}, \quad 
  (\fatL_{2})_{ij} = \sum_{k=1}^N L_{kij}, \\
  (\fatK_{12})_{ij} &= \sum_{k=1}^N K_{ijk}\,u_{\infty i}, \quad
  (\fatK_{21})_{ij} = \sum_{k=1}^N K_{jik}\,v_{\infty i}, \\
  (\fatK_{32})_{ij} &= \sum_{k=1}^N K_{kji}\,u_{\infty k}, \quad 
  (\fatK_{31})_{ij} = \sum_{k=1}^N K_{jki}\,v_{\infty k}.
\end{align}

\emph{The Jacobian $\fatB$ in \eqref{eq:B2def} is not a $\mathbb{W}$-matrix.} To
understand why, note in \eqref{eq:bimolecular_reaction} that the total
sum of molecules is not a preserved quantity. However, if each
$C$-molecule is counted twice, this new weighted sum is in fact
preserved. Indeed, by inspection and after some work we find that
multiplying $\fatB$ by a diagonal matrix with 1 on the diagonal in the first $2N$ rows and 2 in the last $N$ rows we have a
$\mathbb{W}$-matrix. Using the $\mathbb{W}$-matrix
Lemma~\ref{lem:Wmatrix} it therefore follows that small perturbations
around an equilibrium solution are stable. 

Equivalent reaction rates can be defined as follows,
\begin{equation}\label{eq:kldef3}
  k_{eq} = \sum_{i,j,k=1}^N K_{ijk}\,u_{\infty i}\,v_{\infty j}, \quad
  l_{eq} = \sum_{i,j,k=1}^N L_{ijk}\,w_{\infty k}.
\end{equation}
Insert $\fatu(\fatx, t), \fatv(\fatx, t),$ and $\fatw(\fatx, t)=w(\fatx, t)\,\fatwi$ into 
\eqref{eq:RDPDE2_internalstates1}--\eqref{eq:RDPDE2_internalstates3} as in \eqref{eq:RDPDE_internalstates3} and let $\gamma_w=\fate^T\,\fatD\,\fatwi$. 
As expected we recover the reaction-diffusion PDE for $u, v,$ and $w$
\begin{align}\label{eq:bimolintstates1}
  \frac{\partial u}{\partial t}&= \gamma_u\, \Delta u - k_{eq}\,u\,v + l_{eq}\,w, \\
  \frac{\partial v}{\partial t}&= \gamma_v\, \Delta v - k_{eq}\,u\,v + l_{eq}\,w, \label{eq:bimolintstates2}\\
  \frac{\partial w}{\partial t}&= \gamma_w\, \Delta w + k_{eq}\,u\,v - l_{eq}\,w.\label{eq:bimolintstates3}
\end{align}



\section{The internal states approximation of the FPDE}

In this section we start off with the promising observation made in
\cite{MOMMER09} that non-Markovian waiting times can be arbitrarily well
approximated by a set of Markovian waiting times, each associated with
its own \emph{internal state}. In turn, those states are to be visited
according to a certain random walk model which again can be taken as
Markovian, all in all resulting in a computationally quite attractive
modeling framework for subdiffusion and reactions. 

In Section~\ref{subsec:internal_states}, we recapitulate the basic internal
states subdiffusive model and its relation to the FPDE, and in Section~\ref{subsec:asymptotic} we
determine its asymptotic behavior for short and long times. In
Section~\ref{subsec:internal_states_reaction}, we consider the same framework in the
presence of reactions with two different mesoscopic models.
In particular, we
discuss the feasibility of obtaining coarse-grained macroscopic coefficients in a FPDE from
observations of subdiffusive systems. 

 
\subsection{Internal states diffusion system}
\label{subsec:internal_states}

The asymptotic behavior of the waiting time
\eqref{eq:waiting_time_subdiffusion} at large time follows the power
law
\begin{equation}
\psi(t) \approx A_{\alpha}\,\frac{\tau^{\alpha}}{t^{1+\alpha}},
\label{eq:waitingtime_alpha}
\end{equation}
with
\begin{equation}
A_{\alpha} = \frac{\sin(\pi\,\alpha)}{\pi}\,\Gamma(1+\alpha).
\label{eq:Aalpha}
\end{equation}
With a change of variable in the Euler's $\Gamma$ function, the diffusive representation of the totally monotone
function $\frac{1}{t^{1+\alpha}}$ in \eqref{eq:waitingtime_alpha} is
\cite{DESCH88,HADDAR10,HELESCHEWITZ00,STAFFANS94}
\begin{equation}
\frac{1}{t^{1+\alpha}} = \frac{1}{\Gamma(1+\alpha)}\,\int_0^{\infty}s^{\alpha}\,e^{-s\,t}\,ds.
\label{eq:diffusive_representation}
\end{equation}
The diffusive representation of $\frac{1}{t^{1+\alpha}}$ \eqref{eq:diffusive_representation} is 
approximated by using a quadrature formula in $N$ points, with weights $\tilde{\mu}_{i}$ and abscissae $s_{i}$:
\begin{equation}
\frac{1}{t^{1+\alpha}} \simeq \sum\limits_{i=1}^N \tilde{\mu}_{i}\,e^{-s_{i}\,t},
\label{eq:diffusive_approximation1}
\end{equation}
leading to the {\it diffusive approximation} \cite{BLANC15}. 

Our objective is then to approximate the function $F_{ex}(t) =
\frac{1}{t^{1+\alpha}}$ by $F_{approx}(t) = \sum\limits_{i=1}^N
\tilde{\mu}_i\,e^{-s_i t}$ in the time interval $\left[
  t_{min},t_{max}\right]$. One possibility to quantify the error of
the model $\varepsilon_{mod}$ is
\begin{equation}
  \varepsilon_{mod} = \left|\left| \frac{F_{approx}(t)}{F_{ex}(t)} - 1 \right|\right| = \left( \frac{1}{t_{max}-t_{min}}\int_{t_{min}}^{t_{max}}\left| \frac{F_{approx}(t)}{F_{ex}(t)} - 1 \right|^2dt\right)^{1/2}.
  \label{eq:error_model}
\end{equation}
Based on the error \eqref{eq:error_model}, a nonlinear optimization is
shown in \cite{BLANC15} to be a better way to determine the parameters
$\tilde{\mu}_i$ and $s_i$ than Gaussian quadrature. Consequently, this
method is used in all what follows.

Setting 
\begin{equation}
\tau_{i} = 1/s_{i},\quad \tau = \left(A_{\alpha}\,\sum\limits_{i=1}^N \frac{\tilde{\mu}_i}{s_i}\right)^{-1/\alpha},\quad \mu_{i} = \frac{\tilde{\mu}_i}{s_i}\,\left(\sum\limits_{i=1}^N \frac{\tilde{\mu}_i}{s_i}\right)^{-1},
\end{equation}
the approximation of the waiting time \eqref{eq:waitingtime_alpha} is
\begin{equation}\label{eq:waitapprox}
\psi(t) \, \approx \, \sum\limits_{i=1}^N \mu_{i}\,\tau_{i}^{-1}\,e^{-t/\tau_{i}}.
\end{equation}
The jump length variance $\Sigma^2$ \eqref{eq:length_charac_ctrw} is then computed using \eqref{eq:Kalpha}. 

\begin{figure}[htbp]
\begin{center}
\begin{tabular}{c}
\includegraphics[scale=0.35]{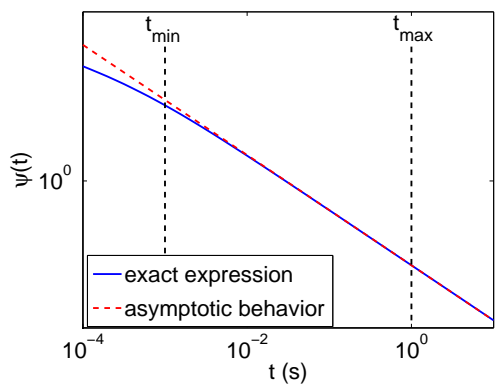}
\end{tabular}
\end{center}
\caption{Section~\ref{subsec:internal_states}. Waiting time PDF \eqref{eq:waiting_time_subdiffusion} (blue solid line), and its asymptotic expansion \eqref{eq:waitingtime_alpha} (red dotted line).
The scales are logarithmic on both axes.}
\label{fig:waiting_time}
\end{figure}

The waiting time PDF \eqref{eq:waiting_time_subdiffusion} and its asymptotic expansion \eqref{eq:waitingtime_alpha} are compared in Figure~\ref{fig:waiting_time}. 
The parameters are those used in the numerical experiments in Section~\ref{sec:exp:config} (Table~\ref{table:parameters}, second set of parameters). 
Anomalous diffusion is expected in the time range of interest $\left[ t_{min},t_{max}\right]$.
In this time interval, the expansion \eqref{eq:waitingtime_alpha} is already accurate as illustrated in Figure~\ref{fig:waiting_time}. 

The CTRW algorithm in Section~\ref{sec:background} is extended to a multistate CTRW (MCTRW) algorithm in \cite{MOMMER09} for a diffusive system with
internal states. One diffusive jump in the algorithm is performed as follows. 
Firstly, a state $i$ is drawn  with probability $\mu_i$. Secondly, the waiting time is sampled from an exponential 
distribution with the PDF $\psi_i(t)=\tau_i^{-1}e^{-t/\tau_i}$.
The PDF of the jump length in an internal state $i$ is $\lambda_i(\fatx)$ as in \eqref{eq:jump_length_gaussian} with variance $\sigma^2_i$. 
Finally, the length of the jump is drawn with the normal distribution of $\lambda_i$. 

The combined PDFs for the jump length and the waiting time of the internal states are
\begin{equation}
\lambda(\fatx) = \sum\limits_{i=1}^N \mu_{i}\,\lambda_{i}(\fatx),\quad \tilde{\psi}(t)=\sum_{i=1}^N\mu_i\,\psi_i(t)=\sum_{i=1}^N\mu_i\tau_i^{-1}e^{-t/\tau_i}.
\label{eq:approx_jumplength}
\end{equation}
The reaction-diffusion system \eqref{eq:analyt_diff} in \cite{MOMMER09} is derived from the MCTRW algorithm and the approximation \eqref{eq:waitapprox}. 
In the $i$th internal state, the diffusion of $u_{i}$ is ordinary, whereas the sum of the concentrations of the internal states at the macroscopic level $U = \fate^T\,\fatu$ diffuses anomalously within the time range
$[t_{min},t_{max}]$. A collection of internal states with waiting times $\psi_i(t)$ are here approximated by one state with the waiting time $\psi(t)$ in \eqref{eq:waitapprox}.



\subsection{Diffusive behavior of the internal states model}
\label{subsec:asymptotic}

To mimic a subdiffusive behavior and to define a Markov process, the waiting time PDF $\psi(t)$ is approximated by a sum of $N$
exponentials \eqref{eq:waitapprox}. Let us assume that $\tau_1 < \tau_2 < \cdots < \tau_N$.
Equation \eqref{eq:waitapprox} implies the following.
\begin{itemize}
\item At small times, the Taylor expansion of $\psi(t)$ is
\begin{equation}
\begin{array}{lll}
\psi(t)  \displaystyle \approx \tilde{\psi}(t)=\sum\limits_{i=1}^N \mu_{i}\,\tau_{i}^{-1}\,e^{-t/\tau_{i}}
\approx \sum\limits_{i=1}^N \mu_{i}\,\tau_{i}^{-1}\,\left(1 - \frac{t}{\tau_{i}}\right)
\mathop{\sim}\limits_{t\rightarrow 0^+}  1 - \frac{t}{\tau_{eq}}
\mathop{\sim}\limits_{t\rightarrow 0^+}  e^{-t/\tau_{eq}},
\end{array}
\label{eq:small_time}
\end{equation}
with $\tau_{eq} = \left(\sum\limits_{i=1}^N
\mu_{i}\,\tau_{i}^{-1}\right)/\left(\sum\limits_{i=1}^N
\mu_{i}\,\tau_{i}^{-2}\right)$. At small time, $\psi(t)$ is therefore
equivalent to a Poisson law with parameter $\tau_{eq}$. Hence,
ordinary diffusion is expected for $t\ll\tau_1$;
\item As we see in \eqref{eq:Ueq}, the macroscopic $U$ satisfies a diffusion equation with a diffusion coefficient $\bgamma(x, t)$ \eqref{eq:gamma} varying in space and time. For large $t$, the long-time diffusion coefficient is given by \eqref{eq:gammaappr}. It does not depend on space and time. Consequently, a return to ordinary diffusion is expected.
\item For $t\in [\tau_1, \tau_N]$, a subdiffusive behavior is expected with subdiffusive exponent $\alpha$.
\end{itemize}
This is in agreement with observations in physical experiments \cite{HofFra,JEON11,KUSUMI05,Saxton07}. 
Numerical illustrations are found in Section~\ref{sec:num_exp}.

\subsection{Internal states reaction-diffusion system}
\label{subsec:internal_states_reaction}

In this section, the mesoscopic model with $N$ internal states of the
participating molecules in Section~\ref{sec:internalstates} is
compared to the macroscopic FPDE models I and II in
Section~\ref{subsec:reaction_modeling}. The annihilation reaction,
monomolecular reactions and bimolecular reactions summarized in
Table~\ref{tab:coefficients} are investigated. For mesoscopic models
with certain reactions, there are corresponding FPDE models but in
other cases the macroscopic level with summation over the internal
states is not so easily expressed as a FPDE.

\subsubsection{Annihilation process}\label{sec:IS_annihilation}

We consider one species $A$ with $N$ different internal states. Each internal state can be annihilated with rate $k_i$
\begin{equation}
A_i \mathop{\longrightarrow}^{k_i} \emptyset, \quad i=1,\cdots,N.
\label{eq:decay_IS}
\end{equation}
The internal states diffusion system \eqref{eq:RDPDE_internalstates1} without $\fatv$ is then modified as follows
\begin{equation}
\frac{\partial \fatu}{\partial t} = \fatD\,\Delta \fatu + \fatA\,\fatu - \fatK_1\,\fatu.
\label{eq:PDE_internalstates_reaction}
\end{equation}

Let $k_i=k/\tau_i$ in \eqref{eq:decay_IS}
as in Table~\ref{tab:coefficients}. 
Then $\fatK_1=k\,\fatT$ and \eqref{eq:PDE_internalstates_reaction} can be written
\begin{equation}
  \frac{\partial \fatu}{\partial t} =\fatT\,(\sigma^2\,\Delta \fatu - k\,\fatu) + \fatA\,\fatu.
\label{eq:PDE_internalstates_new}
\end{equation}
In model I \eqref{eq:FPDE_decayI}, the fractional order temporal derivative $\frac{\partial^{1-\alpha}}{\partial t^{1-\alpha}}$ acts on both the Laplace operator and the reaction term. In the internal states diffusion system \eqref{eq:PDE_internalstates_reaction} and \eqref{eq:PDE_internalstates_new}, the diffusion matrix is scaled by the waiting time matrix $\fatT$ to mimic an anomalous behavior. Scaling the reaction matrix $\fatK_1$ in the same way will approximate the macroscopic model I. The waiting time $\psi_i(t)$ for diffusive and reactive events is the same for all internal states. The relation between $k$ in \eqref{eq:PDE_internalstates_new} and $k_\ast$ in \eqref{eq:FPDE_decayI} is the same as between $\sigma^2$ and $K_\alpha$ in \eqref{eq:Kalpha},
\begin{equation}
\frac{k_\ast}{k} =\frac{\sigma^2}{K_\alpha} = \tau^\alpha.
\label{eq:FPDEmesocoeff}
\end{equation}

For an approximation of model II, take $\fatK$ and $\fatK_1$ to be $k\,\fatI$ in \eqref{eq:formK}. 
This corresponds to an annihilation case in Table~\ref{tab:coefficients} with all $\tau_i=1$. The right hand side of \eqref{eq:PDE_internalstates_reaction} is with this $\fatK_1$
\begin{equation}
  \frac{\partial \fatu}{\partial t} =\sigma^2\,\fatT\,\Delta \fatu + \fatA\,\fatu - k\,\fatu.
\label{eq:PDE_internalstates_RHS2}
\end{equation}
The waiting time $\psi_i(t)$ for diffusive and reactive events is different for the internal states. 

If $\fatK_1$ and $\fatT$ and $\fatK_1$ and $\fatA$ commute then the equation for ${\fatu}$ in \eqref{eq:PDE_internalstates_RHS2} can be written
\begin{equation}
\frac{\partial {\fatu}}{\partial t} = \sigma^2\,e^{\fatK_1t}\,\fatT\,e^{-\fatK_1t}\,\Delta {\fatu} + e^{\fatK_1t}\,\fatA\,e^{-\fatK_1t}\,{\fatu} - k\,\fatu.
\label{eq:PDE_internalstates_reaction_a}
\end{equation}
This equation corresponds to model II in Section~\ref{subsubsec:reaction_modeling_annihilation} at the macroscopic level. 
The relation between $k$ in \eqref{eq:PDE_internalstates_RHS2} and $k_\ast$ in \eqref{eq:FPDE_decayII} is
\begin{equation}
  k_\ast = k.
  \label{eq:FPDEmesocoeff2}
\end{equation}
After a change of variables $\tilde{U} = e^{k_\ast t}\,U$ in the model II FPDE \eqref{eq:FPDE_decayII}, the FPDE \eqref{eq:FPDE_monoII} is obtained. In the same manner, introduce a change of variables $\tilde{\fatu} = e^{\fatK_1\,t}\,\fatu$ in \eqref{eq:PDE_internalstates_reaction}. Then
\begin{equation}
\frac{\partial \tilde{\fatu}}{\partial t} = \sigma^2\,\fatT\,\Delta \tilde{\fatu} + \fatA\,\tilde{\fatu}.
\label{eq:PDE_internalstates_reaction_b}
\end{equation}
A sufficient condition for $\fatK_1$ to commute with $\fatT$ and
$\fatA$ is that $\fatK_1=k\,\fatI$ as is the case in \eqref{eq:PDE_internalstates_RHS2}. If $\fatK_1$ does
not commute with $\fatT$ and $\fatA$, the macroscopic equation for $U=\fate^T\,\fatu$ may not be as simple as \eqref{eq:FPDE_monoII}.

Consider $\bfatu$, the total amount of $A$ in the different states in $\Omega$, defined by
\begin{equation}
  \bfatu(t)=\int_{\Omega}\fatu(\fatx, t)\, d\Omega.
\label{eq:ubardef}
\end{equation}
Integrating \eqref{eq:PDE_internalstates_reaction} and using the Neumann boundary condition in \eqref{eq:PDE_Na} and \eqref{eq:neumannBC} leads to
\begin{equation}
\displaystyle \frac{d \bfatu}{dt}=\fatA\,\bfatu-\fatK_1\,\bfatu.
\label{eq:ubareq}
\end{equation}
The time evolution of $\bfatu$ is then
\begin{equation}
\bfatu(t) = e^{(\fatA-\fatK_1)\,t}\,\bfatu(0).
\label{eq:ubar}
\end{equation}
By Gerschgorin's theorem for the eigenvalues of a matrix, the eigenvalues of $\fatA-\fatK_1$ are strictly in the left half plane.
By \eqref{eq:ubareq}, the equation for the sum over all internal states $\bar{U}=\fate^T\bfatu$ is 
\begin{equation}
\frac{d \bar{U}}{dt} = -\fate^T\,\fatK_1\,\bfatu=-k'(t)\,\bar{U},
\label{eq:totubareq}
\end{equation}
with
\begin{equation}
   k'(t) = \frac{\fate^T\,\fatK_1\,\bfatu}{\fate^T\,\bfatu} = \frac{\fate^T\,\fatK_1\,e^{(\fatA-\fatK_1)\,t}\,\bfatu(0)}{\fate^T\,e^{(\fatA-\fatK_1)\,t}\,\bfatu(0)}
\label{eq:kprim}
\end{equation}
In model I with $\fatK_1=k\fatT$, $k'$ varies in time depending on $\bfatu$, and the kinetics is then anomalous. 
Let $\lambda_1$ be the eigenvalue of $\fatA-\fatK_1$ with maximum real part and $\fats_1$ the corresponding eigenvector.
Then for large $t$
\begin{equation}
k'(t) \approx \frac{\fate^T\,\fatK_1\,\fats_1 e^{\lambda_1 t}}{\fate^T\,\fats_1 e^{\lambda_1 t}} = k\frac{\sum_i s_{1i}/\tau_i}{\sum_i s_{1i}},
\label{eq:kprim2}
\end{equation}
cf. \eqref{eq:gammaappr}. On the contrary, in model II $k'=k$ does not depend on time, and the kinetics is then ordinary. This conclusion agrees with the comments in Section~\ref{subsubsec:reaction_modeling_annihilation}.

\subsubsection{Monomolecular reactions}\label{sec:IS_mono}

We study two species $A$ and $B$ with $N$ different internal states and the reversible reactions in \eqref{eq:unimolecular_reaction}. 
The internal states diffusion system is found in \eqref{eq:RDPDE_internalstates1} and \eqref{eq:RDPDE_internalstates2}.
The aim of this section is to choose the reaction matrices $\fatK$ and $\fatL$ in \eqref{eq:formK} and \eqref{eq:formL} to mimic either 
model I \eqref{eq:FPDE_monoI} or model II \eqref{eq:FPDE_monoII} and to discuss other alternatives. 

The reaction rates are first chosen to be scaled by the waiting time matrix and one state $i$ is transformed to the same state $i$ in
\begin{align}
  A_i \mathop{\rightleftharpoons}\limits_{\ell/\tau_i}^{k/\tau_i} B_i,\; i=1,\ldots,N.
  \label{eq:unimolecular_reaction2}
\end{align}
Then $K_{ii} = k/\tau_{i}$ and $L_{ii} = \ell/\tau_{i}$, cf. Table~\ref{tab:coefficients}. 
With $\fatK=k\,\fatT$ and $\fatL=\ell\,\fatT$ in 
\eqref{eq:RDPDE_internalstates1} and \eqref{eq:RDPDE_internalstates2} we have
\begin{equation}
\begin{array}{l}
\displaystyle \frac{\partial \fatu}{\partial t} = \fatT\,(\sigma^2\,\Delta \fatu - k\,\fatu + \ell\,\fatv)+ \fatA\,\fatu,\\
[10pt]
\displaystyle \frac{\partial \fatv}{\partial t} = \fatT\,(\sigma^2\,\Delta \fatv + k\,\fatu - \ell\,\fatv)+ \fatA\,\fatv.
\end{array}
\label{eq:PDE_internalstates_reaction3I}
\end{equation}
The diffusion and the reactions have the same waiting time as in model I in \eqref{eq:FPDE_monoI}.

The steady state of model I in \eqref{eq:PDE_internalstates_reaction3I} has an analytical solution. Let $\fatui$ span the nullspace of $\fatA$ (see \eqref{eq:uinfexpr}) and
insert the constant solutions in space $u\,\fatui$ and $v\,\fatui$ with $\|\fatui\|_1=1$ as $\fatu$ and $\fatv$ in \eqref{eq:PDE_internalstates_reaction3I}. Then
the right hand sides are
\begin{equation}
\begin{array}{rl}
\displaystyle{- k\,u\,\fatT\,\fatui + \ell\,v\,\fatT\,\fatui}&=0,\\
\displaystyle{k\,u\fatT\,\fatui - \ell\,v\,\fatT\,\fatui}&=0.
\end{array}
\label{eq:PDE_ststates_reaction}
\end{equation}
Both equations in \eqref{eq:PDE_ststates_reaction} are satisfied if $k\,u=\ell\,v$. Hence, the steady state solution is
$u(\fatui\, , (k/\ell)\,\fatui)$, where $u$ depends on the initial data. By mass conservation in Lemma~\ref{lem:Wmatrix} 
\begin{equation}
  \fate^T\fatu(0)+\fate^T\fatv(0)=\fate^T\fatu(t)+\fate^T\fatv(t)=u(\fate^T\fatui+(k/\ell)\fate^T\fatui)=u(1+k/\ell).
\label{eq:masscons}
\end{equation}

The macroscopic reaction coefficients can be regarded as time and space dependent and the macroscopic equations can be expressed without fractional derivatives. 
Introduce small perturbations $\delta\fatu$ and $\delta\fatv$
around the steady state in \eqref{eq:PDE_internalstates_reaction3I}
\[
  \fatu(\fatx, t)=u\,\fatui+\delta\fatu(\fatx, t),\quad \fatv(\fatx, t)=v\,\fatui+\delta\fatv(\fatx, t).
\]
Then for $U=\fate^T\,\fatu=u+\ordo(\|\delta \fatu\|)$ and $V=\fate^T\,\fatv=v+\ordo(\|\delta \fatv\|)$ we derive
\begin{equation}
\begin{array}{l}
\displaystyle \frac{\partial U}{\partial t} = \sigma^2\,\fate^T\,\fatT\,\Delta\delta\fatu - k\,u\,\fate^T\,\fatT\,\fatui- k\,\fate^T\,\fatT\,\delta\fatu + 
                                              \ell\,v\,\fate^T\,\fatT\,\fatui+\ell\,\fate^T\,\fatT\,\delta\fatv=-K'\,U+L'\,V+\ordo(\|\delta\fatu\|),\\
[10pt]
\displaystyle \frac{\partial V}{\partial t} = \sigma^2\,\fate^T\,\fatT\,\Delta\delta\fatv + k\,u\,\fate^T\,\fatT\,\fatui + k\,\fate^T\,\fatT\,\delta\fatu - 
                                              \ell\,v\,\fate^T\,\fatT\,\fatui-\ell\,\fate^T\,\fatT\,\delta\fatv=K'\,U-L'\,V+\ordo(\|\delta\fatv\|),
\end{array}
\label{eq:PDE_internalstates_reaction3inf}
\end{equation}
where
\begin{equation}
\begin{array}{l}
\displaystyle K'=k\,\frac{\fate^T\,\fatT\,(u\,\fatui+\delta\fatu)}{\fate^T\,(u\,\fatui+\delta\fatu)}=k\,\fate^T\,\fatT\,\fatui+\ordo(\|\delta\fatu\|),\\
[10pt]
\displaystyle L'=\ell\,\frac{\fate^T\,\fatT\,(v\,\fatui+\delta\fatv)}{\fate^T\,(v\,\fatui+\delta\fatv)}=\ell\,\fate^T\,\fatT\,\fatui+\ordo(\|\delta\fatv\|).
\end{array}
\label{eq:PDE_internalstates_reaction3coeff}
\end{equation}
The effective macroscopic reaction coefficients are $K'$ and $L'$ when $t\gg 0$, cf. \eqref{eq:kldef} with $\fatK_1=k\fatT$.

With more general $\fatK$ and $\fatL$ in \eqref{eq:unimolecular_reaction}, the structure of the equations is the same as in \eqref{eq:PDE_internalstates_reaction3I} 
but there is no corresponding macroscopic FPDE.
For example, let $\fatK=\fatT\,\fate\,\fatf^T\,\fatT$ and $\fatL=\fatT\,\fate\,\fatg^T\,\fatT$ with $\fatf^T\,\fatT\,\fate>0$ and $\fatg^T\,\fatT\,\fate>0$. Then $\fatK_1=k\,\fatT$ with $k=\fatf^T\,\fatT\,\fate$ and 
$\fatL_1=\ell\,\fatT$ with $\ell=\fatg^T\,\fatT\,\fate$. The equations for $\fatu$ and $\fatv$ are
\begin{equation}
\begin{array}{l}
\displaystyle \frac{\partial \fatu}{\partial t} = \fatT\,(\sigma^2\,\Delta\fatu - k\,\fatu + \fatg\,\fate^T\,\fatT\,\fatv)+ \fatA\,\fatu,\\
[10pt]
\displaystyle \frac{\partial \fatv}{\partial t} = \fatT\,(\sigma^2\,\Delta\fatv + \fatf\,\fate^T\,\fatT\,\fatu - \ell\,\fatv)+ \fatA\,\fatv.
\end{array}
\label{eq:PDE_internalstates_reaction3Igen}
\end{equation}
The macroscopic waiting time for the reactions is not obvious in this case. 

As in \eqref{eq:unimolecular_reaction2}, let one state $i$ be transformed to the same state $i$ and choose the reaction rates to be constant for all states $K_{ii} = k$ and $L_{ii} = \ell$ with $\tau_i=1$ in Table~\ref{tab:coefficients}. Insert $\fatK=k\,\fatI$ and $\fatL=\ell\,\fatI$ into \eqref{eq:RDPDE_internalstates1} and \eqref{eq:RDPDE_internalstates2} to obtain
\begin{equation}
\begin{array}{l}
\displaystyle \frac{\partial \fatu}{\partial t} = \sigma^2\,\fatT\,\Delta \fatu + \fatA\,\fatu- k\,\fatu + \ell\,\fatv,\\
[10pt]
\displaystyle \frac{\partial \fatv}{\partial t} = \sigma^2\,\fatT\,\Delta \fatv + \fatA\,\fatv+ k\,\fatu - \ell\,\fatv.
\end{array}
\label{eq:PDE_internalstates_reaction3II}
\end{equation}
Only the diffusion has different waiting times in different states while the waiting times for the reactions are the same in all states. 
Let $\fatM_r$ and $\fatM_d$ be the reaction matrix and the operator for the diffusion and the change of state in \eqref{eq:PDE_internalstates_reaction3II} 
\begin{equation}
\fatM_r = \left(
\begin{array}{cc}
-k\,\fatI & \ell\,\fatI\\
k\,\fatI & -\ell\,\fatI
\end{array}
\right),\quad
\fatM_d = \left(
\begin{array}{cc}
\sigma^2\,\fatT\,\Delta +\fatA & \bm{0} \\
\bm{0} & \sigma^2\,\fatT\,\Delta +\fatA
\end{array}
\right).
\label{eq:GHdef}
\end{equation}
There is a transformation from \eqref{eq:PDE_internalstates_RHS2} via \eqref{eq:PDE_internalstates_reaction_a} to \eqref{eq:PDE_internalstates_reaction_b}. 
Since $\fatM_r$ and $\fatM_d$ commute in \eqref{eq:GHdef}, 
there is a similar transformation for \eqref{eq:PDE_internalstates_reaction3II}. 
As in \eqref{eq:PDE_internalstates_RHS2},
$\fate^T\fatu$ and $\fate^T\fatv$ approximate model II in \eqref{eq:FPDE_monoII}.

If the diffusion rate is independent of the state but the waiting time for the reactions depends on the state in \eqref{eq:unimolecular_reaction2} then 
$D_{ii}=\sigma^2$ in \eqref{eq:analyt_diff}, $K_{ii}=k/\tau_i$ and $L_{ii}=\ell/\tau_i$. The mesoscopic model is
\begin{equation}
\begin{array}{l}
\displaystyle \frac{\partial \fatu}{\partial t} = \fatT\,(- k\,\fatu + \ell\,\fatv)+ \fatA\,\fatu+\sigma^2\,\Delta \fatu,\\
[10pt]
\displaystyle \frac{\partial \fatv}{\partial t} = \fatT\,(k\,\fatu - \ell\,\fatv)+ \fatA\,\fatv+\sigma^2\,\Delta \fatv.
\end{array}
\label{eq:PDE_internalstates_reaction3III}
\end{equation}
We have the ordinary diffusion but the reactions behave anomalously at the macroscopic level with the waiting time $\psi_i(t)$ depending on the state.


\subsubsection{Bimolecular reactions}

We consider three species $A$, $B$ and $C$ in $N$ different internal states and the chemical reactions
\begin{equation}
A_i + B_i \mathop{\rightleftharpoons}^{\kappa_i}_{\lambda_i} C_i, \quad i=1,\cdots,N.
\label{eq:bimo_IS}
\end{equation}
This is a less general set of reactions than in \eqref{eq:bimolecular_reaction} in that only molecules in the same internal state react with each other. 
The internal states reaction-diffusion system \eqref{eq:RDPDE2_internalstates1}-\eqref{eq:RDPDE2_internalstates3} is then modified as follows for state $i$
\begin{equation}
\begin{array}{lll}
\displaystyle \frac{\partial u_i}{\partial t} & = & \displaystyle D_i\, \Delta u_i + \sum_{j=1}^N A_{ij}\,u_j - K_{iii}\,u_i\,v_i+ L_{iii}\,w_i,\\
[14pt]
\displaystyle \frac{\partial v_i}{\partial t} & = & \displaystyle D_i\, \Delta v_i + \sum_{j=1}^N A_{ij}\,v_j - K_{iii}\,u_i\,v_i+ L_{iii}\,w_i,\\
[14pt]
\displaystyle \frac{\partial w_i}{\partial t} & = & \displaystyle D_i\, \Delta w_i + \sum_{j=1}^N A_{ij}\,w_j + K_{iii}\,u_i\,v_i - L_{iii}\,w_i.
\end{array}\label{eq:bimo_eq}
\end{equation}
The reaction coefficients in \eqref{eq:bimo_eq} are chosen to be $K_{iii}=k/\tau_i$ and $L_{iii}=\ell/\tau_i$. Then at state $i$
\begin{equation}
\begin{array}{lll}
\displaystyle \frac{\partial u_i}{\partial t} & = & \displaystyle \frac{1}{\tau_i}\left(\sigma^2\,\Delta u_i - k\,u_i\,v_i+ \ell\,w_i\right)+ \sum_{j=1}^N A_{ij}\,u_j,\\
[14pt]
\displaystyle \frac{\partial v_i}{\partial t} & = & \displaystyle \frac{1}{\tau_i}\left(\sigma^2\,\Delta v_i - k\,u_i\,v_i+ \ell\,w_i\right)+ \sum_{j=1}^N A_{ij}\,v_j,\\
[14pt]
\displaystyle \frac{\partial w_i}{\partial t} & = & \displaystyle \frac{1}{\tau_i}\left(\sigma^2\,\Delta w_i + k\,u_i\,v_i - \ell\,w_i\right)+ \sum_{j=1}^N A_{ij}\,w_j.
\end{array}\label{eq:bimo_eq2}
\end{equation}
The reactions and the diffusion have the same waiting time in each state in \eqref{eq:bimo_eq2} as in \eqref{eq:PDE_internalstates_reaction3I}. 
Hence, the equations will approximate model I at the macroscopic level in \eqref{eq:FPDE_bimoI}. 
By letting $K_{iii}=k$ and $L_{iii}=\ell$ independent of the internal state in \eqref{eq:bimo_eq}, equations
similar to the model II equations \eqref{eq:PDE_internalstates_reaction3II} are obtained but there is no transformation \eqref{eq:PDE_internalstates_reaction_a} 
with a constant $\fatK_1$.

We have found that certain annihilation, monomolecular, and bimolecular reactions at the mesoscopic level have a macroscopic counterpart as a FPDE.
Furthermore, the long time behavior of the reactions tends to that of ordinary reactions without internal states. This will be confirmed in
numerical examples in the next section.




\section{Numerical experiments}\label{sec:num_exp}

Since this paper focuses on the mesoscopic approximation of subdiffusion, we will investigate only the 1D
reaction-diffusion system. However, the strategy proposed here can be extended straightforwardly to 2D and 3D geometries.

The integration of the internal states systems \eqref{eq:analyt_diff}, \eqref{eq:RDPDE_internalstates1}--\eqref{eq:RDPDE_internalstates2} and \eqref{eq:RDPDE2_internalstates1}--\eqref{eq:RDPDE2_internalstates3} is detailed in Section~\ref{sec:num_modeling}. The general configuration of the numerical experiments is introduced in Section~\ref{sec:exp:config}: physical parameters, numerical parameters of the discretization, and initial conditions. Then the numerical experiments are described. In Section~\ref{sec:without_reaction}, both the diffusive approximation introduced in Section~\ref{subsec:internal_states} and the numerical method introduced in Section~\ref{sec:num_modeling} are verified in comparisons. In Sections~\ref{sec:reaction_decay} and \ref{sec:linear_reaction}, the method is applied to both model I and model II for reactive systems. The differences between the two models are also illustrated. Finally examples of bimolecular reactions, for which no analytical solutions are available, are presented in Section~\ref{sec:nonlinear_reaction}.

\subsection{Numerical modeling}\label{sec:num_modeling}


In order to integrate the internal states systems \eqref{eq:analyt_diff}, \eqref{eq:RDPDE_internalstates1}--\eqref{eq:RDPDE_internalstates2} and \eqref{eq:RDPDE2_internalstates1}--\eqref{eq:RDPDE2_internalstates3}, a uniform grid is introduced with mesh size $h$ and time step $\Delta t$. The
approximation of the exact solution $\fatu$ is denoted by $\fatu_j^n$ at $x_j$ and $t^n$. The Laplace operators involved in the internal states systems \eqref{eq:analyt_diff}, \eqref{eq:RDPDE_internalstates1}--\eqref{eq:RDPDE_internalstates2} and \eqref{eq:RDPDE2_internalstates1}--\eqref{eq:RDPDE2_internalstates3} are discretized using second order centered finite differences as in \eqref{eq:Laplaceapprox}:
\begin{equation}
\Delta \fatu \equiv \frac{\partial^2 \fatu}{\partial x^2} = \frac{1}{h^2}\,\left( \fatu_{j+1} -2\,\fatu_j + \fatu_{j-1} \right).
\label{eq:discrete_laplacian}
\end{equation}
The jump coefficients in \eqref{eq:Laplaceapprox} are $\lambda_{j1}=\lambda_{j2}=1/h^2$ and $\lambda_j=2/h^2$.
The resulting system of ODEs in time can be written
\begin{equation}
\frac{\partial \fatu_j}{\partial t} = F_{\ell}\left(\fatu_{j-1},\fatu_j,\fatu_{j+1}\right) + F_{n\ell}\left(\fatu_j\right),
\label{eq:crank1}
\end{equation}
where $F_{\ell}$ contains the discrete Laplacian and the linear reaction terms and the $F_{n\ell}$ the nonlinear reaction terms. The system of ODEs \eqref{eq:crank1} is discretized using the following finite difference scheme
\begin{equation}
\frac{\fatu_j^{n+1} - \fatu_j^n}{\Delta t} = \frac{1}{2}\left( F_{\ell}\left(\fatu_{j-1}^{n+1},\fatu_j^{n+1},\fatu_{j+1}^{n+1}\right) + F_{\ell}\left(\fatu_{j-1}^n,\fatu_j^n,\fatu_{j+1}^n\right) \right) + F_{n\ell}\left(\fatu_j^n\right).
\label{eq:crank2}
\end{equation}
In the case of an annihilation process and a monomolecular reaction,
there is no nonlinear reaction term. Hence, $F_{n\ell}$ is zero, and
\eqref{eq:crank2} reduces to the Crank-Nicholson scheme. It is
second-order accurate in space and time and it is unconditionally
stable. In the case of a bimolecular reaction, because of the
nonlinearities, the method is first order accurate in time.

\subsection{Configurations}
\label{sec:exp:config}

In order to demonstrate the ability of the present method to be applied to a wide range of problems, we numerically test in Section~\ref{sec:without_reaction} two different sets of parameters, given in Table~\ref{table:parameters}. Then only the first set of parameters is used in Sections~\ref{sec:reaction_decay}, \ref{sec:linear_reaction} and \ref{sec:nonlinear_reaction}. In our code, $\alpha$, $K_{\alpha}$, $t_{min}$, $t_{max}$ are input parameters, and $\sigma^2$, $\tau$, $N$, $\tau_i$, $\mu_i$ are ouput parameters. 
The quadrature coefficients $\tau_i$ and
$\mu_i$ in \eqref{eq:waitapprox} are determined by nonlinear optimization \cite{BLANC15} and the
corresponding model error $\varepsilon_{mod}$
\eqref{eq:error_model} is also given in the table.

The computational domain is $\Omega=[-1,1]$ in Figure~\ref{fig:test1_u} and
$[-10,10]$ in Figure~\ref{fig:test1_msd}, is discretized with $N_x =
128$ grid points. Neumann boundary conditions are used. As the initial
condition, we use a Gaussian $g(x)$, centered at point
$(0\,,0)$ and of variance $\sigma_g^2 = 10^{-3}$,
rather than a Dirac distribution to avoid spurious numerical
artifacts. Moreover, each internal state is initialized by the weight
$\mu_{\ell}$.\\
%
\begin{table}[htbp]
\begin{center}\footnotesize
\begin{tabular}{llll}
 & Parameters & Set $1$ & Set $2$\\
\hline
\rule[-1mm]{0mm}{3mm} Physical parameters & $\alpha$  & $0.5$ & $0.5$\\
\rule[-1mm]{0mm}{3mm}  & $K_{\alpha}$ (m$^2\,$s$^{-\alpha}$) & $0.04$ & $0.04$\\
\rule[-1mm]{0mm}{3mm}  & $t_{min}$ (s) & $10^{-4}$ & $10^{-3}$\\
\rule[-1.5mm]{0mm}{3mm}  & $t_{max}$ (s) & $5\cdot10^{-2}$ & $1$\\
\rule[-1mm]{0mm}{3mm}  & $\sigma^2$ (m$^2$) & $3.49\cdot10^{-4}$ & $7.18\cdot10^{-4}$\\
\rule[-1mm]{0mm}{3mm}  & $\tau$ (s) & $7.62\cdot10^{-5}$ & $3.22\cdot10^{-4}$\\
\rule[-1mm]{0mm}{3mm} Optimization & $N$ & $4$ & $5$\\

\rule[-1mm]{0mm}{3mm} & $\tau_1$ (s) & $9.51\cdot10^{-5}$ & $7.58\cdot10^{-4}$\\
\rule[-1mm]{0mm}{3mm} & $\tau_2$ (s) & $5.40\cdot10^{-4}$ & $3.55\cdot10^{-3}$\\
\rule[-1mm]{0mm}{3mm} & $\tau_3$ (s) & $3.09\cdot10^{-3}$ & $1.66\cdot10^{-2}$\\
\rule[-1mm]{0mm}{3mm} & $\tau_4$ (s) & $2.13\cdot10^{-2}$ & $7.89\cdot10^{-2}$\\
\rule[-1mm]{0mm}{3mm} & $\tau_5$ (s) & - & $4.78\cdot10^{-1}$\\

\rule[-1mm]{0mm}{3mm} & $\mu_1$ & $4.96\cdot10^{-1}$ & $3.23\cdot10^{-1}$\\
\rule[-1mm]{0mm}{3mm} & $\mu_2$ & $2.07\cdot10^{-1}$ & $1.48\cdot10^{-1}$\\
\rule[-1mm]{0mm}{3mm} & $\mu_3$ & $8.80\cdot10^{-2}$ & $6.84\cdot10^{-2}$\\
\rule[-1mm]{0mm}{3mm} & $\mu_4$ & $4.42\cdot10^{-2}$ & $3.24\cdot10^{-2}$\\
\rule[-1mm]{0mm}{3mm} & $\mu_5$ & - & $1.85\cdot10^{-2}$\\
\rule[-1mm]{0mm}{3mm} & $\varepsilon_{mod}$ & $5.25\cdot10^{-2}$ & $2.92\cdot10^{-2}$\\ \hline
\end{tabular}
\end{center}
\caption{Parameters used in the numerical experiments.}
\label{table:parameters}
\end{table}


\subsection{Subdiffusion}\label{sec:without_reaction}

\begin{figure}[htbp]
\begin{center}
\begin{tabular}{cc}
(a) & (b)\\
\includegraphics[scale=0.35]{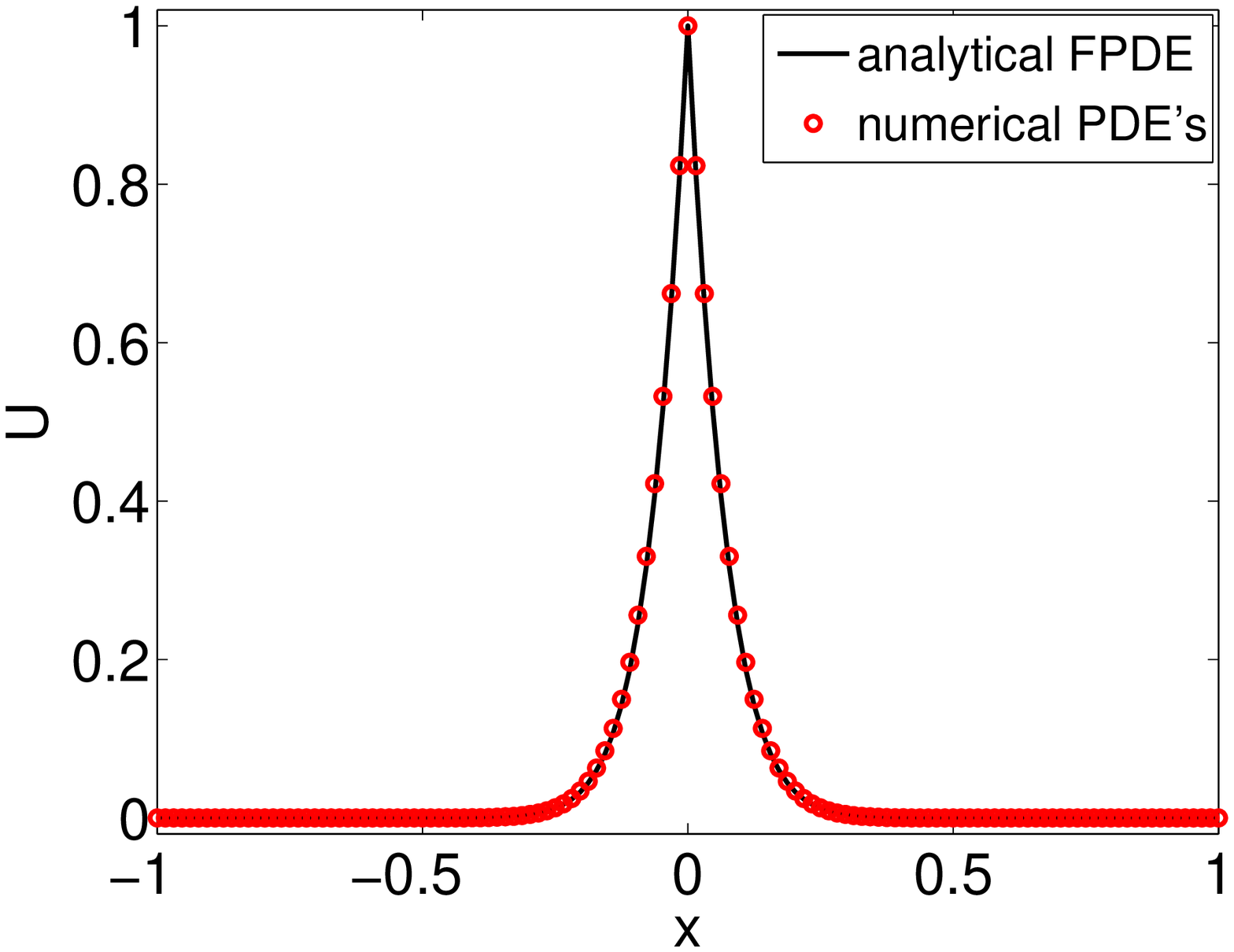} &
\includegraphics[scale=0.35]{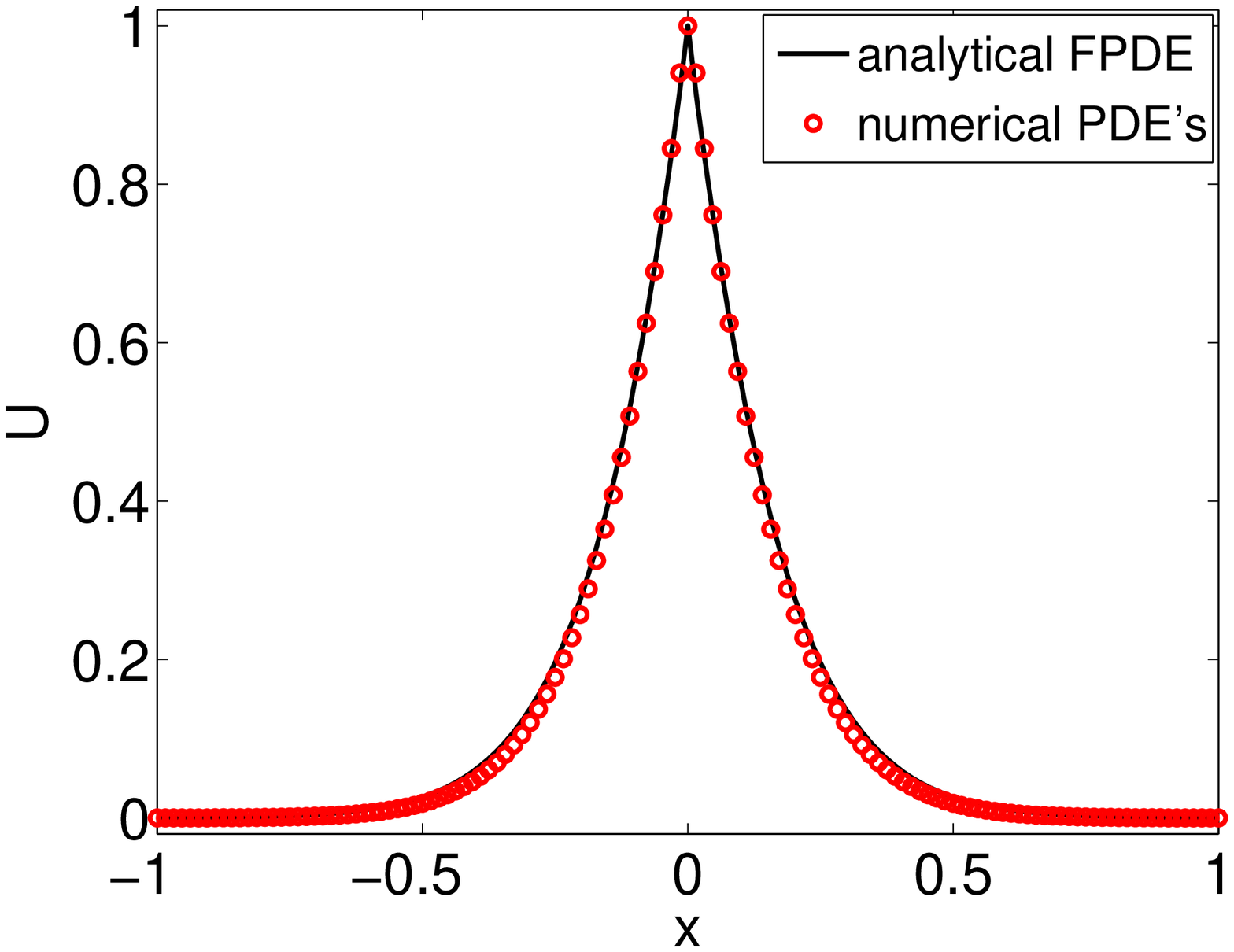}
\end{tabular}
\end{center}
\caption{Section~\ref{sec:without_reaction}. Comparison between the numerical values (circle) and the analytical values (solid line) of the concentration $U=\fate^T\,\fatu$ of $A$. (a): Set 1 of parameters at $t_1 = 5\cdot10^{-3}$ s, (b): Set 2 of parameters at $t_2 = 1.5\cdot 10^{-1}$ s.}
\label{fig:test1_u}
\end{figure}

The aim of the first test is to check the accuracy of the numerical method above when no reaction occurs. The following initial conditions are used
\begin{equation}
\fatu(x,0) = \fatmu\,g(x).
\label{eq:initial_condition_test1}
\end{equation}
Figure \ref{fig:test1_u} compares the numerical solution $U=\fate^T\,\fatu$ obtained with the PDEs in the internal state diffusion system \eqref{eq:analyt_diff} with the analytical solution of the FPDE \eqref{eq:sol_FPDE1_meijer}. Figure~\ref{fig:test1_u}-(a) corresponds to the first set of parameters, given in Table~\ref{table:parameters}, at time $t_1 = 5\cdot10^{-3}$ s and Figure~\ref{fig:test1_u}-(b) corresponds to the second set of parameters at time $t_2 = 1.5\cdot10^{-1}$ s. Note in Table~\ref{table:parameters} that the macroscopic parameters $\alpha$ and $K_{\alpha}$, involved in \eqref{eq:FPDE}, are the same for both sets of parameters. Consequently, they approximate the same macroscopic FPDE, but in two different time intervals $\left[ t_{min},t_{max}\right]$. In both cases, the final time of the simulation in Figure~\ref{fig:test1_u} is inside this interval. Excellent agreement is found between the two solutions.

Two errors should be mentioned here: the modeling error, defined as the difference between the internal states diffusion model and the FPDE model, and the numerical error $\varepsilon_{num}$, resulting from the numerical discretization of the internal states diffusion system. The modeling error is dependent on $\varepsilon_{mod}$ in \eqref{eq:error_model}. The error $\varepsilon_{mod}$ is given in Table~\ref{table:parameters} and the total error $\varepsilon_{tot}$ is measured in Figure~\ref{fig:test1_u} as the difference between the numerical solution obtained with the PDE model (red circles) with the analytical solution obtained with the FPDE model (black solid line). In Figure~\ref{fig:test1_u}-(a) $\varepsilon_{tot} \simeq 2.33\cdot10^{-2}$ and in Figure~\ref{fig:test1_u}-(b) $\varepsilon_{tot} \simeq 3.65\cdot10^{-2}$.

\begin{figure}[htbp]
\begin{center}
\begin{tabular}{cc}
(a) & (b)\\
\includegraphics[scale=0.30]{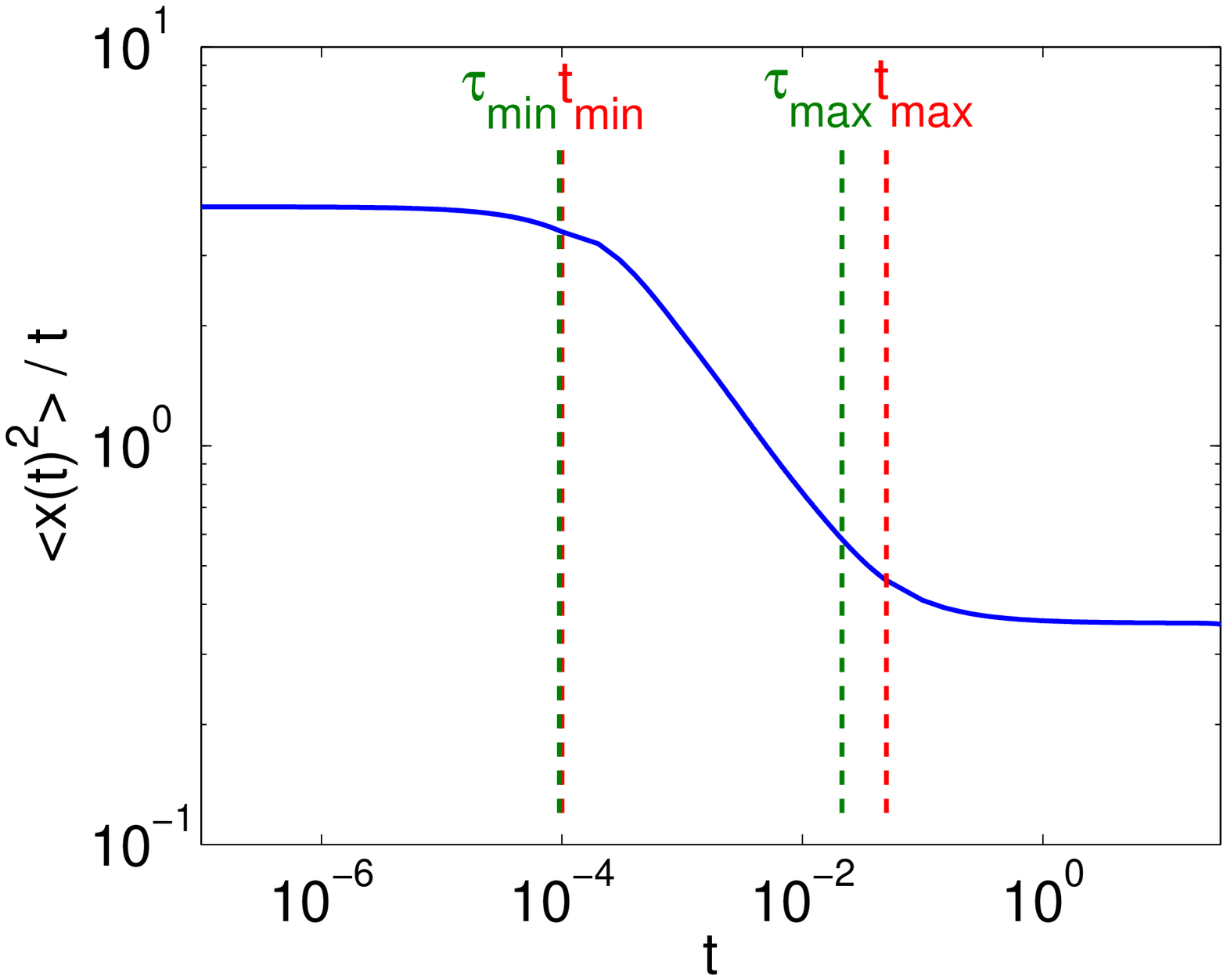} & 
\includegraphics[scale=0.30]{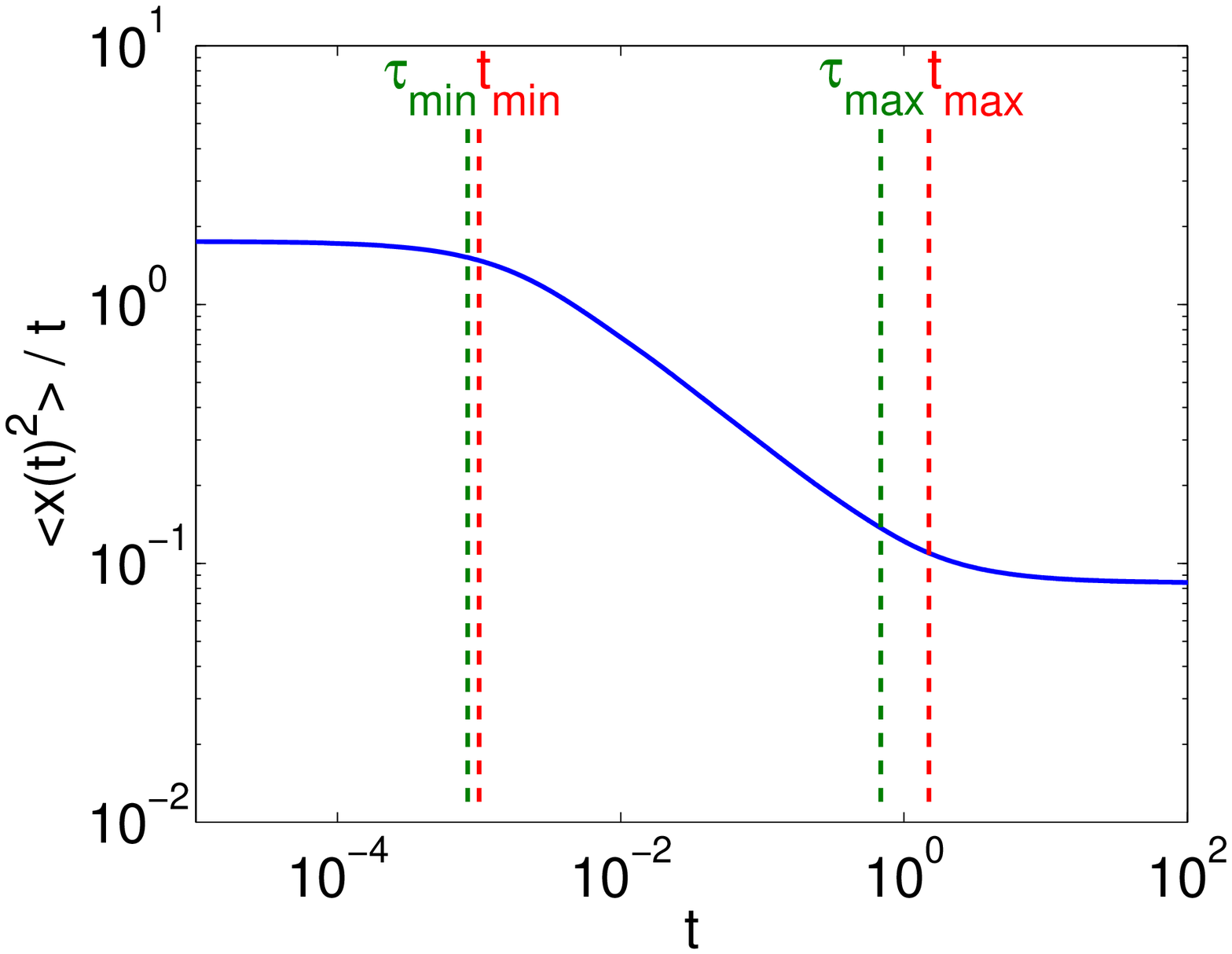}
\end{tabular}
\end{center}
\caption{Section~\ref{sec:without_reaction}. Mean square displacement
  as a function of time. (a): Parameters set 1, (b): Parameters set 2.
  The scales are logarithmic on both axes.}
\label{fig:test1_msd}
\end{figure}

Figure \ref{fig:test1_msd} shows the mean square displacement
divided by the time, which is constant for
ordinary diffusion, see \eqref{eq:MSD_brownian}. According to the analysis in
Section~\ref{subsec:asymptotic}, ordinary diffusion is observed for $t \ll
\min\limits_{\ell=1,...,N}\tau_{\ell}$ and for $t \gg
\max\limits_{\ell=1,...,N}\tau_{\ell}$. 
For intermediate times, the diffusion is anomalous by construction, cf. Section~\ref{subsec:internal_states}. 
The subdiffusive exponent
$\alpha$ is measured by linear regression: $\alpha \simeq 0.4997$ in Figure~\ref{fig:test1_msd}-(a) and 
$\alpha \simeq 0.5454$ in Figure~\ref{fig:test1_msd}-(b). In both cases
the measured subdiffusive exponent is close to the theoretical one.\\


\subsection{Annihilation process}\label{sec:reaction_decay}

\begin{figure}[htbp]
\begin{center}
\begin{tabular}{cc}
Model I & Model II\\
\includegraphics[scale=0.35]{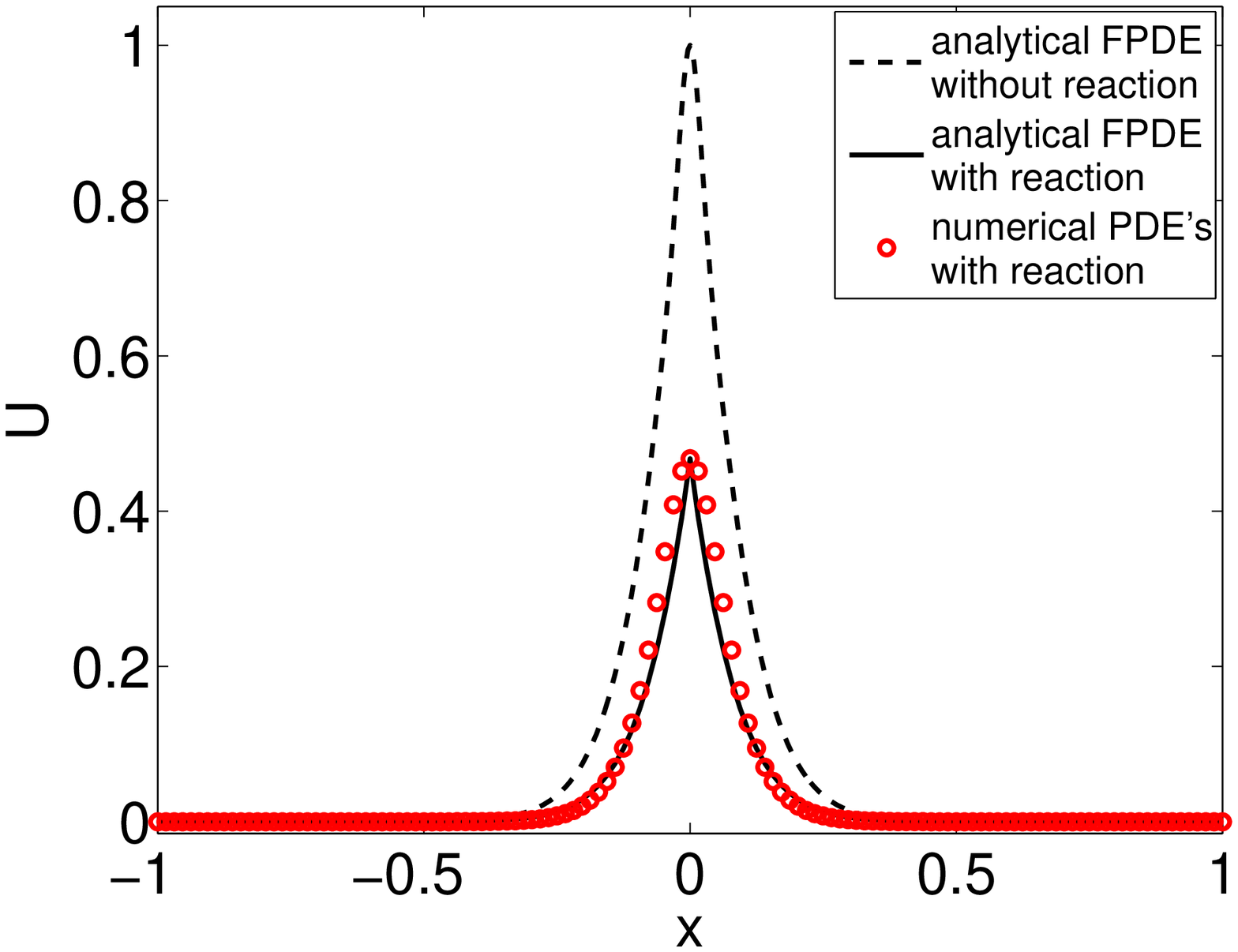} &
\includegraphics[scale=0.35]{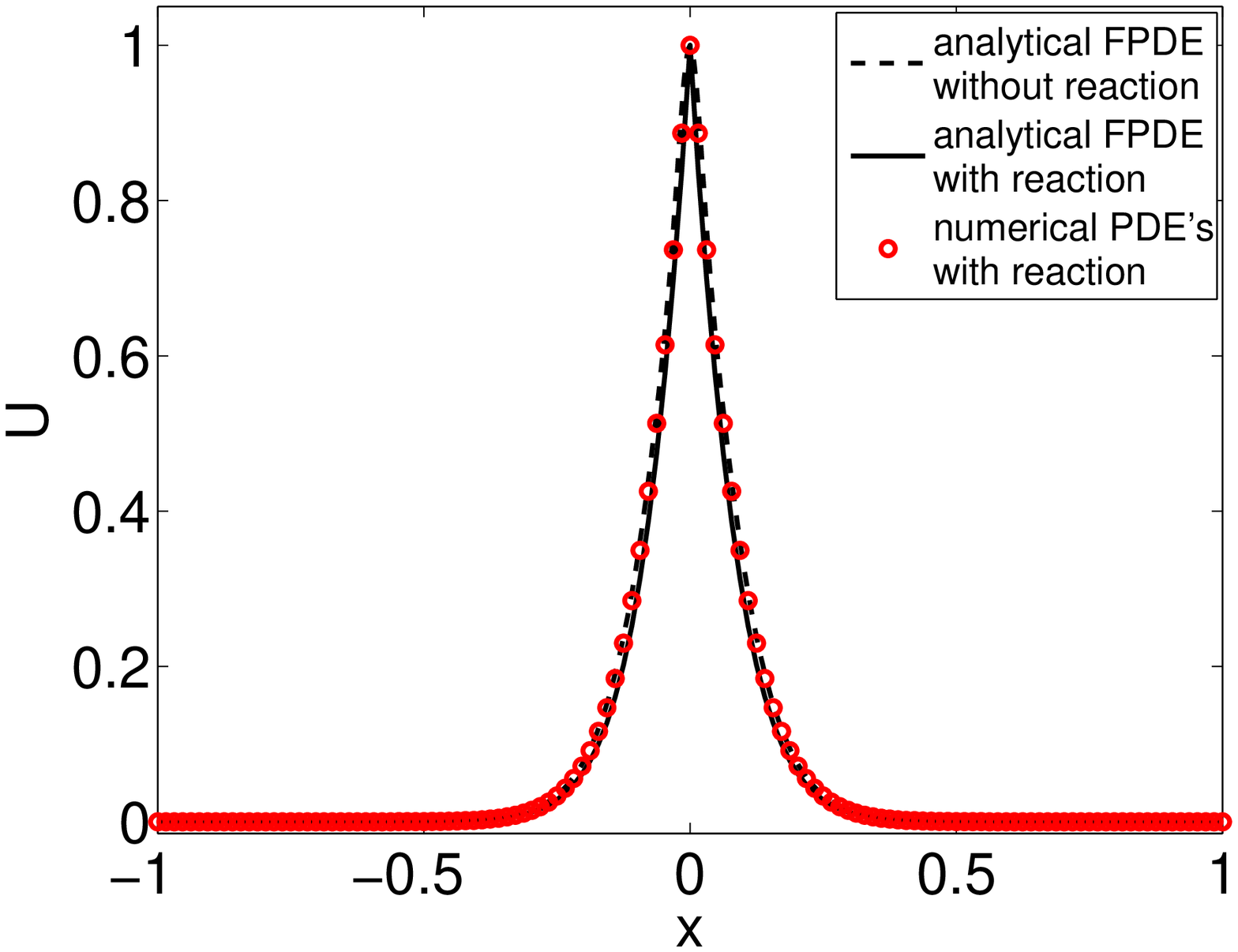}\\
\includegraphics[scale=0.35]{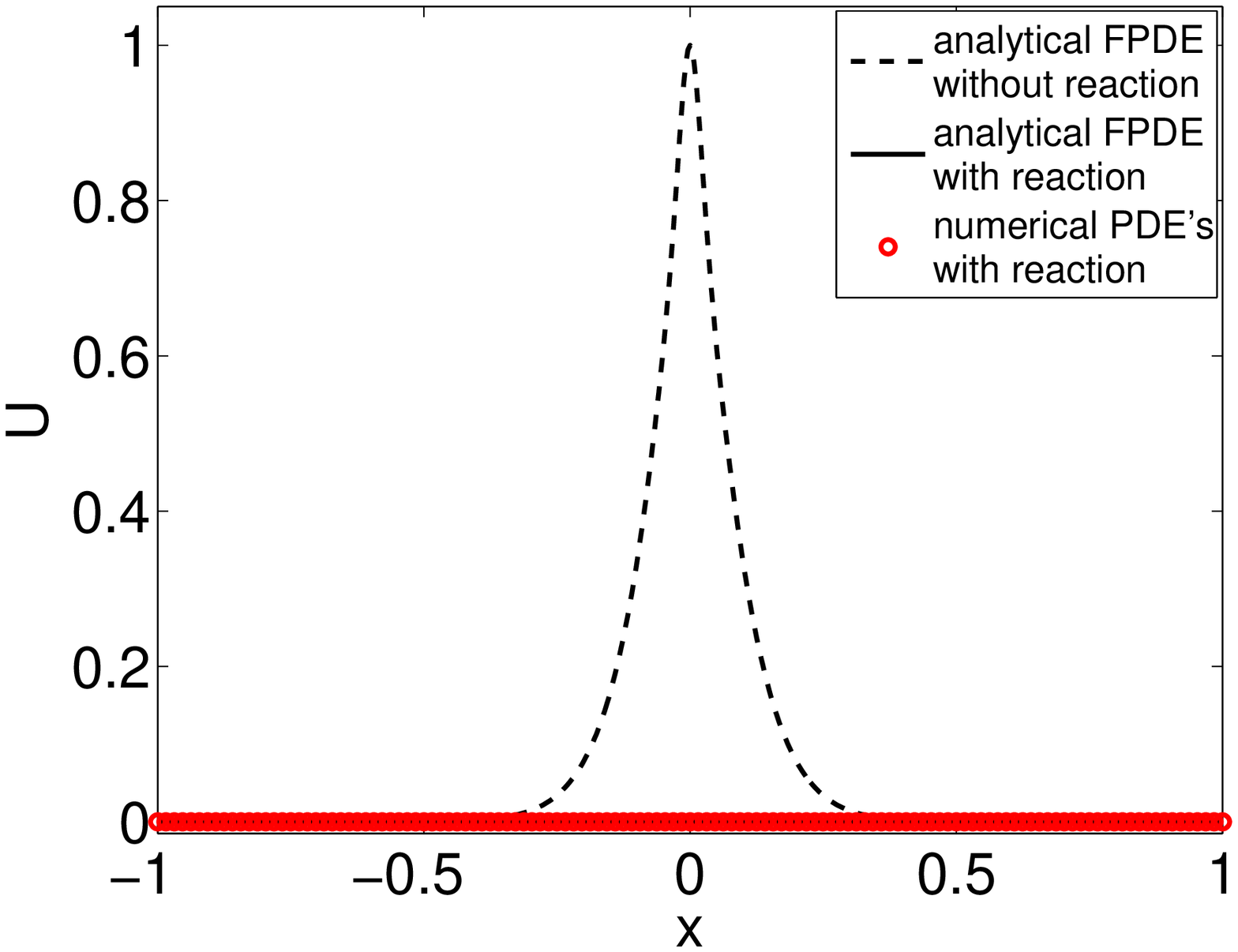} &
\includegraphics[scale=0.35]{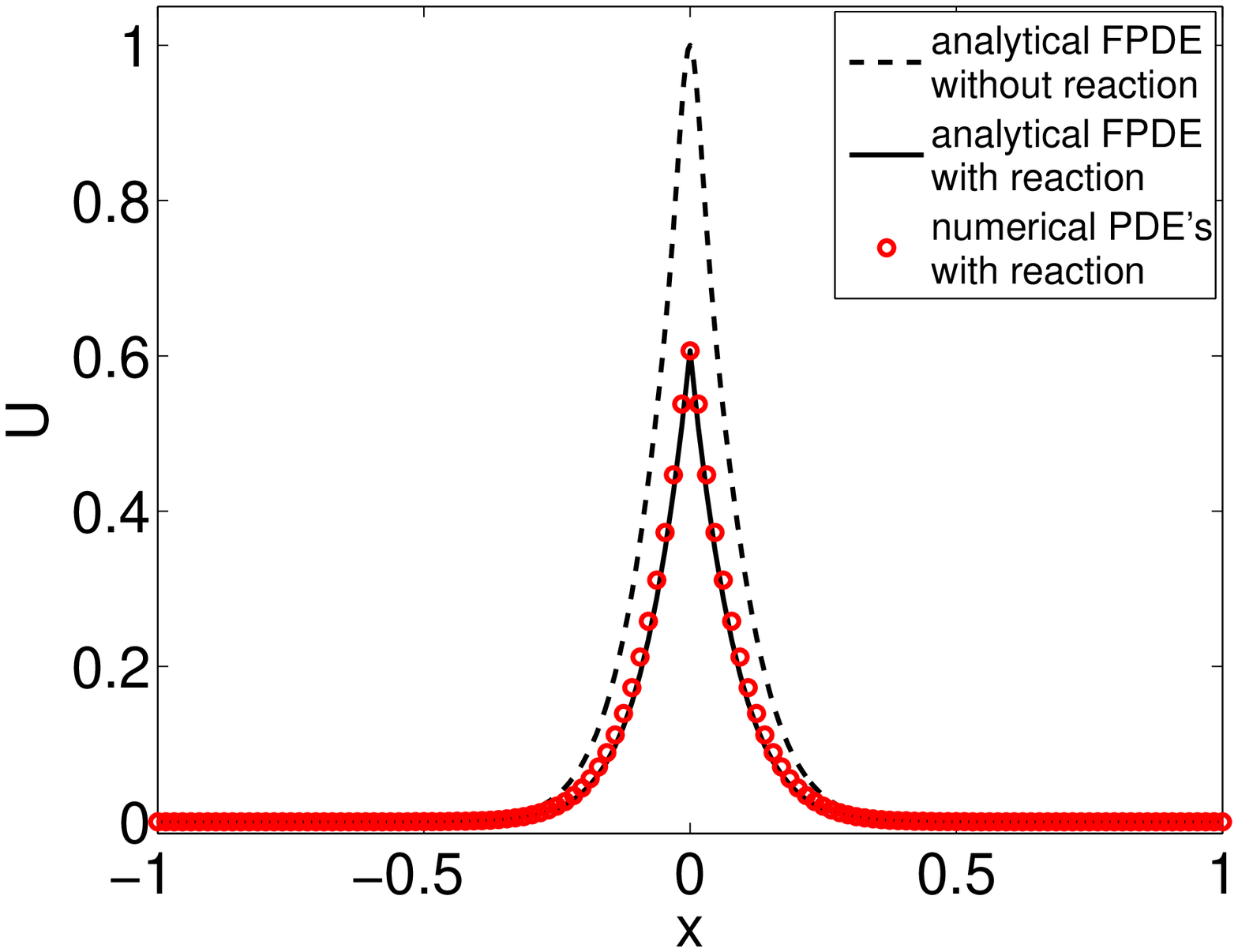}
\end{tabular}
\end{center}
\caption{Section~\ref{sec:reaction_decay} Comparison between the numerical values (circles) and the analytical values (black line) of the conecntration $U=\fate^T\,\fatu$ of $A$ at time $t = 10^{-2}$ s. Top: First run, bottom: Second run.}
\label{fig:test2_u}
\end{figure}

\begin{table}[htbp]
\begin{center}\footnotesize
\begin{tabular}{l|cc||cc}
& \multicolumn{2}{c||}{First run} & \multicolumn{2}{c}{Second run}\\
[5pt]
\rule[-1.5mm]{0mm}{3mm} & Model I & Model II & Model I & Model II\\
\hline
\rule[-1mm]{0mm}{4mm} $k$ & $0.1$ & $8.727\cdot10^{-4}$ s$^{-1}$ & $5.729\cdot10^{-3}$ & $50$ s$^{-1}$\\
\rule[-1mm]{0mm}{3mm} $k_\ast$ & $8.727\cdot10^{-4}$ s$^{-\alpha}$ & $8.727\cdot10^{-4}$ s$^{-1}$ & $50$ s$^{-\alpha}$ & $50$ s$^{-1}$
\end{tabular}
\end{center}
\caption{Reaction rates of annihilation process.}
\label{table:parameters_annihilation}
\end{table}
The purpose of the second test is to investigate the accuracy of the numerical method in the case of an annihilation process. 
Two different models are possible in Section~\ref{subsubsec:reaction_modeling_annihilation}
\begin{equation}
\begin{array}{ll}
\mbox{Model I}\quad & A_i \mathop{\longrightarrow}\limits^{k/\tau_i} \emptyset,\\
[10pt]
\mbox{Model II}\quad & A_i \mathop{\longrightarrow}\limits^{k} \emptyset.
\end{array}
\end{equation}
The reaction coefficient $k$ used in the simulations and the
corresponding macroscopic reaction rate $k_\ast$
\eqref{eq:FPDEmesocoeff}--\eqref{eq:FPDEmesocoeff2} are given in
Table~\ref{table:parameters_annihilation}. The rate $k$ has been
chosen such that the two models have the same macroscopic reaction
rates in \eqref{eq:FPDEmesocoeff} and \eqref{eq:FPDEmesocoeff2}. Note
the order-of-magnitude difference between the value of $k$ in model I
and II, which is explained by the fact that the units of $k$ and
$k_\ast$ are are different. The initial conditions are the same as in
the previous test. Figure~\ref{fig:test2_u} compares the numerical
solution $U=\fate^T\,\fatu$ obtained with the internal state
reaction-diffusion system \eqref{eq:PDE_internalstates_reaction} with
the analytical solutions of the FPDEs in
\eqref{eq:sol_FPDE_decayI_fox} and \eqref{eq:sol_FPDE_decayII_fox} at
time $t_1 = 10^{-2}$ s. For both models, there is excellent agreement
between the two solutions. The analytical solution is also given in
the case where no reaction occurs ($k = 0$). For a small macroscopic
reaction rate (first run, top of Figure \ref{fig:test2_u}), no
reaction takes place yet in the case of model II. For large
macroscopic reaction rate (second run, bottom of Figure
\ref{fig:test2_u}), all particles have disappeared in model I. The
differences between the two models are clearly illustrated in
Figure~\ref{fig:test2_u}.\\

\begin{figure}[htbp]
\begin{center}
\begin{tabular}{cc}
Model I & Model II\\
\includegraphics[scale=0.30]{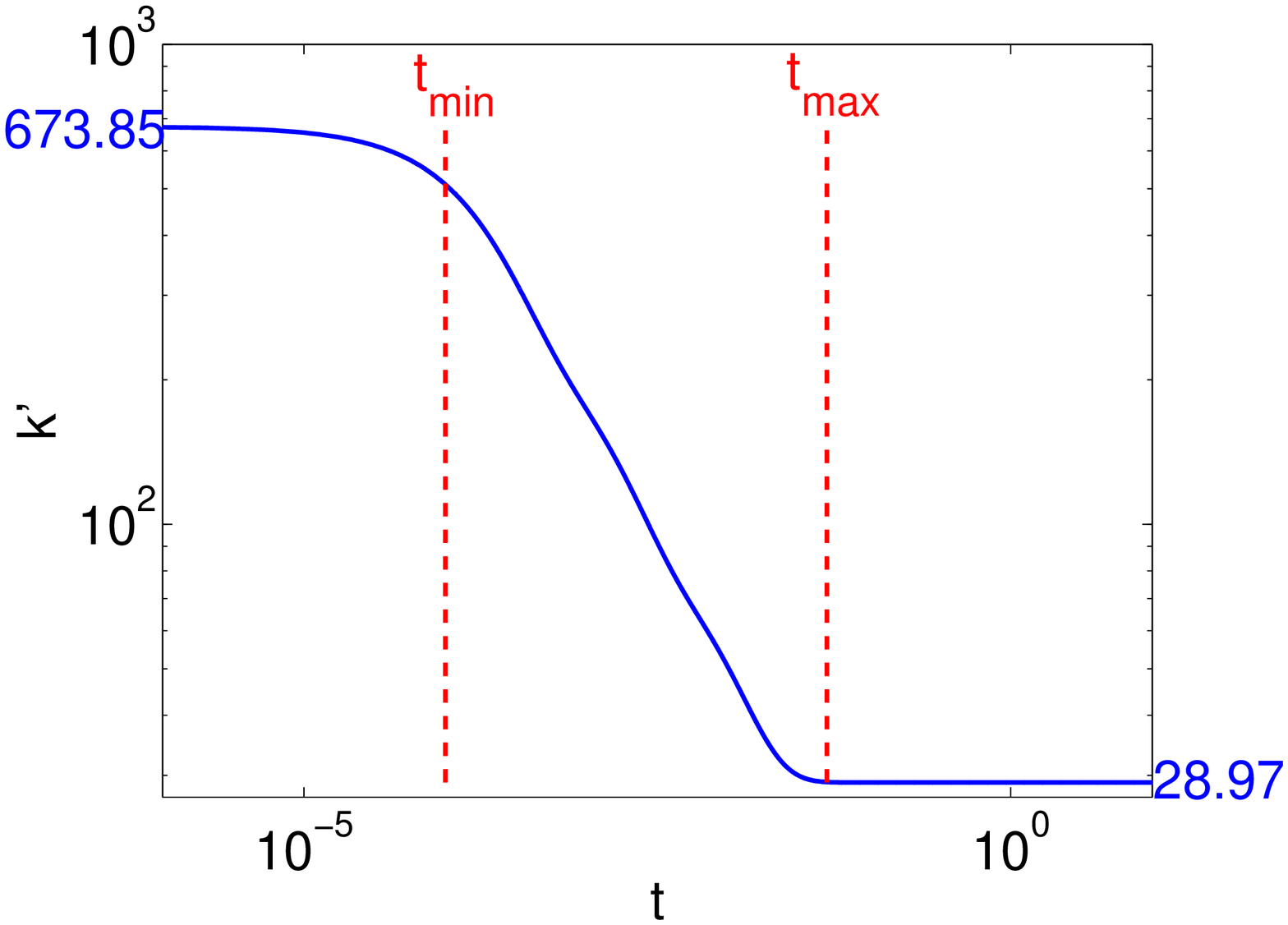} &
\includegraphics[scale=0.30]{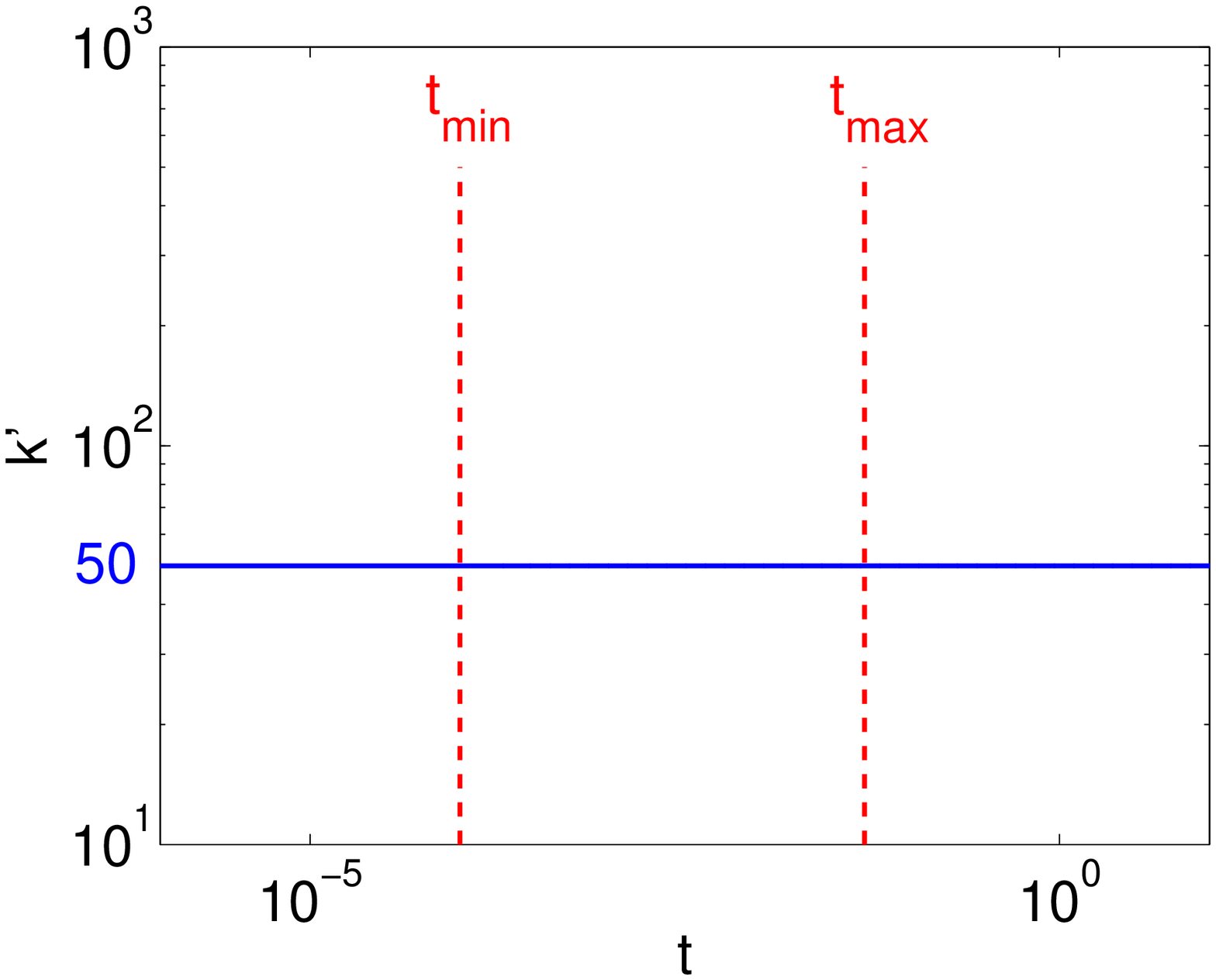}
\end{tabular}
\end{center}
\caption{Section~\ref{sec:reaction_decay} $k'$ in \eqref{eq:kprim} in terms of time. Left: Model I, right: Model II.
The scales are logarithmic on both axes.}
\label{fig:test2_kprim}
\end{figure}

\begin{figure}[htbp]
\begin{center}
\begin{tabular}{ccc}
short time & intermediate time & long time\\
\includegraphics[scale=0.27]{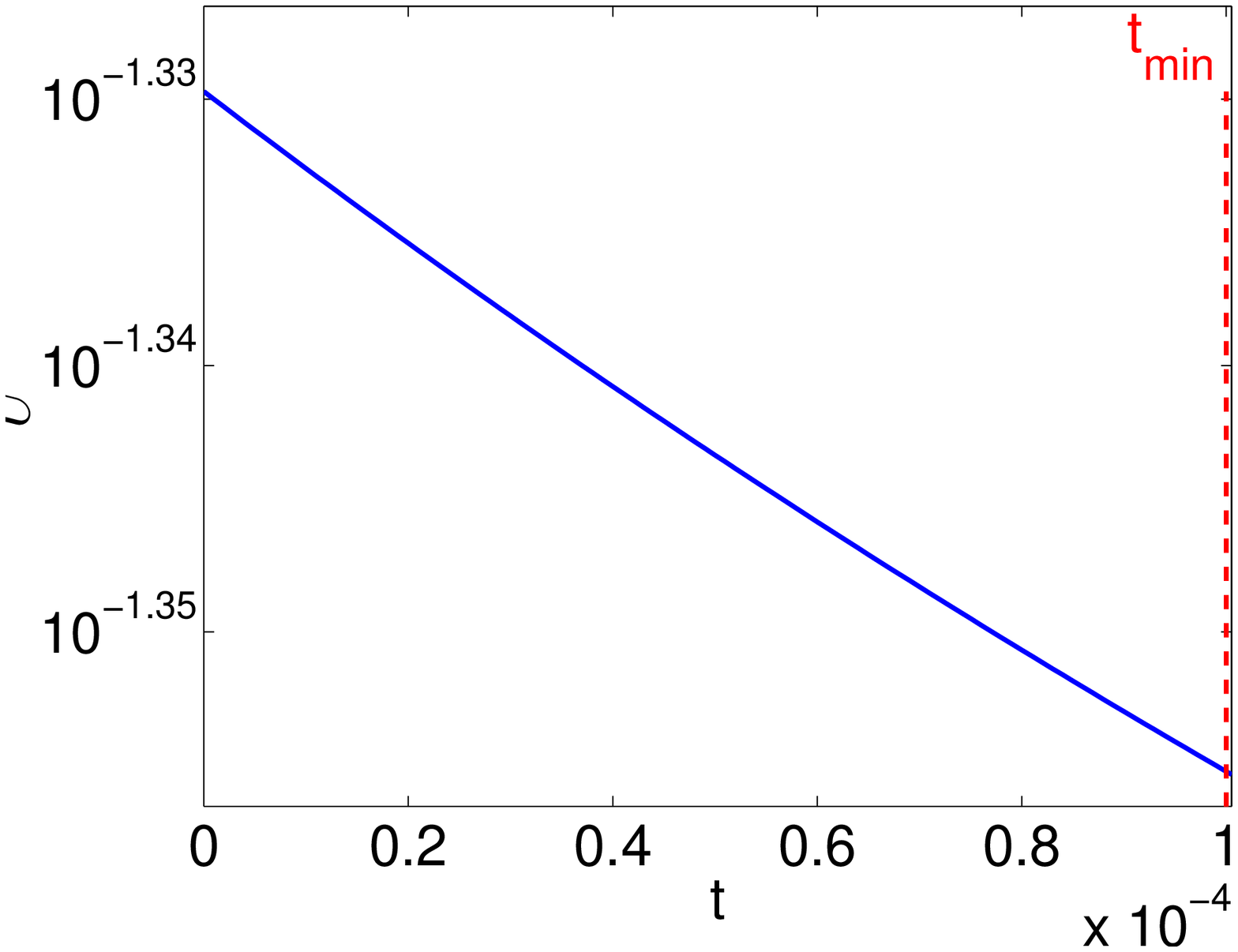} &
\includegraphics[scale=0.27]{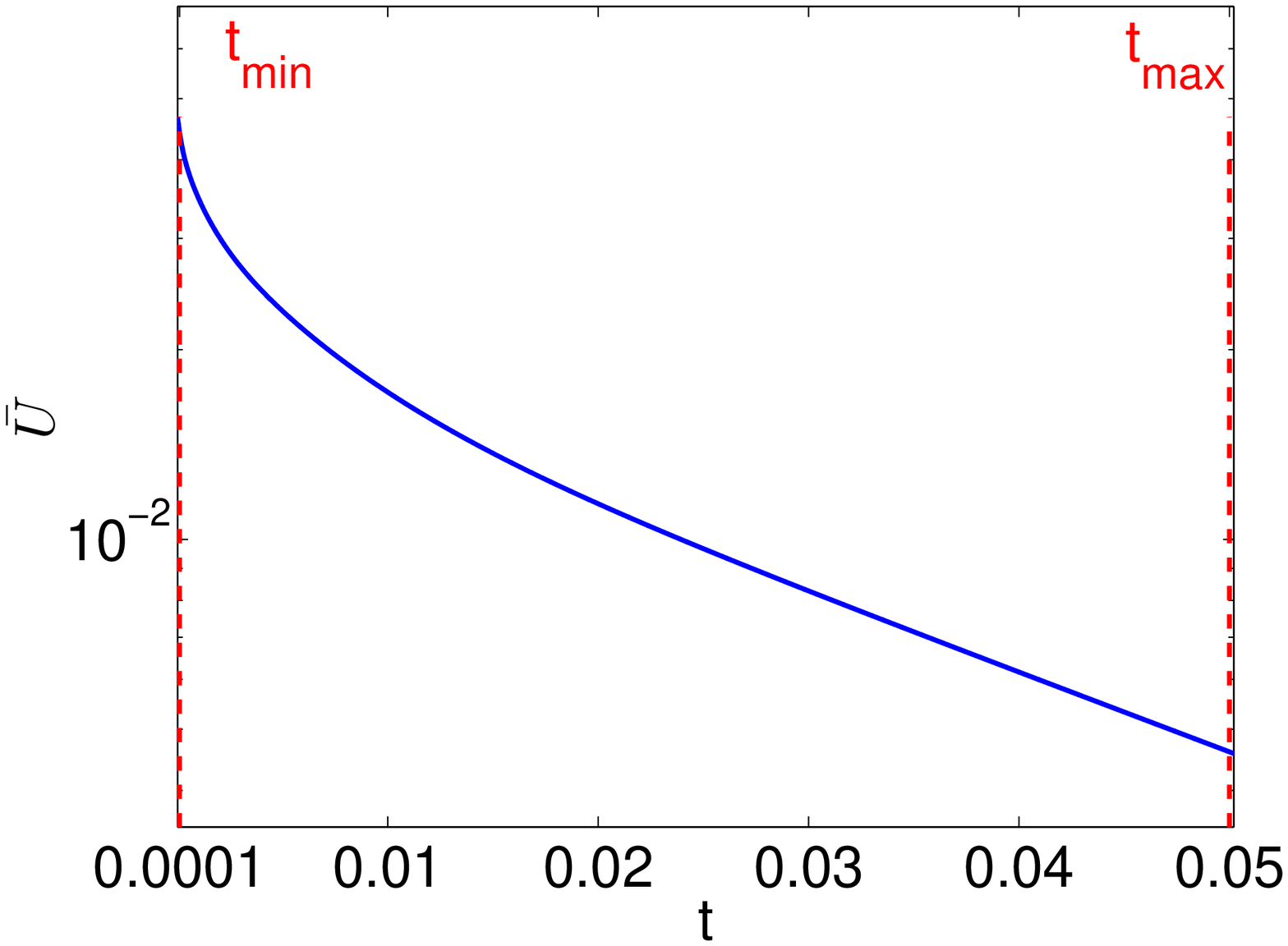} &
\includegraphics[scale=0.27]{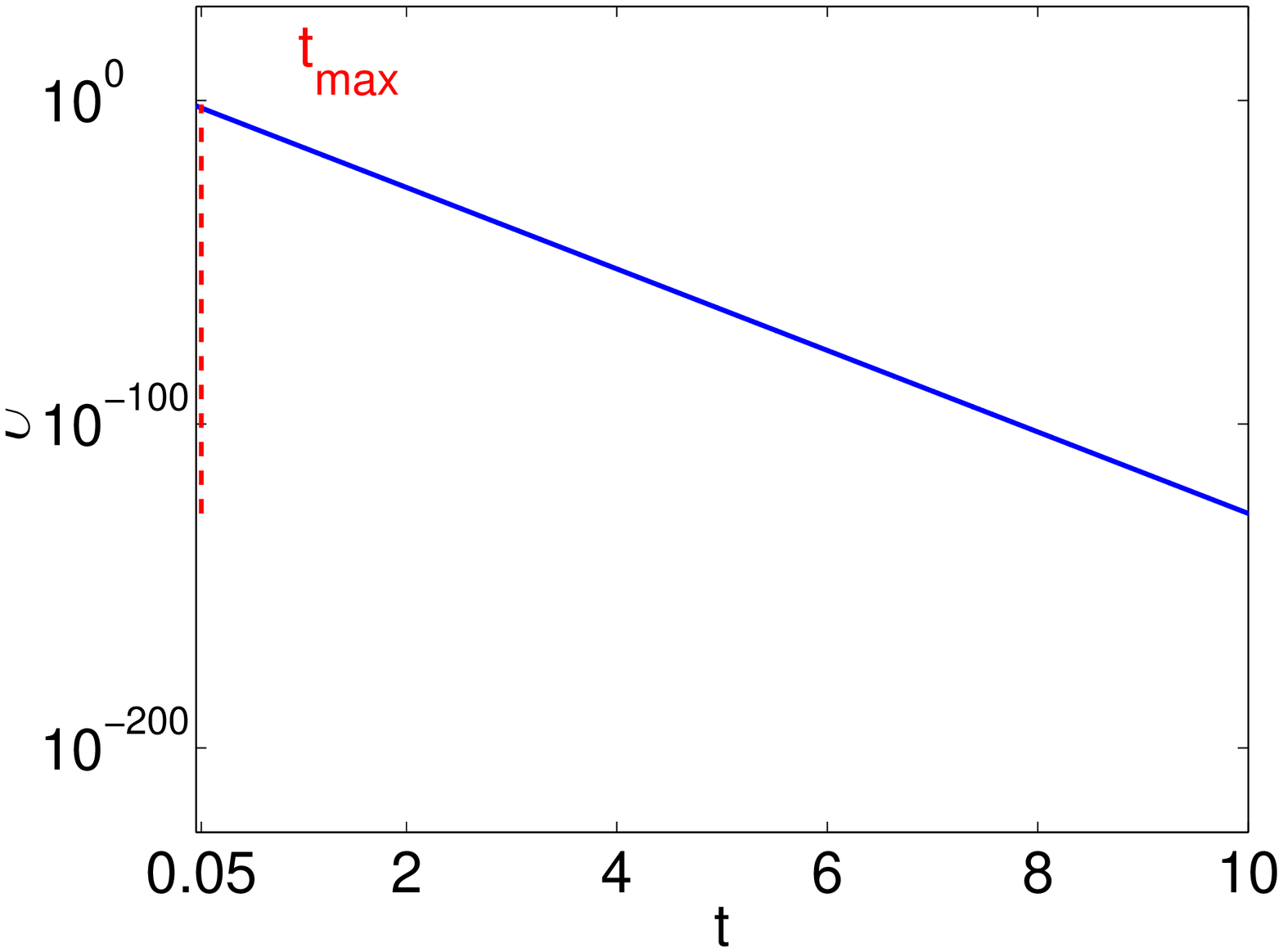}\\
\includegraphics[scale=0.27]{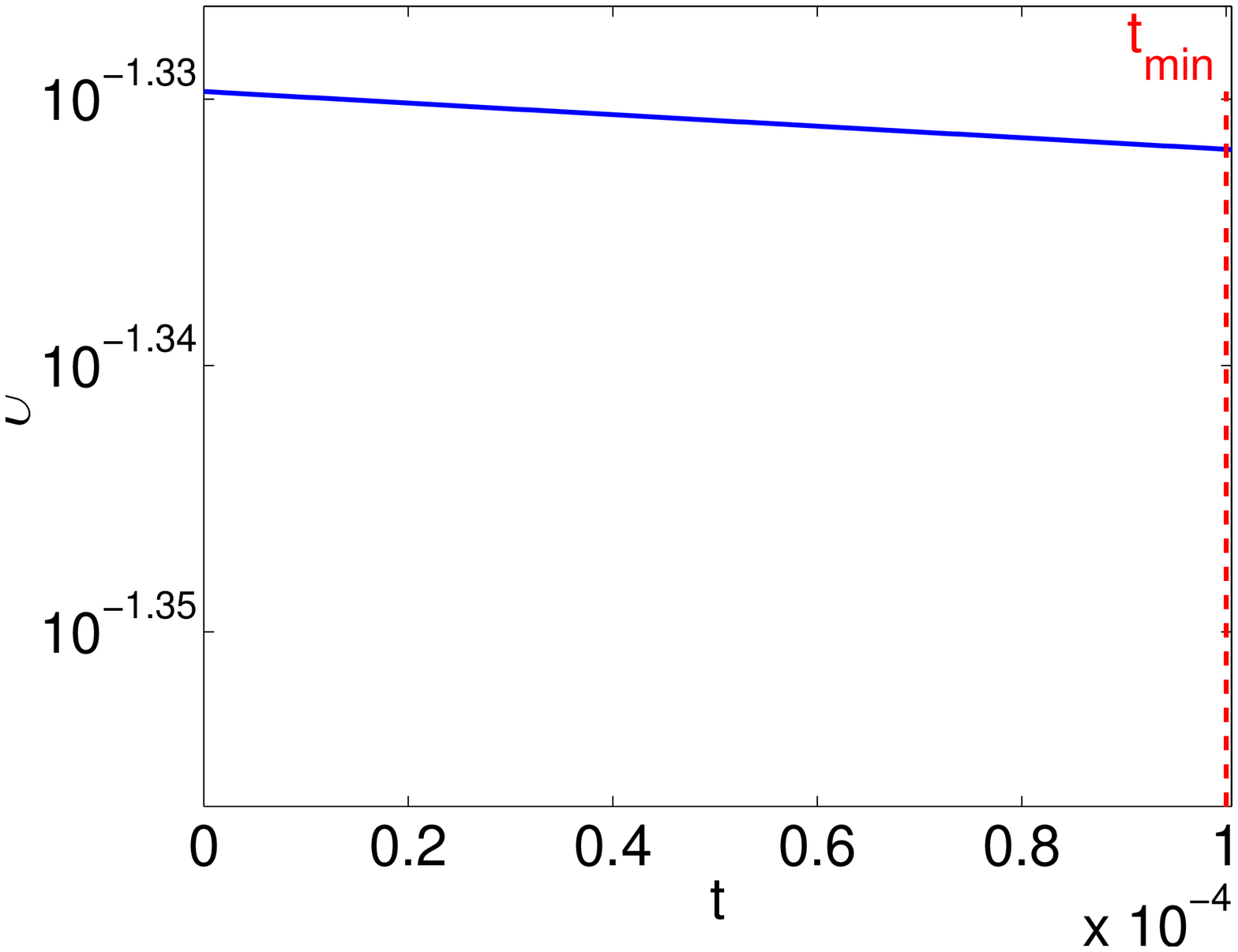} &
\includegraphics[scale=0.27]{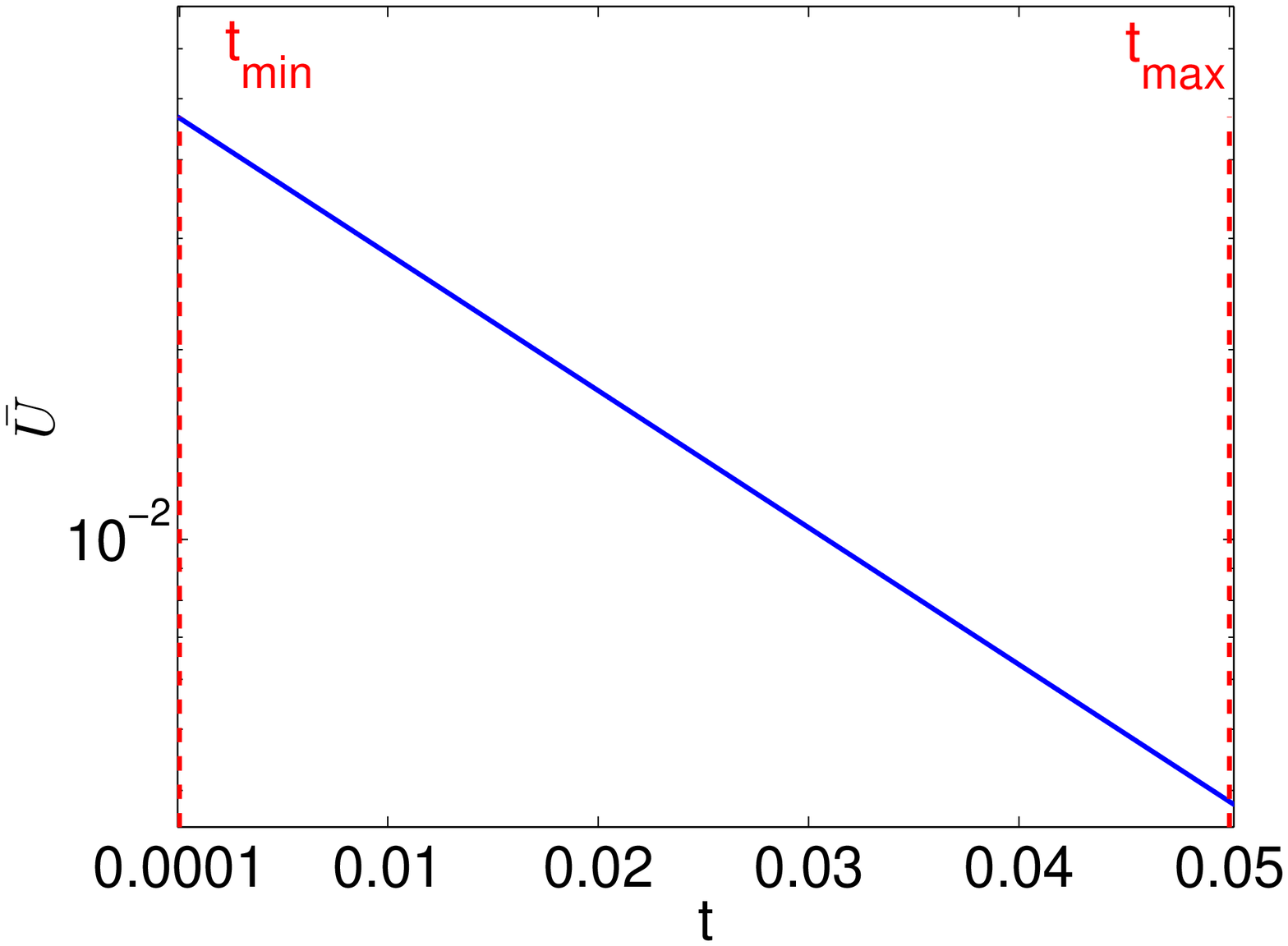} &
\includegraphics[scale=0.27]{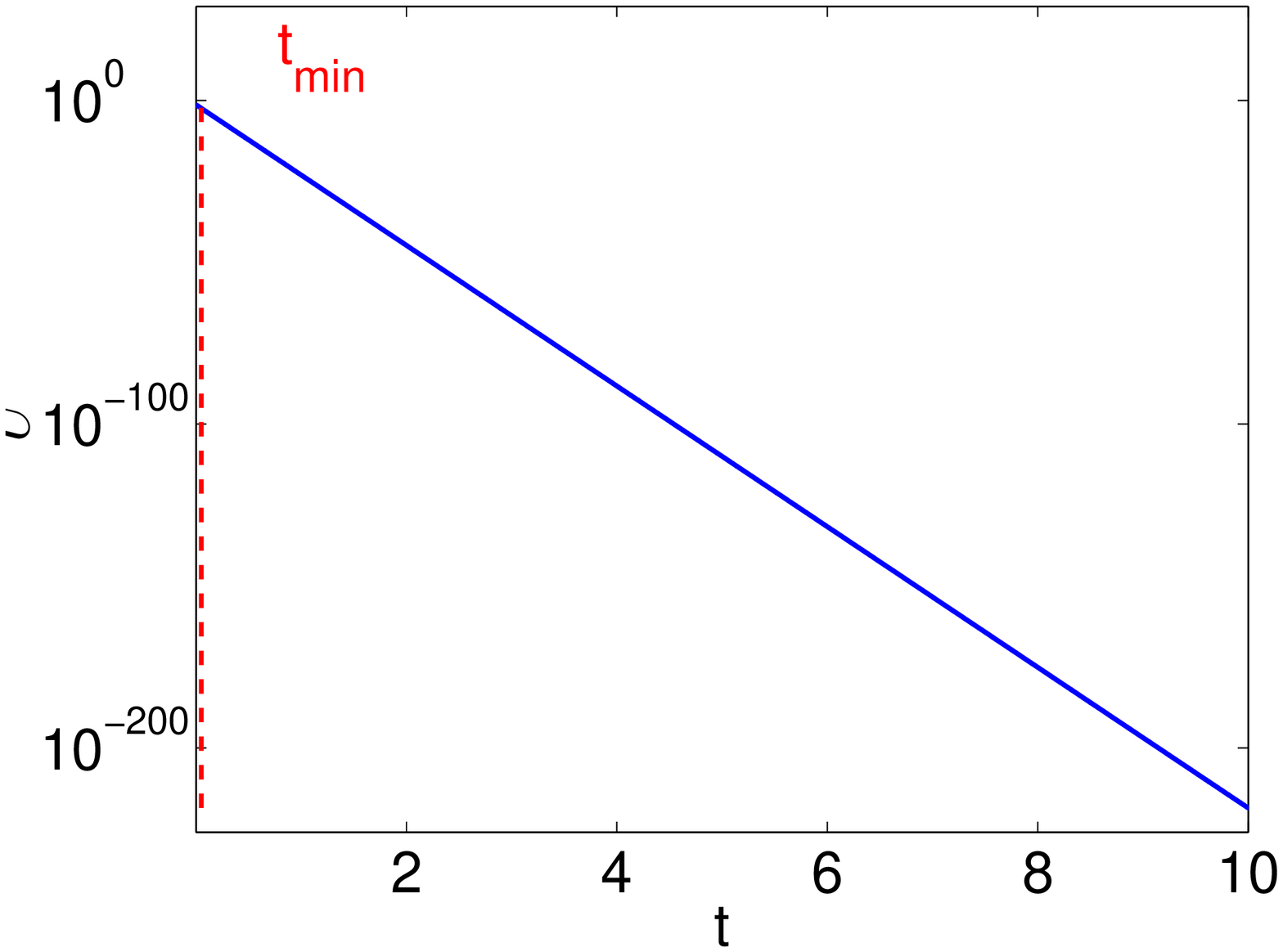}
\end{tabular}
\end{center}
\caption{Section~\ref{sec:reaction_decay}. Total amount of $A$ \eqref{eq:totubareq} in terms of time. Upper: Model I, lower: Model II.
The scale is logarithmic on the $y$ axis.}
\label{fig:test2_Ubar}
\end{figure}

The time development of $k'$ \eqref{eq:kprim} is plotted in Figure~\ref{fig:test2_kprim}. In model I (first run in Table~\ref{table:parameters_annihilation}), $k'$ is almost constant at small $t$ and at large $t$, meaning that the kinetics is ordinary there but at intermediate $t$ $k'$ varies in time implying that the kinetics is anomalous. This is illustrated in Figure~\ref{fig:test2_Ubar}, where the time evolution of the total amount of $A$, i.e. $\bar{U}$ in \eqref{eq:totubareq}, is found. Exponential decay characterizing ordinary kinetics is observed for $t<t_{min}$ and $t>t_{max}$. The exponential parameters are measured by linear regression in Figure~\ref{fig:test2_Ubar} resulting in $k'(0)\approx 28.97$ and $k'_\infty\approx 674.08$ when $t\rightarrow\infty$. In both cases, the measured exponential parameter is close to the theoretical one in \eqref{eq:kprim} (blue values in Figure~\ref{fig:test2_kprim}). For $t_{min}<t<t_{max}$, $\bar{U}$ does not decrease exponentially. In model II (second run in Table~\ref{table:parameters_annihilation}), Figure~\ref{fig:test2_kprim} shows that $k'=k$ does not depend on time and exponential decay is observed for all times in Figure~\ref{fig:test2_Ubar}. The measured exponential parameter $k'\approx 50.00$ is close to the theoretical one. These observations agree with the analysis in Sections~\ref{subsubsec:reaction_modeling_annihilation} and \ref{sec:IS_annihilation}.


\subsection{Monomolecular reaction}\label{sec:linear_reaction}

\begin{figure}[htbp]
\begin{center}
\begin{tabular}{cc}
Model I & Model II\\
\includegraphics[scale=0.35]{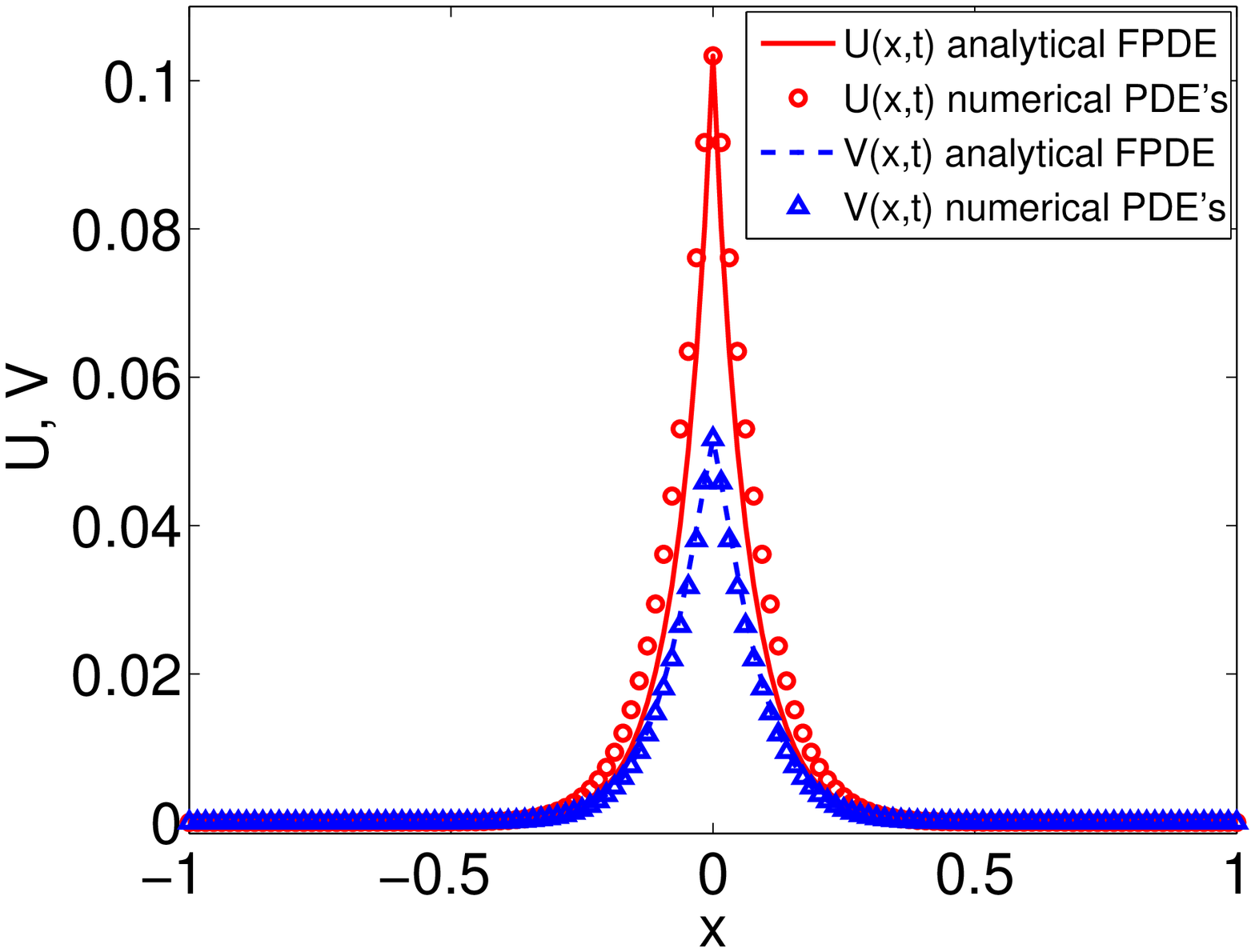} &
\includegraphics[scale=0.35]{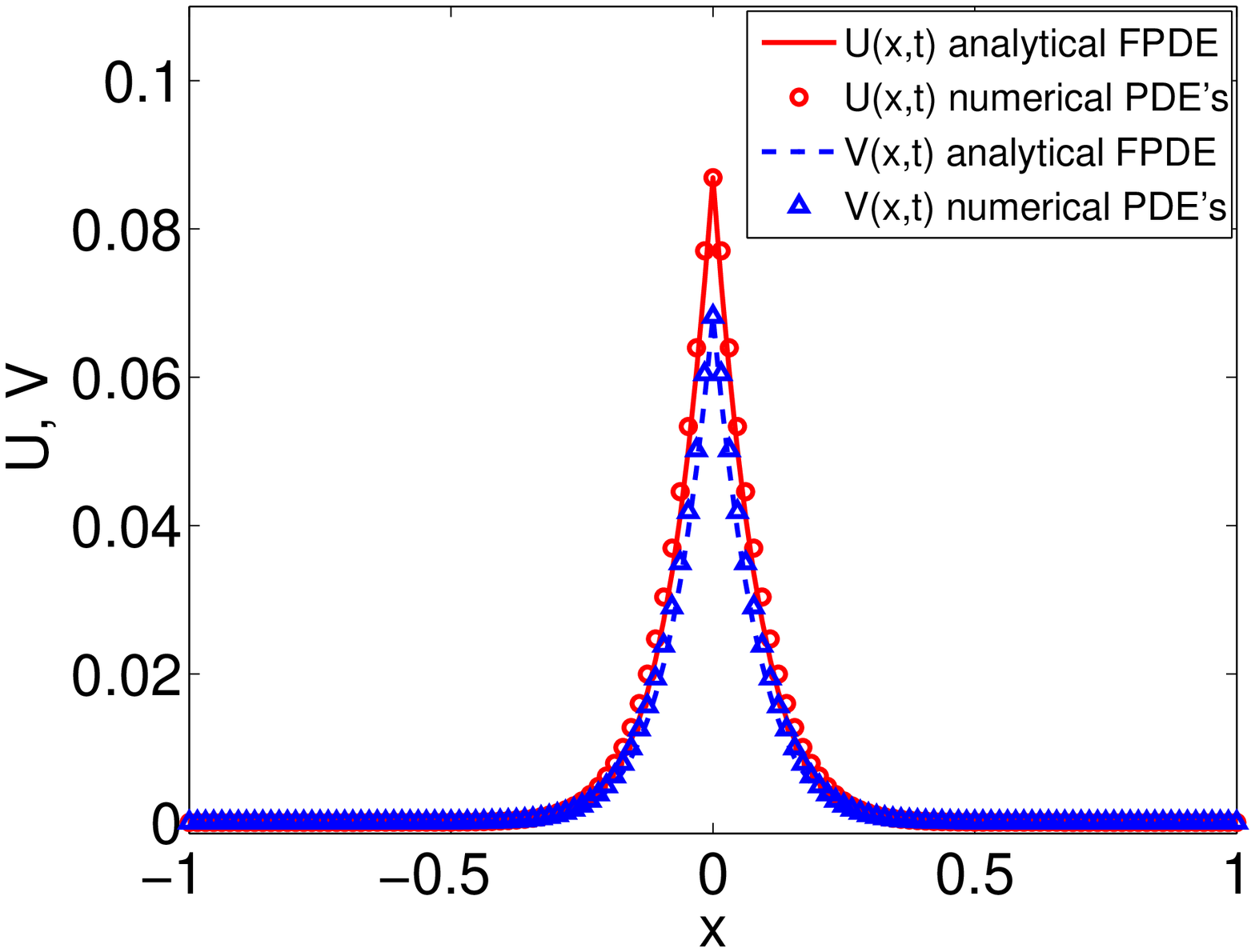}
\end{tabular}
\end{center}
\caption{Section~\ref{sec:linear_reaction}. Comparison between the numerical values (circles and triangles) and the analytical values (solid line and dashed line) of the concentration $U=\fate^T\,\fatu$, $V=\fate^T\,\fatv$ of $A$ and $B$ at time $t_1 = 10^{-2}$ s.}
\label{fig:test3_u}
\end{figure}

\begin{table}[htbp]
\begin{center}\footnotesize
\begin{tabular}{l|cc}
\rule[-1.5mm]{0mm}{3mm} & Model I & Model II\\
\hline
\rule[-1mm]{0mm}{4mm} $k$ & $1.719\cdot10^3$ & $15$ s$^{-1}$\\
\rule[-1mm]{0mm}{3mm} $\ell$ & $3.437\cdot10^3$ & $30$ s$^{-1}$\\
\rule[-1mm]{0mm}{3mm} $k_\ast$ & $15$ s$^{-\alpha}$ & $15$ s$^{-1}$\\
\rule[-1mm]{0mm}{3mm} $\ell_\ast$ & $30$ s$^{-\alpha}$ & $30$ s$^{-1}$\\
\rule[-1mm]{0mm}{3mm} $u_\infty$ (mol$\,$m$^{-1}$) & $8.698\cdot10^{-3}$ & $8.698\cdot10^{-3}$\\
\rule[-1mm]{0mm}{3mm} $v_\infty$ (mol$\,$m$^{-1}$) & $4.349\cdot10^{-3}$ & $4.349\cdot10^{-3}$\\
\rule[-1mm]{0mm}{3mm} $k_{eq}$ (s$^{-1}$) & $1.047\cdot10^6$ & $15$\\
\rule[-1mm]{0mm}{3mm} $\ell_{eq}$ (s$^{-1}$) & $2.094\cdot10^6$ & $30$
\end{tabular}
\end{center}
\caption{Monomolecular reversible reaction. Reaction rates and theoretical steady state.}
\label{table:parameters_monomolecular}
\end{table}
The numerical method is here applied to monomolecular reversible reactions. Consider the reactions
\begin{equation}
\begin{array}{ll}
\mbox{Model I}\quad & A_i \mathop{\rightleftharpoons}\limits^{k/\tau_i}_{\ell/\tau_i} B_i,\\
[10pt]
\mbox{Model II}\quad & A_i \mathop{\rightleftharpoons}\limits^{k}_{\ell} B_i,
\end{array}
\end{equation}
where the coefficients $k$ and $\ell$ used in the simulations are given in Table~\ref{table:parameters_monomolecular}. As discussed in Section~\ref{sec:IS_mono}, these cases correspond to the macroscopic FPDEs \eqref{eq:FPDE_monoI} and \eqref{eq:FPDE_monoII} with macroscopic reaction rates $k_\ast$ and $\ell_\ast$ given in Table~\ref{table:parameters_monomolecular}. The rates $k$ and $\ell$ have been chosen such that the two models have the same macroscopic reaction rates $k_\ast$ and $\ell_\ast$ and the same steady states. Figure \ref{fig:test3_u} then illustrates the differences between the two models.
The following initial conditions are used
\begin{equation}
\fatu(x,0) = \fatv(x,0) = \fatmu\,g(x).
\label{eq:initial_condition_test3}
\end{equation}
Figure~\ref{fig:test3_u} compares the numerical solutions of $U=\fate^T\,\fatu$ and $V=\fate^T\,\fatv$ obtained with the internal state reaction-diffusion system \eqref{eq:RDPDE_internalstates1}--\eqref{eq:RDPDE_internalstates2} with the analytical solutions of the FPDEs in Section~\ref{sec:analytical_mono}. The difference between the models is obvious in Figure~\ref{fig:test3_u}-(a) corresponding to model I at time $t_1 = 10^{-2}$ s and Figure~\ref{fig:test3_u}-(b) with the model II at the same time. For both models, good agreement is found between the mesoscopic PDE and the macroscopic FPDE solutions.

\begin{figure}[htbp]
\begin{center}
\begin{tabular}{cc}
Model I & Model II\\
\includegraphics[scale=0.30]{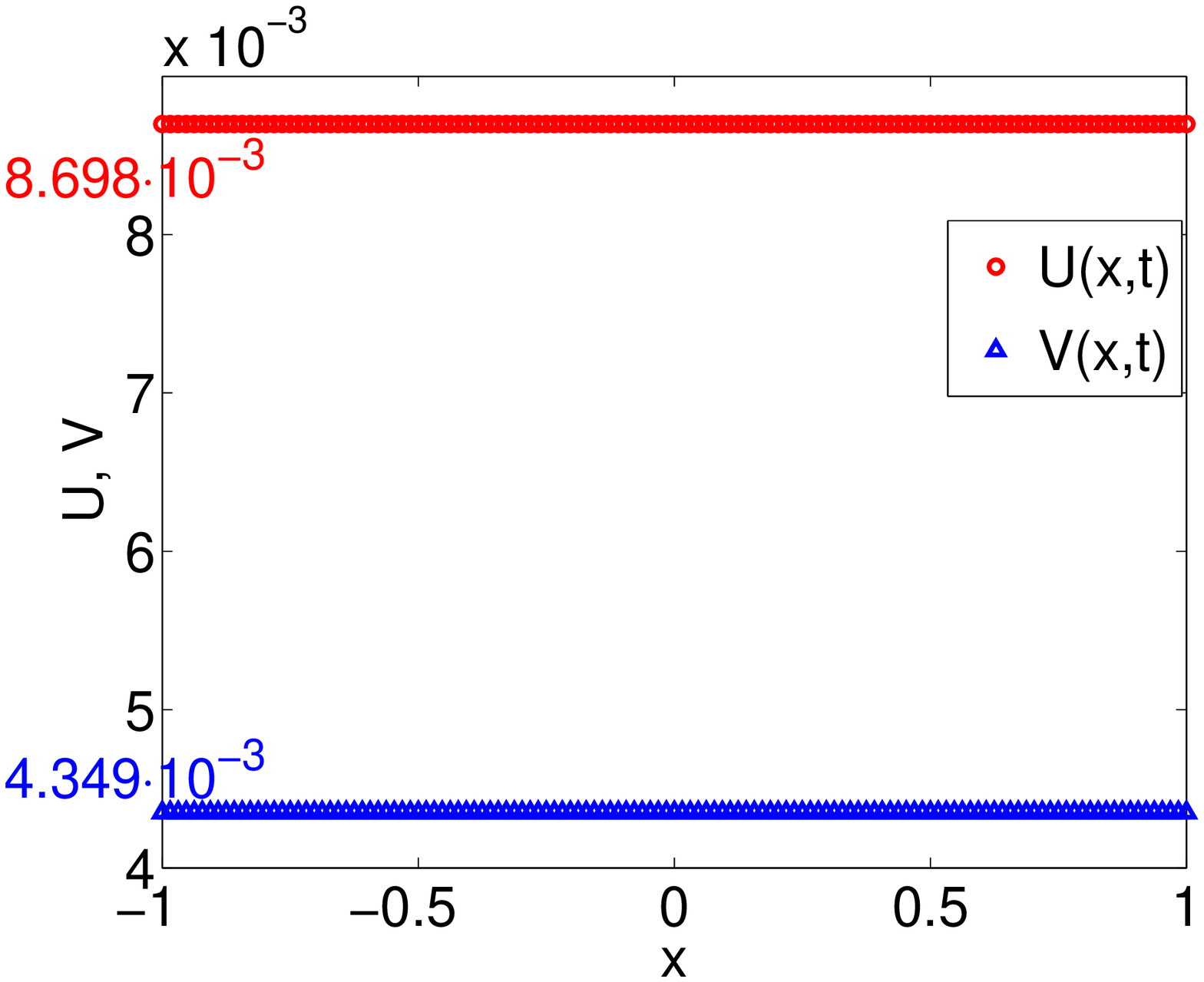} &
\includegraphics[scale=0.30]{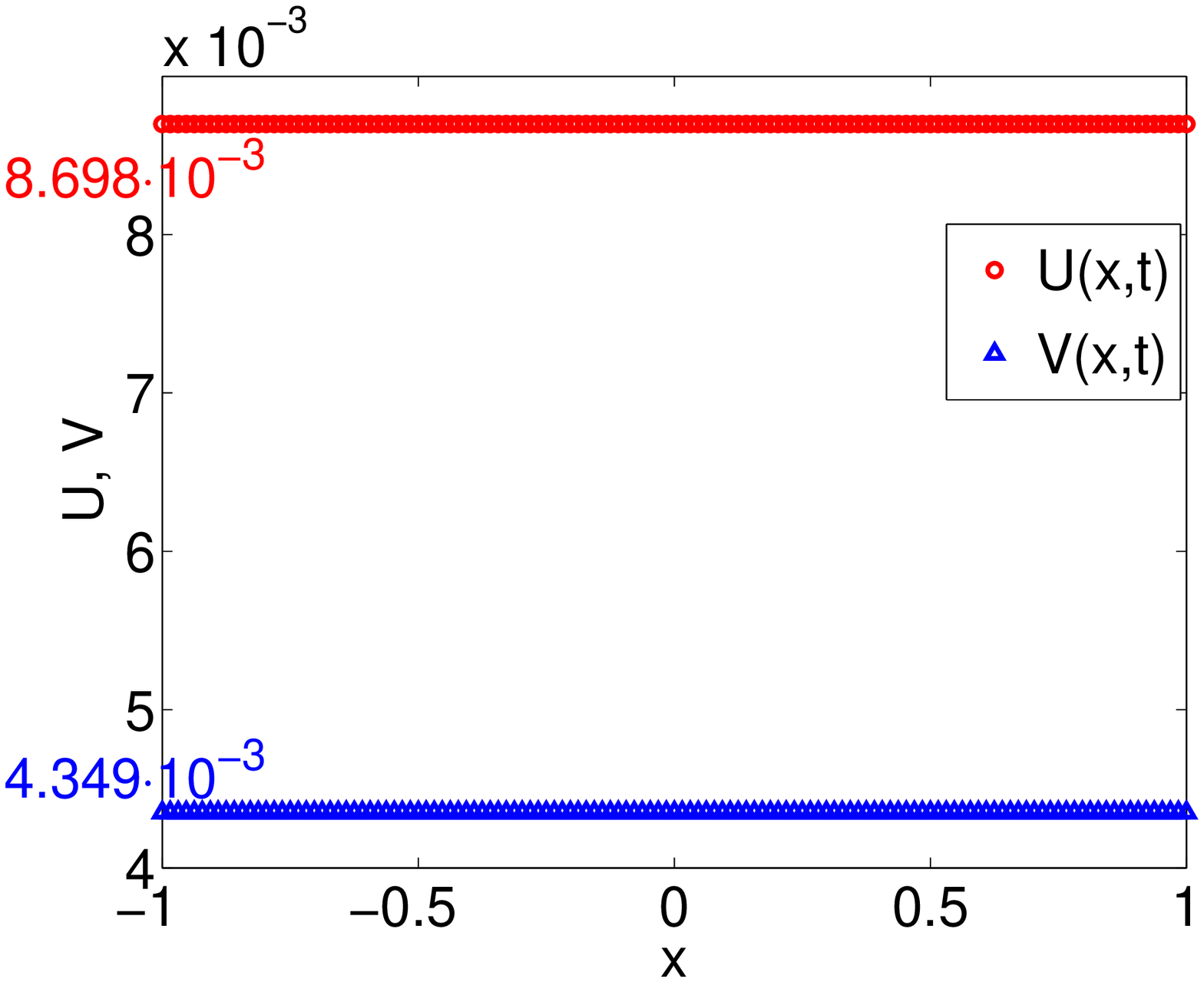}
\end{tabular}
\end{center}
\caption{Section~\ref{sec:linear_reaction}. Numerical values of the concentration $U=\fate^T\,\fatu$, $V=\fate^T\,\fatv$ of $A$ and $B$ at time $t_2 = 10$ s.}
\label{fig:test3_u_ss}
\end{figure}

The numerical solution of the internal states reaction-diffusion system at time $t_2 = 10$ s is depicted in Figure~\ref{fig:test3_u_ss} when the steady state is reached. We note that the numerical values of the steady states are close to the theoretical ones, given by the kernel of the matrix $\fatB$ in \eqref{eq:Bdef}. Using \eqref{eq:kldef}, the equivalent reactions rates $k_{eq}$ and $\ell_{eq}$ are computed, cf Table~\ref{table:parameters_monomolecular}, and the property $k_{eq}\,u_\infty = \ell_{eq}\,v_\infty$ when $t\rightarrow\infty$ in \eqref{eq:RDPDE_internalstates3_inf1} and \eqref{eq:RDPDE_internalstates3_inf2} is verified.



\subsection{Bimolecular reaction}\label{sec:nonlinear_reaction}

\begin{table}[htbp]
\begin{center}\footnotesize
\begin{tabular}{l|ccc}
\rule[-1.5mm]{0mm}{3mm} & Reaction I & Reaction II & Reaction III\\
\hline
\rule[-1mm]{0mm}{4mm} $k$ & $37.5$ m$\,$mol$^{-1}$ & $37.5$ m$\,$mol$^{-1}$ & $7500$ m$\,$mol$^{-1}\,$s$^{-1}$\\
\rule[-1mm]{0mm}{3mm} $\ell$ & $0.25$ & $0.25$ & $5$ s$^{-1}$\\
\rule[-1mm]{0mm}{3mm} $u_\infty$ (mol$\,$m$^{-1}$) & $8.204\cdot10^{-3}$ & $5.826\cdot10^{-3}$ & $3.772\cdot10^{-3}$\\
\rule[-1mm]{0mm}{3mm} $v_\infty$ (mol$\,$m$^{-1}$) & $6.573\cdot10^{-3}$ & $4.195\cdot10^{-3}$ & $2.141\cdot10^{-3}$\\
\rule[-1mm]{0mm}{3mm} $w_\infty$ (mol$\,$m$^{-1}$) & $1.582\cdot10^{-3}$ & $3.959\cdot10^{-3}$ & $6.014\cdot10^{-3}$\\
\rule[-1mm]{0mm}{3mm} $k_{eq}$ (m$\,$mol$^{-1}\,$s$^{-1}$) & $2.450\cdot10^3$ & $1.187\cdot10^5$ & $2.483\cdot10^1$\\
\rule[-1mm]{0mm}{3mm} $\ell_{eq}$ (s$^{-1}$) & $8.351\cdot10^1$ & $7.328\cdot10^2$ & $10^2$\\
\end{tabular}
\end{center}
\caption{Bimolecular reversible reaction. Reaction rates and theoretical steady states.}
\label{table:parameters_bimolecular}
\end{table}
The purpose of the last example is to establish whether the numerical methods presented in this paper can be used to handle more complex reactions. As an example, we consider the following bimolecular reversible reactions
\begin{equation}
\begin{array}{ll}
\mbox{Reaction I}\quad & A_i + B_i \mathop{\rightleftharpoons}\limits^{k/\tau_i}_{\ell/\tau_i} C_i,\\
[10pt]
\mbox{Reaction II}\quad & A_i + B_j \mathop{\rightleftharpoons}\limits^{\frac{k}{2}\,\left(\frac{1}{\tau_i} + \frac{1}{\tau_j}\right)}_{\ell/\tau_k} C_k,\\
[10pt]
\mbox{Reaction III}\quad & A_i + B_i \mathop{\rightleftharpoons}\limits^{k}_{\ell} C_i.
\end{array}
\end{equation}
The reaction rates are found in Table~\ref{tab:coefficients} with $\tau_{ij}$ for reaction II chosen as in \eqref{eq:CollinsKimball2} 
with $\theta=1/2$ since $A$ and $B$ diffuse with the same speed. The rates of reactions I and III are as in model I and II in Sections~\ref{sec:IS_annihilation} 
and \ref{sec:IS_mono} but reaction II is more general.
 The initial conditions are 
\begin{equation}
\fatu(x,0) = \frac{1}{2}\,\fatmu\,g(x),\quad
\fatv(x,0) =  \frac{1}{4}\,\fatmu\,g(x),\quad
\fatw(x,0) =  \fatmu\,g(x).
\label{eq:initial_condition_test4}
\end{equation}
The numerical solutions $U=\fate^T\,\fatu$, $V=\fate^T\,\fatv$ and $W=\fate^T\,\fatw$ of the internal states reaction-diffusion system are shown in Figure \ref{fig:test3_u} at time $t_1 = 10^{-2}$ s (top) and at time $t_2 = 10$ s (bottom) when the steady state is reached. For reaction I, we expect that the internal states model in \eqref{eq:bimo_eq2} approximates the macroscopic FPDE model I \eqref{eq:FPDE_bimoI}. No analytical solution of the FPDE is available in this case. 
For reaction II and reaction III, the macroscopic level with summation over the internal states is not so easily expressed as a FPDE. Using \eqref{eq:mat_bimolecular}--\eqref{eq:kldef3}, the equivalent reactions rates $k_{eq}$ and $\ell_{eq}$ are computed in Table~\ref{table:parameters_bimolecular} and the property ${k_{eq}\,u_\infty\,v_\infty = \ell_{eq}\,w_\infty}$ when $t\rightarrow\infty$ in \eqref{eq:bimolintstates1}--\eqref{eq:bimolintstates3} is verified.

\begin{figure}[htbp]
\begin{center}
\begin{tabular}{ccc}
Reaction I & Reaction II & reaction III\\
\includegraphics[scale=0.27]{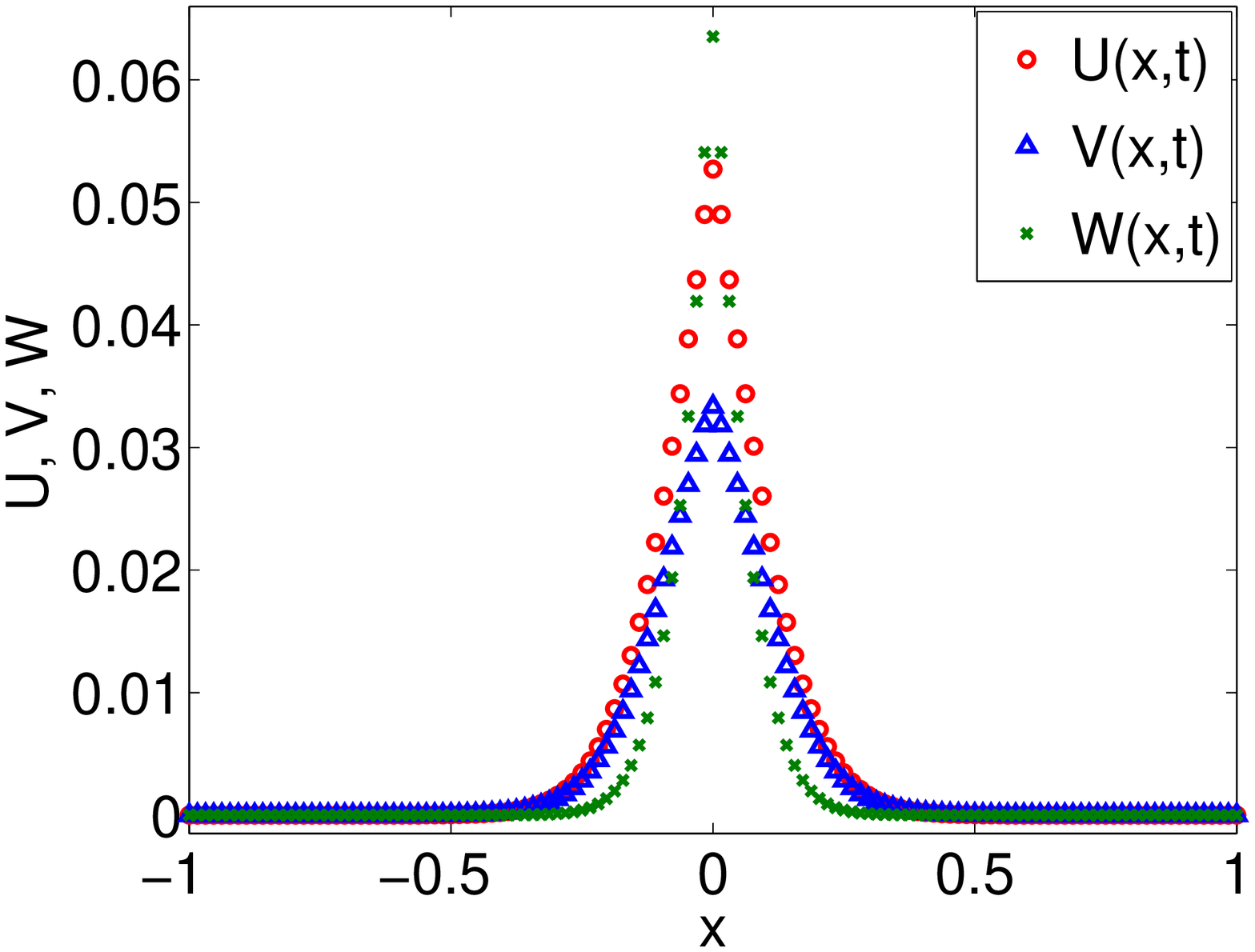} &
\includegraphics[scale=0.27]{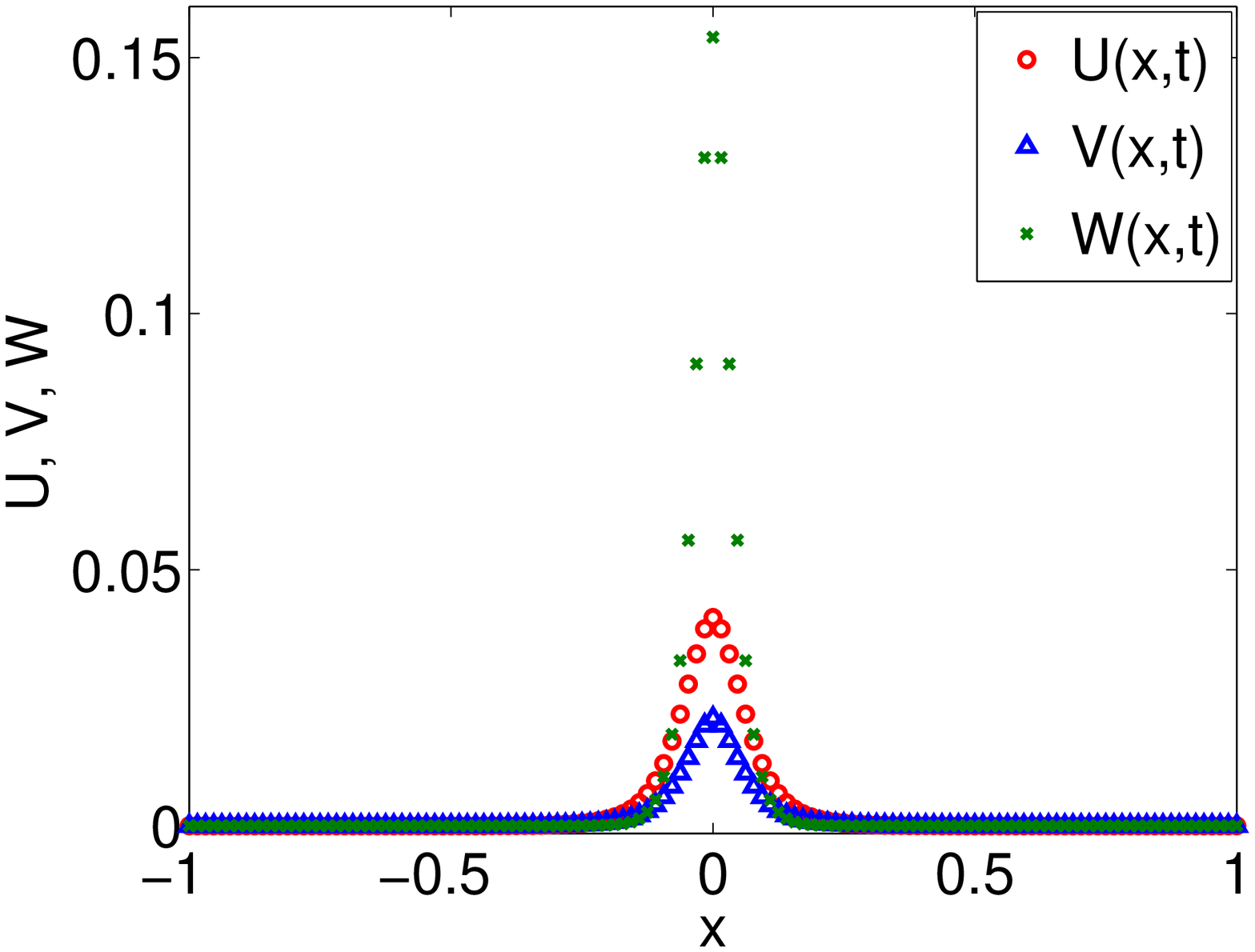} &
\includegraphics[scale=0.27]{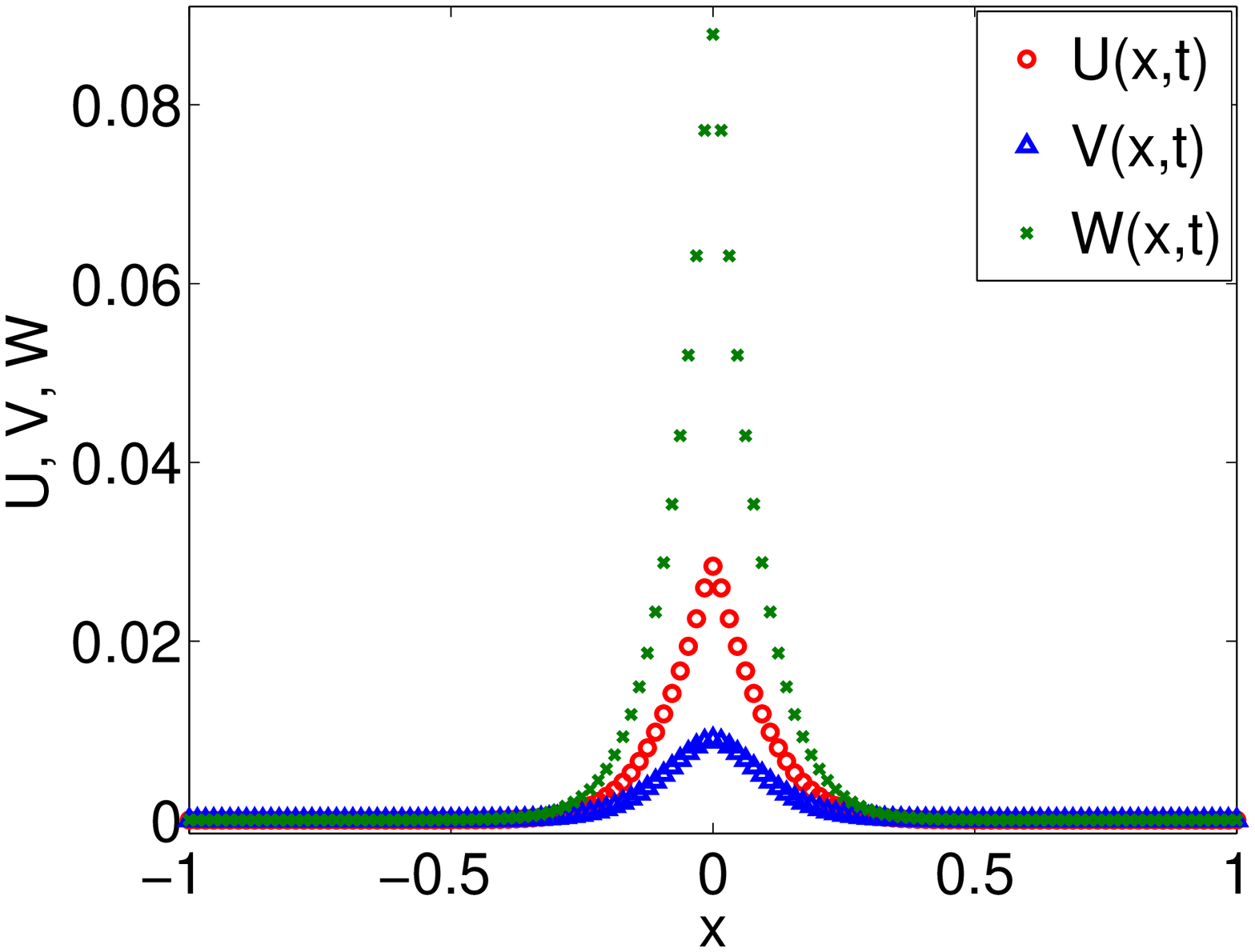}\\
\includegraphics[scale=0.27]{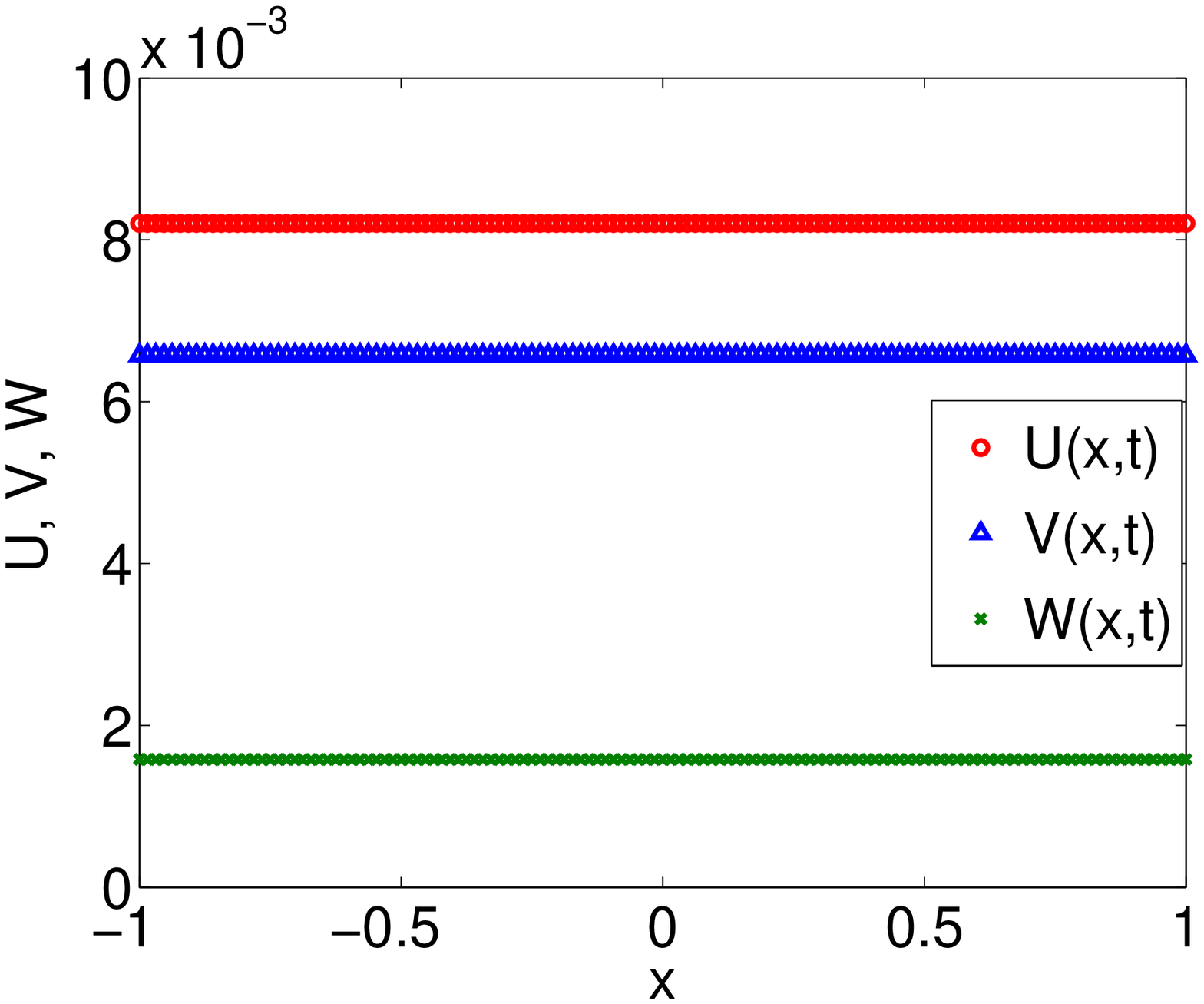} &
\includegraphics[scale=0.27]{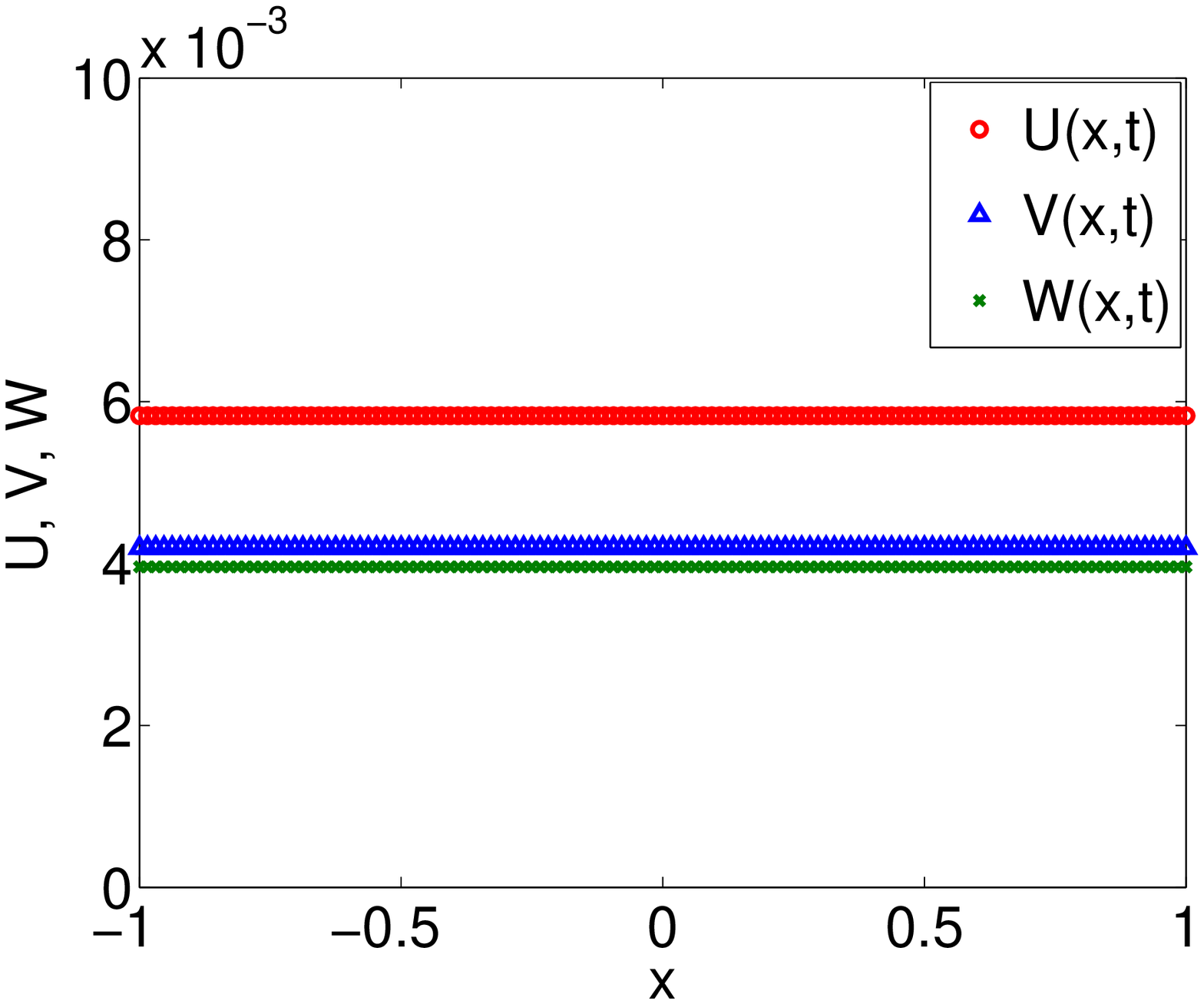} &
\includegraphics[scale=0.27]{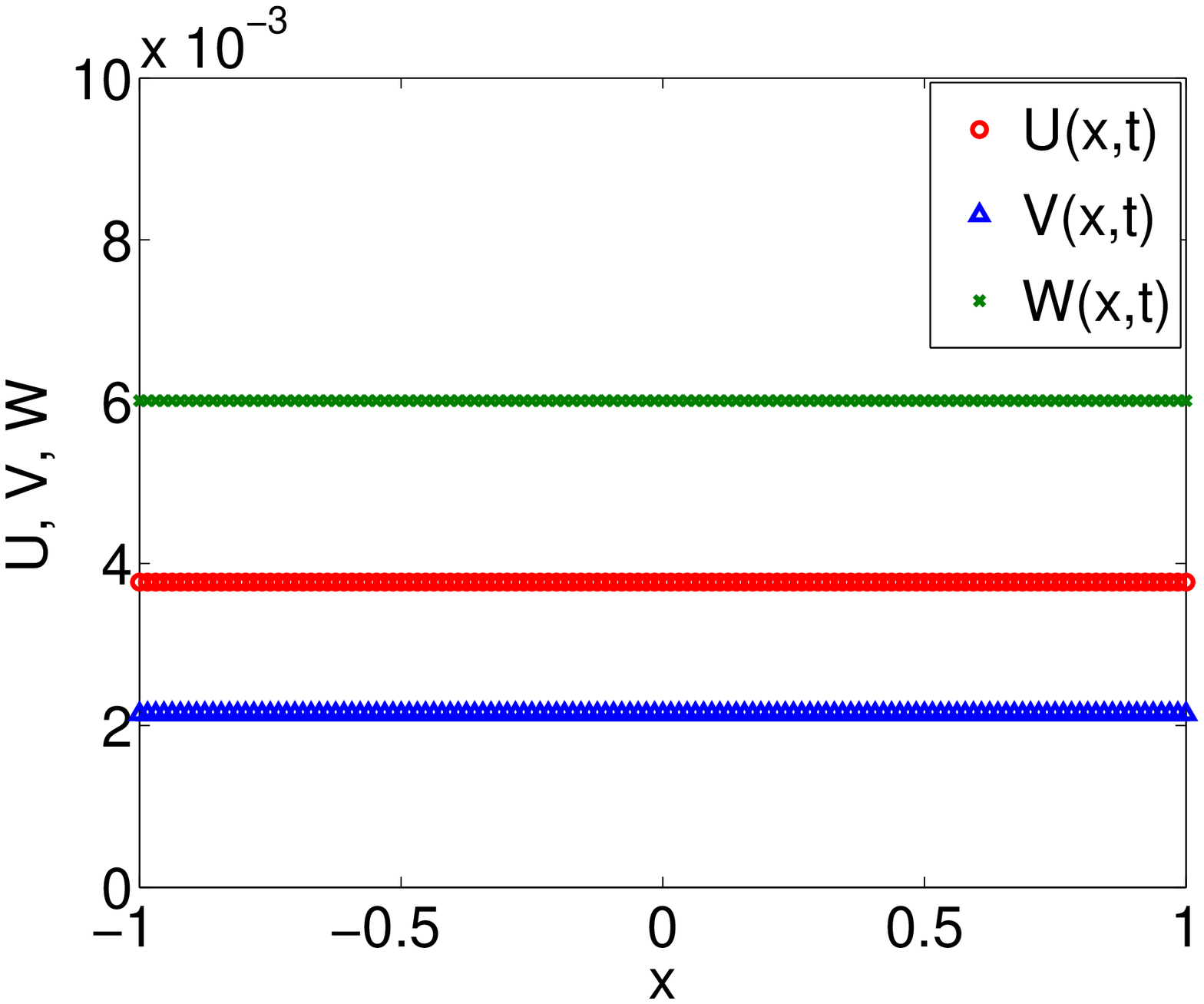}
\end{tabular}
\end{center}
\caption{Section~\ref{sec:nonlinear_reaction}. Numerical values of the concentration $U=\fate^T\,\fatu$, $V=\fate^T\,\fatv$, $W=\fate^T\,\fatw$ of $A$, $B$ and $C$ at time $t_1 = 1\cdot 10^{-2}$ s (top) and $t_2 = 10$ s (bottom).}
\label{fig:test4_u}
\end{figure}
%


\section{Conclusion}

A numerical method is presented here for simulating fractional-in-time
reaction-diffusion equations. A diffusive representation transforms
the function $\frac{1}{t^{\alpha+1}}$, involved in the mesoscopic CTRW
on a lattice, into a continuum of decreasing exponentials,
approximated by quadrature formulae. The CTRW model is then replaced
by an approximation, much more tractable numerically. Contrary to the
approach used in \cite{MOMMER09}, the coefficients of the diffusive
approximation are determined by a nonlinear optimization procedure,
leading to a smaller number of internal states. At the macroscopic
level, the internal states diffusion system thus obtained corresponds
to the fractional diffusion equation in a chosen time interval. In
contrast to the FPDE model, the diffusion in the internal states model
is ordinary at small and large times, but it is anomalous at
intermediate times. This behavior can also be observed in crowded
systems of hard-spheres due to caging effects, and hence the model
used herein may actually be a better model for microscopic crowding
than the traditional FPDE description.

The internal states model for diffusion in \cite{MOMMER09} is here
extended to account for chemical reactions. On the macroscopic FPDE
level, two different models for reactions with subdiffusion are
investigated. In model I the fractional derivative acts on both on the
standard diffusion term and the reaction term, whereas in model II the
fractional derivative acts only on the diffusion term. Both
macroscopic FPDE models correspond to a mesoscopic internal states
model with particular reaction coefficients. However, the opposite is
not true; mesoscopic models with general reactions may not have a
simple interpretation at the macroscopic level. In model I, the
reactions are subdiffusion controlled, that is the reaction kinetics
is ordinary at small and large times, whereas it is anomalous at
intermediate times. In model II, the reaction kinetics is ordinary for
all times. Which one of these models provides a better description of
a reaction system subject to subdiffusion does not have a simple
answer. In either case, the present work provids a theoretical
foundation for practical and efficient mesoscopic simulations.








\section*{Acknowledgment}
This work was supported by the Swedish strategic research programme eSSENCE, the UPMARC Linnaeus center of Excellence, and the NIH grant for StochSS with number 1R01EB014877-01.


\appendix
\section{Special functions} \label{annexe:special_functions}

The Fox-H function is defined as a Mellin-Barnes integral \cite{FOX61}
\begin{equation}
H^{m,n}_{p,q}\left[ z \left| \begin{array}{cccc}
(a_1,A_1) & (a_2,A_2) & \cdots & (a_p,A_p) \\
(b_1,B_1) & (b_2,B_2) & \cdots & (b_q,B_q)
\end{array} \right.  \right] = \frac{1}{2\,i\,\pi}\,\oint_L \frac{\prod\limits_{j=1}^m\Gamma(b_j - B_j\,s)\times\prod\limits_{j=1}^n\Gamma(1 - a_j + A_j\,s)}{\prod\limits_{j=m+1}^q\Gamma(1 - b_j + B_j\,s)\times\prod\limits_{j=n+1}^p\Gamma(a_j - A_j\,s)}\,z^s\,ds,
\label{eq:fox_function}
\end{equation}
where $L$ is a certain contour separating the poles of the two factors in the numerator. The special case for which the Fox-H function reduces to the Meijer-G function is $A_j = B_k = C$, $C > 0$ for $j=1,\cdots,p$ and $k=1,\cdots,q$ \cite{srivastava84}
\begin{equation}
\begin{array}{ll}
\displaystyle H^{m,n}_{p,q}\left[ z \left| \begin{array}{cccc}
(a_1,C) & (a_2,C) & \cdots & (a_p,C) \\
(b_1,C) & (b_2,C) & \cdots & (b_q,C)
\end{array} \right.  \right] &
\displaystyle = \frac{1}{2\,i\,\pi\,C}\,\oint_L \frac{\prod\limits_{j=1}^m\Gamma(b_j - s)\times\prod\limits_{j=1}^n\Gamma(1 - a_j + s)}{\prod\limits_{j=m+1}^q\Gamma(1 - b_j + s)\times\prod\limits_{j=n+1}^p\Gamma(a_j - s)}\,z^{s/C}\,ds,\\
[30pt]
& \displaystyle = \frac{1}{C}\,G^{m,n}_{p,q}\left[ z^{1/C} \left| \begin{array}{c}
a_1,\cdots,a_p \\
b_1,\cdots,b_p
\end{array} \right.  \right].
\end{array}
\label{eq:fox_to_meijer}
\end{equation}


\bibliographystyle{acm}
\bibliography{allbiblio}

\providecommand{\noopsort}[1]{} \providecommand{\doi}[1]{\texttt{doi:#1}}
  \providecommand{\available}[1]{Available at \texttt{#1}}
  \providecommand{\availablet}[2]{Available at \texttt{#2}}
\begin{thebibliography}{10}

\bibitem{BESB13}
{\sc Berkowitz, Y., Edery, Y., Scher, H., and Berkowitz, B.}
\newblock Fickian and non-{F}ickian diffusion with bimolecular reactions.
\newblock {\em Phys. Rev. E 87\/} (2013), 032812.

\bibitem{BLANC15}
{\sc Blanc, E.}
\newblock Approximation of the diffusive representation by decreasing
  exponential functions.
\newblock Tech. Rep. 2015-009, Department of Information Technology, Uppsala
  University, 2015.

\bibitem{blinov2004bionetgen}
{\sc Blinov, M.~L., Faeder, J.~R., Goldstein, B., and Hlavacek, W.~S.}
\newblock Bio{N}et{G}en: software for rule-based modeling of signal
  transduction based on the interactions of molecular domains.
\newblock {\em Bioinformatics 20}, 17 (2004), 3289--3291.

\bibitem{Bray03}
{\sc Bray, D.}
\newblock Molecular prodigality.
\newblock {\em Science 299}, 5610 (2003), 1189--1190.

\bibitem{CAPUTO67}
{\sc Caputo, M.}
\newblock Linear models of dissipation whose {Q} is almost frequency
  independent, part 2.
\newblock {\em Geophys. J. R. Astr. Soc. 13\/} (1967), 529--539.

\bibitem{CoKi}
{\sc Collins, F.~C., and Kimball, G.~E.}
\newblock Diffusion-controlled reaction rates.
\newblock {\em J. Colloid. Sci. 4\/} (1949), 425--437.

\bibitem{DESCH88}
{\sc Desch, W., and Miller, R.}
\newblock Exponential stabilization of {V}olterra integral equations with
  singular kernels.
\newblock {\em J. Int. Eq. Appl. 1}, 3 (1988), 397--433.

\bibitem{urdme}
{\sc Drawert, B., Engblom, S., and Hellander, A.}
\newblock {URDME}: a modular framework for stochastic simulation of
  reaction-transport processes in complex geometries.
\newblock {\em BMC Syst. Biol. 6}, 1 (2012), 76.

\bibitem{EndyBrent01}
{\sc Endy, D., and Brent, R.}
\newblock Modelling cellular behaviour.
\newblock {\em Nature 409\/} (2001), 391--395.

\bibitem{EnFeHeLo}
{\sc Engblom, S., Ferm, L., Hellander, A., and L\"{o}tstedt, P.}
\newblock Simulation of stochastic reaction-diffusion processes on unstructured
  meshes.
\newblock {\em SIAM J.~Sci.~Comput. 31\/} (2009), 1774--1797.

\bibitem{FaEl}
{\sc Fange, D., and Elf, J.}
\newblock Noise induced {M}in phenotypes in \textit{{E}. coli}.
\newblock {\em PLoS Comput. Biol. 2}, 6 (2006), e80.

\bibitem{FOX61}
{\sc Fox, C.}
\newblock The {G} and {H} functions as symmetrical {F}ourier kernels.
\newblock {\em Trans. Am. Math. Soc. 98}, 3 (1961), 395--429.

\bibitem{gillespie}
{\sc Gillespie, D.~T.}
\newblock A general method for numerically simulating the stochastic time
  evolution of coupled chemical reactions.
\newblock {\em J.~Comput.~Phys. 22}, 4 (1976), 403--434.

\bibitem{Gill77}
{\sc Gillespie, D.~T.}
\newblock Master equations for random walks with arbitrary pausing time
  distributions.
\newblock {\em Phys. Lett. 64A\/} (1977), 22--24.

\bibitem{GLOCKLE95}
{\sc Gl{\"o}ckle, W., and Nonnenmacher, T.}
\newblock A fractional calculus approach to self-similar protein dynamics.
\newblock {\em Biophys. J. 68}, 1 (1995), 46--53.

\bibitem{GLOCKLE91}
{\sc Gl{\"o}ckle, W.~G., and Nonnenmacher, T.~F.}
\newblock Fractional integral operators and {F}ox functions in the theory of
  viscoelasticity.
\newblock {\em Macromol. 24}, 24 (1991), 6426--6434.

\bibitem{HADDAR10}
{\sc Haddar, H., Li, J.~R., and Matignon, D.}
\newblock Efficient solution of a wave equation with fractional-order
  dissipative terms.
\newblock {\em J. Comput. Appl. Math. 234}, 6 (2010), 2003--2010.

\bibitem{mesoRD}
{\sc Hattne, J., Fange, D., and Elf, J.}
\newblock Stochastic reaction-diffusion simulation with {MesoRD}.
\newblock {\em Bioinformatics 21}, 12 (2005), 2923--2924.

\bibitem{HausKehr}
{\sc Haus, J.~W., and Kehr, K.~W.}
\newblock Diffusion in regular and disordered lattices.
\newblock {\em Phys. Rep. 150}, 5 (1987), 263--406.

\bibitem{HELESCHEWITZ00}
{\sc Heleschewitz, D.}
\newblock {\em Analyse et simulation de syst\`eme diff\'erentiels
  fractionnaires et pseudo-diff\'erentiels lin\'eaires sous repr\'esentation
  diffusive}.
\newblock PhD thesis, ENST, France, 2000.

\bibitem{HLW06}
{\sc Henry, B.~J., Langlands, T. A.~M., and Wearne, S.~L.}
\newblock Anomalous diffusion with linear reaction dynamics: From continuous
  time random walks to fractional reaction-diffusion equations.
\newblock {\em Phys. Rev. E 74\/} (2006), 031116.

\bibitem{steps}
{\sc Hepburn, I., Chen, W., Wils, S., and Schutter, E.}
\newblock {STEPS}: efficient simulation of stochastic reaction-diffusion models
  in realistic morphologies.
\newblock {\em BMC Syst Biol 6}, 35 (2012).

\bibitem{HofFra}
{\sc H{\"o}fling, F., and Franosch, T.}
\newblock Anomalous transport in the crowded world of biological cells.
\newblock {\em Rep. Progr. Phys. 76\/} (2013), 046602.

\bibitem{HORNUNG05}
{\sc Hornung, G., Berkowitz, B., and Barkai, N.}
\newblock Morphogen gradient formation in a complex environment: An anomalous
  diffusion model.
\newblock {\em Phys. Rev. E 72\/} (2005), 041916.

\bibitem{JEON11}
{\sc Jeon, J.-H., Tejedor, V., Burov, S., Barkai, E., Selhuber-Unkel, C.,
  Berg-S\o{}rensen, K., Oddershede, L., and Metzler, R.}
\newblock \textit{In Vivo} anomalous diffusion and weak ergodicity breaking of
  lipid granules.
\newblock {\em Phys. Rev. Lett. 106\/} (2011), 048103.

\bibitem{VanKampen}
{\sc {\noopsort{Kampen}}{~van Kampen, N.~G.}}
\newblock {\em Stochastic Processes in Physics and Chemistry}, 2nd~ed.
\newblock Elsevier, Amsterdam, 2004.

\bibitem{MONTROLL73}
{\sc Kenkre, V.~M., Montroll, E.~W., and Shlesinger, M.~F.}
\newblock Generalized master equations for continuous-time random walks.
\newblock {\em J. Stat. Phys. 9}, 1 (1973), 45--50.

\bibitem{KLAFTER12}
{\sc Klafter, J., Lim, S., and Metzler, R.}
\newblock {\em Fractional Dynamics: Recent Advances}.
\newblock World Scientific, 2012.

\bibitem{Kur70}
{\sc Kurtz, T.~G.}
\newblock Solutions of ordinary differential equations as limits of pure jump
  {M}arkov processes.
\newblock {\em J.~Appl.~Prob. 7\/} (1970), 49--58.

\bibitem{KUSUMI05}
{\sc Kusumi, A., Nakada, C., Ritchie, K., Murase, K., Suzuki, K., Murakoshi,
  H., Kasai, R.~S., Kondo, J., and Fujiwara, T.}
\newblock Paradigm shift of the plasma membrane concept from the
  two-dimensional continuum fluid to the partitioned fluid: {H}igh-speed
  single-molecule tracking of membrane molecules.
\newblock {\em Annu. Rev. Biophys. Biomol. Struct. 34\/} (2005), 351--378.

\bibitem{Lawson2013}
{\sc Lawson, M.~J., Drawert, B., Khammash, M., Petzold, L., and Yi, T.-M.}
\newblock Spatial stochastic dynamics enable robust cell polarization.
\newblock {\em PLoS Comput. Biol. 9}, 7 (2013), e1003139.

\bibitem{LOMHOLT07}
{\sc Lomholt, M.~A., Zaid, I.~M., and Metzler, R.}
\newblock Subdiffusion and weak ergodicity breaking in the presence of a
  reactive boundary.
\newblock {\em Phys. Rev. Lett. 98\/} (2007), 200603.

\bibitem{pysb}
{\sc Lopez, C.~F., Muhlich, J.~L., Bachman, J.~A., and Sorger, P.~K.}
\newblock Programming biological models in {P}ython using {PySB}.
\newblock {\em Molecular Systems Biology 9}, 1 (2013).

\bibitem{MAINARDI00}
{\sc Mainardi, F., Raberto, M., Gorenflo, R., and Scalas, E.}
\newblock Fractional calculus and continuous-time finance {II}: the
  waiting-time distribution.
\newblock {\em Physica A: Stat. Mech. Appl. 287}, 3–4 (2000), 468--481.

\bibitem{MARQUEZLAGO12}
{\sc Marquez-Lago, T.~T., Leier, A., and Burrage, K.}
\newblock Anomalous diffusion and multifractional {B}rownian motion: simulating
  molecular crowding and physical obstacles in systems biology.
\newblock {\em IET Syst. Biol. 6}, 4 (2012), 134--142.

\bibitem{METZLER00}
{\sc Metzler, R., and Klafter, J.}
\newblock The random walk's guide to anomalous diffusion: a fractional dynamics
  approach.
\newblock {\em Phys. Rep. 339}, 1 (2000), 1--77.

\bibitem{MILLER93}
{\sc Miller, K.~S., and Ross, B.}
\newblock {\em An Introduction to the Fractional Calculus and Fractional
  Differential Equations}.
\newblock Wiley, New-York, 1993.

\bibitem{MOMMER09}
{\sc Mommer, M.~S., and Lebiedz, D.}
\newblock Modeling subdiffusion using reaction diffusion systems.
\newblock {\em SIAM J. Appl. Math. 70}, 1 (2009), 112--132.

\bibitem{MONTROLL65}
{\sc Montroll, E.~W., and Weiss, G.~H.}
\newblock Random walks on lattices. {II}.
\newblock {\em J. Math. Phys. 6}, 2 (1965), 167--181.

\bibitem{Saxton02}
{\sc Saxton, M.~J.}
\newblock Chemically limited reactions on a percolation cluster.
\newblock {\em J. Chem. Phys. 116\/} (2002), 203--208.

\bibitem{Saxton07}
{\sc Saxton, M.~J.}
\newblock A biological interpretation of transient anomalous subdiffusion. {I}.
  {Q}ualitative model.
\newblock {\em Biophys. J. 92\/} (2007), 1178--1191.

\bibitem{SEKI03}
{\sc Seki, K., Wojcik, M., and Tachiya, M.}
\newblock Fractional reaction-diffusion equation.
\newblock {\em The Journal of chemical physics 119}, 4 (2003), 2165--2170.

\bibitem{Shkilev09}
{\sc Shkilev, V.~P.}
\newblock Effect of microscopic inhomogeneity of the medium on the
  reaction-diffusion front velocity.
\newblock {\em J. Exp. Theor. Phys. 108}, 2 (2009), 356--363.

\bibitem{Shkilev14}
{\sc Shkilev, V.~P.}
\newblock Comment on "{A}nomalous versus slowed-down {B}rownian diffusion in
  the ligand-binding equilibrium".
\newblock {\em Biophys. J. 106\/} (2014), 2541--2543.

\bibitem{SOKOLOV06}
{\sc Sokolov, I.~M., Schmidt, M. G.~W., and Sagu\'es, F.}
\newblock Reaction-subdiffusion equations.
\newblock {\em Phys. Rev. E 73\/} (2006), 031102.

\bibitem{SOULA14a}
{\sc Soula, H., Car{\'e}, B., Beslon, G., and Berry, H.}
\newblock {A}nomalous versus slowed-down {B}rownian diffusion in the
  ligand-binding equilibrium.
\newblock {\em Biophys. J. 105\/} (2013), 2064--2073.

\bibitem{srivastava84}
{\sc Srivastava, H., and Manocha, H.}
\newblock {\em A treatise on generating functions}.
\newblock Ellis Horwood {S}eries in {M}athematics and its {A}pplications. E.
  Horwood, 1984.

\bibitem{STAFFANS94}
{\sc Staffans, O.~J.}
\newblock Well-posedness and stabilizability of a viscoelastic equation in
  energy space.
\newblock {\em Trans. Amer. Math. Soc. 345}, 2 (1994), 527--575.

\bibitem{SBSK}
{\sc Stefan, M.~I., Bartol, T.~M., Sejnowski, T.~J., and Kennedy, M.~B.}
\newblock Multi-state modeling of biomolecules.
\newblock {\em PLoS Comput. Biol. 10}, 9 (2014), e1003844.

\bibitem{Sturrock2013}
{\sc Sturrock, M., Hellander, A., Aldakheel, S., Petzold, L., and Chaplain, M.}
\newblock The role of dimerisation and nuclear transport in the {Hes1} gene
  regulatory network.
\newblock {\em Bull. Math. Biol.\/} (2013), 1--33.

\bibitem{YADAV06}
{\sc Yadav, A., and Horsthemke, W.}
\newblock Kinetic equations for reaction-subdiffusion systems: Derivation and
  stability analysis.
\newblock {\em Phys. Rev. E 74\/} (2006), 066118.

\bibitem{YANG03}
{\sc Yang, H., Luo, G., Karnchanaphanurach, P., Louie, T.-M., Rech, I., Cova,
  S., Xun, L., and Xie, X.~S.}
\newblock Protein conformational dynamics probed by single-molecule electron
  transfer.
\newblock {\em Science 302}, 5643 (2003), 262--266.

\bibitem{YAL04}
{\sc Yuste, S.~B., Acedo, L., and Lindenberg, K.}
\newblock Reaction front in an {$A+B\rightarrow C$} reaction-subdiffusion
  process.
\newblock {\em Phys. Rev. E 69\/} (2004), 036126.

\bibitem{YLL07}
{\sc Yuste, S.~B., Lindenberg, K., and Ruiz-Lorenzo, J.~J.}
\newblock Subdiffusion limited reactions.
\newblock In {\em Anomalous {T}ransport: {F}oundations and {A}pplications\/}
  (Weinheim, 2007), R.~Klages, G.~Radons, and I.~M. Sokolov, Eds., Wiley-{VCH},
  pp.~3--33.

\end{thebibliography}

\end{document}